\documentclass[11pt, reqno]{amsart}

\usepackage{amsmath}
\usepackage{amssymb}

\usepackage{amsfonts}

\usepackage{float}
\usepackage{caption}
\usepackage{subfig}

\usepackage{mathtools}

\usepackage[all]{xy}
\usepackage{graphicx}
\usepackage{epstopdf,auto-pst-pdf,pst-node,epsfig}

\newtheorem{theorem}{Theorem}[section]

\newtheorem{defi}[theorem]{Definition}
\newtheorem{lemma}[theorem]{Lemma}

\newtheorem{rem}[theorem]{Remark}

\def\binom#1#2{{#1}\choose{#2}}

\def\slfrac#1#2{\hbox{\kern.1em %
 \raise.5ex\hbox{\the\scriptfont0 #1}\kern-.11em %
 /\kern-.15em\lower.25ex\hbox{\the\scriptfont0 #2}}}

\newcommand{\eqn}[1]{(\ref{#1})}
\newcommand{\hsp}{\hspace*{\parindent}}

\newcommand{\eeq}{\end{equation}}
\newcommand{\beql}[1]{\begin{equation}\label{#1}}
\newcommand{\bsq}{{\vrule height .9ex width .8ex depth -.1ex }}

\font\phvr=phvr at 11pt
\newcommand\he[1]{\mbox{\phvr  #1}}

\newcommand{\tZ}{\tilde{Z}}

\newcommand{\pt}{\partial}

\newcommand{\ep}{\epsilon}

\newcommand{\llog}{{\rm Log~}}
\newcommand{\ZZ}{{\mathbb Z}}
\newcommand{\RR}{{\mathbb R}}

\newcommand{\PP}{{\mathbb P}}
\newcommand{\QQ}{{\mathbb Q}}
\newcommand{\CC}{{\mathbb C}}

\newcommand{\bx}{{\bf x}}

\newcommand{\hD}{{\he{D}}}
\newcommand{\hG}{{\he{G}}}
\newcommand{\hGZ}{{\he{F}_2}}
\newcommand{\hH}{{\he{H}}}

\newcommand{\sQ}{{\rm Q}}

\newcommand{\sA}{{\mathcal A}}
\newcommand{\sB}{{\mathcal B}}

\newcommand{\sG}{{\mathcal G}}

\newcommand{\sM}{{\mathcal M}}
\newcommand{\sN}{{\mathcal N}}
\newcommand{\sO}{{\mathcal O}}

\newcommand{\sR}{{\mathcal R}}

\newcommand{\sU}{{\mathcal U}}
\newcommand{\sV}{{\mathcal V}}

\newcommand{\sW}{{\mathcal W}}
\newcommand{\sX}{{\mathcal X}}

\newcommand{\sDD}{{\mathcal D}}

\newcommand{\tsM}{{\tilde{\mathcal M}}}
\newcommand{\tsN}{{\tilde{\mathcal N}}}

\newcommand{\vep}{{\varepsilon}}

\newcommand{\sgn}{{\rm sgn}}

\usepackage[T1]{fontenc}
\usepackage{mathrsfs}

\usepackage{epsfig}

\usepackage{color}

\title{The Lerch Zeta Function\\  III. Polylogarithms and Special Values}
\subjclass[2000]{Primary: 11M35,  Secondary: 33B30}

\author{Jeffrey C. Lagarias}
\thanks{The research of the first author was supported by NSF grants 
DMS-1101373 and DMS-1401224  and  that of the second author by NSF grant
DMS-1101368.}
\address{Department of Mathematics, University of Michigan,
Ann Arbor, MI 48109-1043,USA}
\email{lagarias@umich.edu}

\author{Wen-Ching Winnie Li}
\address{Department of Mathematics, Pennsylvania State University,
University Park, PA 16802-8401, USA}
\email{ wli@math.psu.edu}

\date{November 16, 2015}
\begin{document}

\begin{abstract}
This paper studies algebraic and analytic structures
associated with the Lerch zeta function,
  complex variables
viewpoint taken in part II.
The  { Lerch transcendent }
$\Phi (s,z,c) := \sum_{n=0}^\infty \frac{z^n}{ (n+c)^{s}}$ is
obtained from the Lerch zeta function $\zeta(s, a, c)$
by the change of variable $z = e^{2 \pi i a}$.
We show that it analytically continues  to a maximal
domain of holomorphy in three complex variables $(s, z, c)$, as
a multivalued function defined over the base manifold
$\CC \times (\PP^{1} (\CC) \smallsetminus \{0, 1, \infty \} )\times
( \CC \smallsetminus \ZZ )$
and compute the monodromy
functions describing the multivaluedness.
For  positive integer values  $s=m$ and  $c=1$  this function
is closely  related to the classical $m$-th order polylogarithm $Li_{m}(z)$.
We study its behavior as a function of two variables
$(z,c)$ for ``special values''  where $s=m$ is an integer.
For  $m \ge 1$ we show that it is  a  one-parameter deformation
of $Li_{m}(z)$,
which satisfies a linear ODE, depending on $c \in \CC$, of order $m+1$  of Fuchsian
type on the Riemann sphere.  We determine the associated
$(m+1)$-dimensional monodromy representation, which is a non-algebraic
deformation of the monodromy of $Li_m(z)$. \end{abstract}
\maketitle
\tableofcontents

\setlength{\baselineskip}{1.2\baselineskip}

%
%
%

\section{Introduction}
\hsp

In this paper we study the
{\em Lerch transcendent}  $\Phi(s,z,c)$,
defined by
\beql{101}
\Phi(s, z, c) :=  \sum_{n=0}^\infty \frac{z^n}{ (n+c)^{s}},
\eeq
which is obtained from the
Lerch zeta function
\beql{102aa}
\zeta(s, a,c) = \sum_{n=0}^{\infty} e^{2 \pi i n a}(n+c)^{-s}
\eeq
under the change of variable $z := e^{2 \pi i a}$.
The Lerch transcendent $\Phi(s, z, c)$ is  called by some  authors
the ``Lerch zeta function''   (e.g.  Oberhettinger \cite{Ob56}), although
 $\zeta(s,a,c)$  is the function originally
studied by Lerch \cite{Le1887} in 1887.
One  obtains by double specialization at $z=1$ and $c=1$
 the Riemann zeta function
 \beql{101h}
 \zeta(s) = \Phi(s, 1, 1) = \sum_{n=0}^{\infty} \frac{1}{(n+1)^s},
 \eeq
and this expansion is valid in the half-plane  $Re(s)>1$.

  In his 1900 problem list Hilbert \cite{Hi00}  raised
  a question related to
   the Lerch transcendent.
 This question appears  just after  the 18-th problem,
 perhaps intended as a prologue to  several of the subsequent problems.
  Hilbert remarked  that functions that satisfy
algebraic partial differential equations form a class of ``significant functions'' ,
but that a number of  important functions seem not to belong to this class.
He wrote:

\begin{quotation}
  {\em The function of
the two variables $s$ and $x$ defined by the infinite series
$$ \zeta(s, x) = x + \frac{x^2}{2^s} +\frac{x^3}{3^s}+ \frac{x^4}{4^s} + ...$$
which stands in close relation with the function $\zeta(s)$, probably
satisfies no algebraic partial differential equation. In the
investigation of this question the functional equation
$$ x \frac {\partial \zeta(s, x)}{\partial x} = \zeta( s - 1, x) $$
will have to be used.}
\end{quotation}

\noindent
The function $\zeta(s, x)$ is sometimes called  {\em Jonqui\'{e}re's function}
because it was studied in 1889  by  de Jonqui\'{e}re \cite{Jo1889}.
\footnote{Under the substitution $x= e^{2 \pi i a}$ it has also been called the {\em periodic zeta function} 
(Apostol \cite[Sec. 12.7]{Ap76}).}
It  is obtained as
$\zeta(s, x) = x\Phi(s, x, 1)$, where $\Phi(s, x,1)$ is
from the specialization of
  the Lerch transcendent at value $c=1$.
  In 1920  Ostrowski \cite{Ost20} justified Hilbert's assertion by
  proving that $\zeta(s, x)$  satisfies no algebraic differential
  equation.
    Further work done on this question is discussed in Garunk\v{s}tis
 and Lauren\v{c}ikas \cite{GL99}.

 \smallskip
 The  Lerch transcendent $\Phi(s,z, c)$, which has an extra variable $c$,
  circumvents Hilbert's objection
  and belongs to Hilbert's class of ``significant functions'' .
  This comes about as follows.
  We introduce the two linear partial differential operators with polynomial coefficients
  \beql{108c}
  \hD_{\Phi}^{-} :=  z \frac{\partial}{\partial z} + c, \quad \text{and} \quad  \, \hD_{\Phi}^{+} := \frac{\partial}{\partial c}.
  \eeq
  One can show that the  Lerch transcendent satisfies two
  independent  ladder relations
   \beql{109}
 \hD_{\Phi}^{-} \, \Phi(s, z, c)=\, \Phi(s-1, z, c),
\eeq
and
\beql{110}
 \quad  \hD_{\Phi}^{+} \, \Phi(s, z, c) = -s \, \Phi(s+1, z, c),
\eeq
see Theorem \ref{th23} below.
By combining these operators, one finds
that the Lerch transcendent  satisfies  a linear partial
 differential equation  with
 polynomial coefficients,
\beql{107}
 \quad (\hD_{\Phi}^{-} \hD_{\Phi}^{+})\Phi(s, z, c) :=
\left( z \frac{\partial}{\partial z}\frac{\partial}{\partial c}  + c \frac{\partial}{\partial c} \right)
 \Phi(s, z, c) = - s \Phi (s, z, c),
\eeq
so that  it  is a ``significant function'' in Hilbert's sense.
We comment more on this linear PDE  below.

 The  Lerch transcendent  yields   classical polylogarithms
 under suitable specialization of its variables (up to an
 inessential factor).
Taking  $s=m$ a positive integer, and further taking
$c=1$, yields a  function
closely related to the $m$-th order (Euler) polylogarithm
$$
Li_m(z) := \sum_{n=1}^\infty \frac{z^n}{n^{m}},
$$
namely
\beql{101a}
Li_m(z) = z\Phi(m, z, 1).
\eeq
We therefore define the function
\beql{101b}
Li_s(z, c) : = z \Phi(s, z, c) = \sum_{n=0}^{\infty} \frac{z^{n+1}}{ (n+c)^{s}},
\eeq
and call it
the {\em extended polylogarithm}. This function  interpolates
all polylogarithms via the  parameter $s$;
 simultaneously it gives a deformation of the polylogarithm with
 deformation parameter $c$.
For nonpositive integers $s=-m \le 0$
it is known that the resulting function $Li_{-m}(z, c)$
is a  rational function of the two variables $(z,c)$; we term it the
 {\em negative polylogarithm} of  order $-m$.

The Lerch transcendent has a connection with mathematical physics.
The  first author will show in \cite{Lag15} that
the Lerch zeta function has a fundamental association with  the real Heisenberg group.
The relation is  visible in the commutation relations
\beql{107a}
 \hD_{\Phi}^{+} \hD_{\Phi}^{-} -  \hD_{\Phi}^{-} \hD_{\Phi}^{+} =I
\eeq
satisfied by the  operators \eqref{109} and \eqref{110}. 
We may reformulate the  linear PDE  that the Lerch transcendent satisfies,
using the modified operator
\begin{equation}\label{111a}
\Delta_{\Phi} := \frac{1}{2}\Big(\hD_{\Phi}^{+} \hD_{\Phi}^{-} + \hD_{\Phi}^{-}\hD_{\Phi}^{+}) =
      z \frac{\partial}{\partial z}\frac{\partial}{\partial c}  + c \frac{\partial}{\partial c} +\frac{1}{2}I.
\end{equation}
One may then rewrite the PDE \eqref{107} in the eigenfunction form
\beql{107b}
\Delta_{\Phi} \Phi(s, z, c) = -(s- \frac{1}{2}) \Phi(s, z, c).
\eeq
This  linear partial differential  operator 
  formally has the $xp$-form suggested as a possible
  form of a Hilbert-Polya operator encoding the zeta zeros as eigenvalues
   (Berry and Keating \cite{BK99}, \cite{BK99b});   more details are given  in Section \ref{sec9}.

The work of this  paper determines new basic analytic properties of this function,
which may bring insight to its specializations such as the 
Riemann zeta function and the polylogarithms. 
We construct an analytic continuation of the Lerch transcendent in all
three variables, revealing  its fundamental character as a  multivalued function.
We give an  exact determination of
its multivaluedness, specified by monodromy functions,
 and  determine  the effect of this multivalued analytic continuation on the
partial differential equations and difference equations above.
An important feature is that this  analytic continuation does not extend to
certain sets
over a base manifold which we term {\em singular strata}; these form
 the branch locus for the  multivaluedness. In  the $(z, c)$-variables these
are points where either $c$ is a nonpositive integer and/or  where $z=1$ or $z=0$.
In particular this three-variable analytic continuation  omits  the specialization to the Riemann zeta function given in \eqref{101h},
which  occurs at  the singular stratum point $(z, c)=(1,1)$.
The Lerch transcendent does possess  additional analytic
continuations in fewer variables 
valid on some singular strata outside the analytic
continuation in three variables; for example the Hurwitz zeta
function $\zeta(s, c)=\sum_{n=0}^{\infty} (n+c)^{-s}$ arises on (one branch of) the singular stratum $z=1$ with $c$ variable.
These additional analytic continuations typically
include meromorphic continuations in the $s$-variable to all $s \in \CC$.
However  for many parameter ranges these functions in fewer variables are not continuous limits of  the three-variable analytic continuations.

A particular goal of the paper is to understand the
 relation of this analytic continuation to the
multivalued structure associated with  the polylogarithm.
The specialization to obtain the polylogarithm  takes $c=1$
(and also $s=m \ge 1$, a positive integer),  which 
lies on a singular  stratum over the Lerch zeta manifold $\sM$ in the $(a, c)$-variables described in
Section \ref{sec11}. In part II however we showed that the three-variable
analytic continuation has removable singularities at
the points  $c=n$ for $n \ge 1$, and hence  extends to a larger manifold $\sM^{\sharp}$.
The polylogarithm case $c=1$  is therefore covered  in this extended analytic continuation
when projected to the manifold $\sN^{\sharp}$ in the $(z, c)$-variables described below.

 The Lerch transcendent has the new feature that its analytic continuation
 introduces   a new  singular stratum consisting of the   $z=0$ manifold,
which is a branch locus around which  it is multivalued, whose monodromy must be determined. This singular stratum is not directly
visible in the $(a, c)$ variables used in the Lerch zeta function.
Upon specializing two variables to obtain  the polylogarithm in the $z$-variable,
we obtain  a new determination  of  its  monodromy structure,
and we also obtain  an interesting one-parameter deformation
of the  polylogarithm in the $c$-variable.

An additional  reason  for interest in this $c$-deformation,
apart from containing the polylogarithms,  concerns
the behavior of  functional equations satisfied by the Lerch transcendent.
A major property of the
Lerch zeta function $\zeta(s, a, c)$  is that it satisfies
three-term and four-term functional equations relating
certain  linear combinations of functions at parameter values $s$ to those at $1-s$.
The functional equation of the Riemann zeta function  $\zeta(s)$ and $\zeta(1-s)$
can be derived for these functional equations, proceeding by
a limiting process to a singular stratum yielding functional equations for the  Hurwitz zeta function (when $Re(s)>0$)
and for the  periodic zeta function (when $Re(s)<1)$,
and from these recovering that of the Riemann zeta function, see Apostol  \cite[Chap. 12]{Ap76}.
The Lerch transcendent inherits multivalued versions of these functional equations,
well-defined for all non-integer values of the $c$-parameter,
but they { fail to extend} to the parameter values corresponding to polylogarithms, as we explain in
Section \ref{sec13}.

The following Sections \ref{sec11}-\ref{sec15}   discuss the  results of this paper in
general terms;  and the  main results  are stated in detail  in Section 2.

\subsection{Analytic continuation in three complex variables}\label{sec11}

  We establish an analytic continuation of the Lerch transcendent
  in three complex variables $(s, z, c)$ as a multivalued function of the variables $(z,c)$,
  which are entire functions of $s$.
  The special choice $s=n$, for $n$ a positive integer, gives  a one-parameter deformation of the
 $n$-th order polylogarithm.

 The existence of the analytic continuation follows using  the results of Part II,
  which gave a multivalued analytic continuation of the Lerch zeta function $\zeta(s, a, c)$,
  to a covering manifold of the manifold
  $$
  \sM= \{ (s, a, c)  \in \CC \times (\CC \smallsetminus \ZZ) \times
(\CC \smallsetminus \ZZ)\}. \,
  $$
  That paper also gave an extended analytic continuation to a covering manifold of
  $$
  \sM^{\sharp}= \{ (s, a, c)  \in \CC \times (\CC \smallsetminus \ZZ) \times
(\CC \smallsetminus \ZZ_{\le 0})\}.
  $$
 In this paper  we lift  this continuation to the Lerch transcendent $\Phi(s, z, c)$  using the multivalued inverse
  change of variable
  $$a= \frac{1}{2 \pi i } (\log z),$$
  to the map $z= e^{2 \pi i a}.$ We choose a branch\footnote{Please note the conventions on logarithms given at
  the end of Section \ref{sec2}.}
 of the logarithm that agrees with the usual definition in the upper-half plane, and on its boundary
  takes $\log(1)=0, \log(-1)= \pi i$.  Our choice of base point on $\sM^{\sharp}$ is $(s_0, a_0, c_0) = (\frac{1}{2}, \frac{1}{2}, \frac{1}{2})$
  and with this choice of branch the point $z_0=-1$ lifts to $a_0=\frac{1}{2}$.
 The base manifold
  now becomes
 $$
\sN :=\{ (s, z,c) \in \CC \times (\PP^1(\CC)\smallsetminus \{0,1, \infty\}) \times
(\CC \smallsetminus \ZZ)\} \,
$$
and, in the extended analytic continuation case,
$$
\sN^{\sharp}  :=\{ (s, z,c) \in \CC \times (\PP^1(\CC)\smallsetminus \{0,1, \infty\}) \times
(\CC \smallsetminus \ZZ_{\le 0} )\} .
$$
The formulas for the monodromy functions describing
the multivaluedness become correspondingly
  more complicated than those in Part II, and give  a representation $\rho_s$
 for fixed $s \in \CC$ of the fundamental group
  $$
\pi_1(\sN, \bx'_0) \simeq  \pi_1(\PP^1(\CC)\smallsetminus \{0,1, \infty\}, z_0) \times
\pi_1 (\CC \smallsetminus \ZZ, c_0)
$$
with specified base point $\bx'_0= ( s_0, z_0, c_0)= (\frac{1}{2},  -1, \frac{1}{2})$, acting on a direct sum
vector space spanned by the monodromy functions\footnote{The base point can be moved to $\bx_{s}^{'}= (s, -1, \frac{1}{2})$
since the manifold $\sN$ has a product structure
splitting off the $s$-coordinate, in which it is simply-connected.}. This vector space is generally infinite
dimensional, but it degenerates at  values where $s$ is an integer, see
Theorem \ref{th24}.

Two main features of the analytic continuation are:
\begin{enumerate}
\item[(1)] The Lerch transcendent becomes single-valued on a certain  covering manifold of $\sN$,
which is a regular  covering (i.e. Galois covering) with solvable covering group.
\item[(2)] The Lerch transcendent continues to satisfy
the two independent differential-difference equations \eqn{109} and \eqn{110}  on the covering manifold.
\end{enumerate}

We distinguish between  the {\em regular stratum}, which are all those parameter values where  the
analytic continuation exists, and {\em singular strata} which correspond to 
parameter values at which the analytic continuation gives a branch point.
That is,  the  analytic continuation does not apply to  certain
(complex) codimension one 
 {\em singular strata} of
$(z,c)$ parameter values, which include all values where the Lerch zeta function formally
becomes the Hurwitz zeta function or Riemann zeta function. The singular strata values are all $(z, c)$ such that
either $z=0$ or $z=1$ or $c$ is a nonpositive integer, or both. 
On certain  singular strata, continuous limiting values may exist for restricted
ranges of the $s$-parameter, as illustrated in some results in  Part I, e.g.  \cite[Theorem 2.3]{LL1}.

 In parallel to results given in Part II, the analytic continuation above has
removable singularities at positive integer values of $c$, and extends to an analytic  continuation over
the larger manifold $\sN^{\sharp}$ above.

\subsection{Specializations and Fuchsian ODE's}\label{sec12}

We study  consequences of this analytic continuation
for functions of fewer variables obtained by specializing the variables.
These specializations  include the $n$-th order polylogarithms, corresponding to a
``nonsingular'' specialization at $c=1$, and $s=n \ge 1$ a positive integer.
The classical specializations, giving rise to
the Hurwitz zeta function or Riemann zeta function, approach  singular strata
where the analytic continuation breaks down. Here limiting values do exist
for some ranges of the singular strata parameters, and a better understanding
of the nature of these degenerations seems  of
particular interest.

First,  we consider specialization to the  point $c=1$.
This value is  a ``nonsingular'' value for the extended analytic continuation to $\sN^{\sharp}$.
We  deduce the complete multivalued analytic continuation of the extended polylogarithm.
This covers  the case of  Hilbert's example function in the variables $(s, x)$ above.
As already noted, this  specialization loses the algebraic PDE property.

  Second, we consider the specialization of variables that treats $s$ as a constant.
  This specialization retains the linear PDE
  property in the $(z, c)$-variables, but loses the differential-difference equation property
  that depends on variation in $s$. In the case of $s=-m$ a non-positive integer,
   the functions are rational functions of two variables $(z, c)$, which are polynomial in
   the variable $c$, and which remain well-defined on certain singular strata in $c$ and $z$.
   For positive integers, i.e. $s=m \in \ZZ_{\ge 1}$,
 this  specialization gives a one-parameter deformation, with
 deformation parameter $c$,  of the classical
 $m$-th order polylogarithm $Li_{m}( z)$, which corresponds to taking $c=1$.

Third, specializing to  integer values  $s=m \ge 1$
and additionally specializing  $c \in \CC$ to be fixed, the
 specialized function $Li_{m}(z, c)$  satisfies
 a  linear  ordinary differential equation  of order $m+1$
in the $z$-variable. This ODE is  of Fuchsian type with regular singular
points at $\{ 0, 1, \infty\}$ for all values of the $c$-parameter;
 in particular  this differential equation is defined for    singular stratum  parameter values
$c \in \ZZ_{\le 0}$. We  determine the monodromy representation of the
fundamental group $\pi_1( \PP^1\smallsetminus \{0, 1, \infty\}, -1)$
for this equation as a function of the
deformation parameter $c$. The monodromy is unipotent for $c \in \ZZ$, is
quasi-unipotent for $c \in \QQ$, and otherwise  lies in a  Borel subgroup of
$GL(m+1, \CC)$ but is not quasi-unipotent.
A second interesting feature is that this  deformation of the monodromy varies
continuously  on the  regular stratum , but
has  discontinuous behavior  of the monodromy
representation at the
singular strata values $c \in \ZZ_{\le 0}$.

It is known that in the case $c=1$
 a mixed Hodge structure can be attached to the collection of
polylogarithms for all $n \ge 1$,  viewed as pro-unipotent connection over
$\PP^{1} \smallsetminus \{ 0, 1, \infty\}$ as  described in
 Bloch~\cite{Bl91}.
  We have not addressed the question whether a mixed Hodge structure can be associated
 to the singular strata cases $c \in \ZZ_{\le 0}$, where the
 monodromy is  unipotent.

\subsection{Functional Equations}\label{sec13}

It is  well known that the Lerch zeta function $\zeta(s, a, c)$ satisfies  a functional equation
relating parameter values  $s$ and $1-s$, which can be given in an asymmetric three-term form
given by  Lerch \cite{Le1887} as
$$
 \zeta (1-s,a,c) =  (2 \pi )^{-s} \Gamma (s)
 \left\{ e^{\frac{\pi is}{2}} e^{- 2 \pi iac} \zeta(s,1-c,a) + e^{-
\frac{\pi is}{2}} e^{ 2 \pi ic (1-a)} \zeta(s,c,1-a) \right\},
$$
and called {\em Lerch's transformation formula}, cf. \cite[Theorem 5.1]{LL1}.
It also satisfies two  symmetric four-term functional equations noted by Weil \cite[p. 57]{We76},
 discussed in Part I and Part II.  To state these,  define the two functions $L^{\pm}(s, a, c)$ by
$$
L^{\pm}(s, a, c) := \zeta(s, a, c) \pm e^{-2\pi i a} \zeta(s, 1-a, 1-c) = \sum_{n \in \ZZ} (\sgn (n))^k \frac{e^{2 \pi i n a}}{(n+c)^s}
$$
in which  $\pm := (-1)^k$.
Secondly, define their completions (adding archimedean Euler factors) by
$$
\hat{L}^{\pm}(s, a, c) ;= \pi^{\frac{s+ k}{2}} \Gamma( \frac{s+k}{2}) L^{\pm}(s, a, c),
$$
with $k \in \{0, 1\}$ and superscripts $\pm = (-1)^{k}$. Here \cite[Theorem 2.1]{LL2}  states that the two symmetrized four term functional equations
\beql{131a}
\hat{L}^{+}(s, a, c) = e^{- 2 \pi i a c} \hat{L}^{+} (1-s, 1-c, a)
\eeq
and
\beql{131b}
\hat{L}^{-}(s, a, c) = i e^{-2 \pi i   a c} \hat{L}^{-}(1-s, 1-c, a)
\eeq
hold for all $(s, a, c)$ in a simply-connected domain, the {\em fundamental polycylinder}
$$
\Omega := \{ 0 < Re(s) <1\} \times \{ 0 < Re(a) <1\} \times \{ 0 < Re(c)<1\}.
$$
Written out, the first functional equation \eqref{131a} has four terms,
and  states,
\begin{eqnarray*}
\pi^{-\frac{s}{2}} \Gamma(\frac{s}{2})\,
\Big(  \zeta(s, a, c) \,\,+ \,\,e^{-2\pi i a} \zeta(s, 1-a, 1-c)\Big)
&& \\
 = \,\, e^{-2 \pi i a c}\,\pi^{\frac{s-1}{2}}\Gamma(\frac{1-s}{2})
\Big(\zeta(1-s, 1-c, a) &+& e^{2\pi i c} \zeta(1-s, c, 1-a)\Big).
\end{eqnarray*}
These functional equations are {\em non-local} in the $(a, c)$-variables,
and it is important that they
leave the fundamental polycylinder $\Omega$ invariant.
One can combine two four-term  functional equations
in such a way as to recover the three-term functional equation for the Lerch zeta
function valid on the fundamental polycylinder.

The four term functional equations \eqref{131a}, \eqref{131b} project to  functional equations for the Lerch transcendent,
but become complicated to state since they  involve exponential and logarithmic
changes of variable.
If we  let $z= e^{2 \pi i a}$ then we find the projected image of the fundamental polycylinder
$\Omega$ to the manifold $\sN$ is
the simply-connected domain
$$
\sDD_{\sN} := \{ 0 < Re(s) <1\} \times \{ z \in \CC \smallsetminus \RR_{\ge 0} \} \times \{ 0 < Re(c)<1\}.
$$
If we set
$$
L_{\Phi}^{\pm}(s, z, c) := \Phi(s, z, c) \pm \frac{1}{z} \Phi(s, \frac{1}{z}, 1-c),
$$
with all variables in  $\sDD_{\sN}$ and complete these two functions  with appropriate  archimedean Euler factors  as above
then we will obtain four-term   functional equations with all four terms lying in $\sDD_{\sN}$.
We then obtain under analytic continuation four-term functional equations
in multivalued functions of shape
\beql{133a}
\hat{L}_{\Phi}^{+}(s, z, c) = z^{-c} \hat{L}_{\Phi}^{+}(1-s, e^{-2\pi i c}, \frac{1}{2 \pi i} \log z),
\eeq
and
\beql{133b}
\hat{L}_{\Phi}^{-}(s, z, c) = i z^{-c} \hat{L}_{\Phi}^{-}(1-s, e^{-2\pi i c}, \frac{1}{2 \pi i} \log z),
\eeq
provided that  correct choices are made of branches of all the multivalued functions
on each  side of the equation.
In  principle the results of this paper permit simultaneous determination of
the  multivaluedness of the four terms  in the  functional equation
following paths in  $\sN$
 starting from the fundamental polycylinder $\Omega$,
but this paper does not carry out
such a determination.

We conclude the topic of functional equations by pointing out  two  important issues, 
which remain to be resolved.

(1)
The two four-term  functional equations for the Lerch transcendent
 are well-defined  on the manifold $\sN$ but are not well-defined
 on the  extended manifold $\sN^{\sharp}$.
 The extended manifold $\sN^{\sharp}$ glues in the  integer values $c= n \ge 1$
in  the  Lerch transcendent $(s, z, c)$ parameters, and these extra
 values include  exactly the value $c=1$ relevant to studying polylogarithms.
This obstruction to  extension occurs because always at least  one of
the four terms in the functional equation  lies  on a genuine singular stratum.
 A consequence is that  four-term functional equations  do not appear
 when studying  the polylogarithm  itself.
Problems also occur with extending the  three-term functional equations.

It is possible that further  information  can be extracted
 from these functional equations   at the polylogarithm values,
if one approaches these points  along specific paths for restricted ranges of parameter values.
  In the $c$-deformation of  the polylogarithm we study,
 the  multivalued  functional equations  relating values at $s$ and $1-s$  ``turn on''  when $c$ takes a non-integer value.
 Perhaps some modified functional equations in fewer variables survive in the limit as
 a value $c = n \ge 1$ is approached for suitable ranges of the $s$-variable, because continuous limits
 to singular strata exist for  some range of $s$, as shown in part I.
One may also ask whether there is a ``vanishing cycle'' interpretation for
 some of this limit  behavior.

(2) On the other hand, at the polylogarithmic points $c=m \ge 1$, new functional equations
appear. Polylogarithms $Li_m(z)$ are well known to satisfy functional equations of quite different
shape, specific to each value of $m$, sometimes relating different values of $m$ shifted by integers.
These functional equations are relevant to geometry and physics,
and  relate these functions at different values  of $z$.
For the Euler dilogarithm  there is a well known functional equation found by Spence \cite{Sp1809}  in 1809,
often given in the form
\begin{equation}\label{SFE}
Li_2(\frac{x}{1-x} \frac{y}{1-y}) = Li_2(\frac{x}{1-y}) + Li_2( \frac{y}{1-x}) - Li_2(x) - Li_2(y) - \log(1-x)\log(1-y),
\end{equation}
see Lewin \cite[Sec. 1.2]{Lew91}. Since $Li_1(x) = -\log (1-x)$ this functional equation relates
polylogarithms with the two different $s$-parameter values  $s=2$ and $s=1$. These shifts
in the $s$-parameter  are different from
the $s$ to $1-s$ parameter shift in the four term functional equation.
The functional equation \eqref{SFE} can  be transformed into the well known $5$-term functional equation
for the  {\em Rogers dilogarithm}, given by $L_2(x) = Li_2(x) + \frac{1}{2} \log z \log(1-z)$, which is
$$
L_2(x) + L_2(y) - L_2(xy) = L_2(\frac{x- xy}{1- xy}) + L_2( \frac{y- xy}{1-xy}),
$$
cf. Rogers \cite{Rog07}, Zagier \cite{Zag88}, \cite{Zag07}. These particular functional equations
have  an important  relation  to three dimensional
geometry, specifically to   Cheeger-Cherns-Simons invariants of
hyperbolic $3$-manifolds, cf. Dupont \cite{Dup87}, Neumann \cite{Neu04}.
The functional equations and related ones for higher polylogarithms
seem specific to  integer values $s=m >1$  of the $s$-parameter,
and are not known to  survive
deformation in  $s$.  Nonetheless  one may ask what is the fate of the functional equations of the dilogarithm
 under the $c$-deformation presented in this paper, which in \eqref{SFE} involve
 the integer values $s=2$ and $s=1$.

  We remark that there are geometric generalizations of the dilogarithm
to higher dimensional cases, which aim to preserve  functional equations having
a geometric meaning, see Gel'fand and MacPherson \cite{GM82}, Hain and MacPherson \cite{HM90}.
These higher-dimensional generalizations have more variables
but seem not directly relatable  to higher polylogarithms $Li_{n}(x)$.

\subsection{Prior work}\label{sec14}

 There is a long history of work on analytic continuation of  the Lerch transcendent.
  After Lerch's 1887 work, in
 1889  Jonqui\`{e}re \cite{Jo1889} studied the two variable function
 $\zeta(s, x) := \sum_{n=1}^{\infty} \frac{x^n}{n^{s}},$
obtaining various contour integral representations and a functional
equation, with $s$ and $x$ allowed to take some complex values;
this is the function considered by Hilbert \cite{Hi00}.
 In 1906 Barnes \cite{Ba06} studied the Lerch transcendent
with some restrictions on its variables, and noted some aspects of  its multivalued
nature.    In the period 1900-2000 there  was much further  work
on these functions obtaining  analytic continuations
in two of the  variables, omitting one of either $a$ or $x$ (resp. $c$ for the Lerch zeta function)
while pursuing other objectives, such as functional equations, which we pass over here.

Concerning analytic continuation in three variables,
in   2000  Kanemitsu, Katsurada and Yoshimoto \cite{KKY00} obtained an
  analytic continuation of the Lerch transcendent in three variables to a single-valued
 function on various large domains in $\CC^3$. These authors  also obtained formulas for
 special values at negative integers,  related to those given below
 in \S5.   They did not  address the issue of further analytic continuation to a multivalued function.
  In 2008 Guillera and Sondow \cite{GS08} also give a single-valued analytic continuation of
 $\Phi(s, z, c)$ for certain ranges of $(s, z, c)$, mostly restricting $c$ to be real-valued.
 Very recently Costin and
Garoufalidis \cite{CG09}   obtained a multivalued analytic
continuation for  the function $\zeta(x, s)$, calling it the   ``fractional polylogarithm'' 
and denoting it
$Li_{\alpha}(x) =  \sum_{n=1}^{\infty} \frac{x^n}{n^{\alpha}}$
in variables $(\alpha, x)$
on a cover of $\CC \times ( \PP^1 ( \CC ) \smallsetminus \{0,1,\infty \} ) $;
such a continuation appears here as a special case of Theorem~\ref{Nth23}.
  Vepstas \cite{Ve07}  also  obtained results applicable to
analytic continuation of the fractional polylogarithm.

The detailed multivalued nature of the Lerch zeta function $\zeta(s, a, c)$ itself
in all variables appears to have been first worked out in  Part II (\cite{LL2}).
We note that an old approach of Barnes \cite{Ba06}
 might be further developed
 to  effect an analytic continuation
of the Lerch transcendent in three variables.

 Polylogarithms have their own independent history, as generalizations of the logarithm,
 and trace back to work of Euler \cite{E736}, cf. \cite[Sect. 2.4]{Lag13}.
 Much classical work on them is presented in the book of   Lewin \cite{L1}
 and in the volume \cite{L2}.
  The appearance of the dilogarithm in many contexts in mathematics
 and physics is described in Zagier \cite{Zag88}, \cite{Zag07}
 and Oesterl\'{e} \cite{Oes92}. It appears in the computation of
 volumes of hyperbolic tetrahedra,
 and from there to define invariants
 of hyperbolic manifolds, related to  its functional equations, see Neumann \cite{Neu04}.
 Polylogarithms appear in the theory of motives, in
iterated integrals and  mixed motives,
see the discussion in Bloch \cite{Bl91} and Hain \cite{Ha91}.
Generalized polylogarithms given by iterated integrals are considered in
Minh et all \cite{MPV00}
and Joyner \cite{Joy10}.
They appear in Beilinson's conjectures on special values of $L$-functions,
in defining regulators (\cite{Be84}),
cf.  Beilinson and Deligne \cite{BD94}, Huber and Wildeshaus \cite{HW98}.
Geometric versions of polylogarithms
have been formulated (Goncharov \cite{Gon94}, \cite{Gon95}, Cartier \cite{Ca02}).

The multivaluedness of the polylogarithms
encodes period data, and also data on mixed
Hodge structures. In addition $p$-adic and $\ell$-adic analogues
of polylogarithms have been introduced
and studied (Coleman \cite{Co82}, Besser \cite{Be02},
 Furusho \cite{Fur04},
\cite{Fur07}, Nakamura and Wajtkowiak \cite{NW02},
and Wojtkowiak \cite{Woj04}, \cite{Woj05a}, \cite{Woj05b}.)
In another direction, an  exponentiated  quantum deformation of
the dilogarithm, the {\em quantum dilograrithm},
which satisfies a deformed functional equation, was proposed by Faddeev and Kashaev \cite{FK94})
in 1994.  It has since been much studied, see the survey of Kashaev and Nakanishi \cite{KN11}.
Certain dilogarithm identities play a role in  integrable models and  in
conformal field theory (Nahm et al \cite{NRT93}, Kirillov \cite{Kir94}, \cite{Kir95}).
Motivic realizations of polylogarithms are discussed in  Wildeshaus \cite{Wil97}.

There has been much other work on   the Lerch zeta function and Lerch transcendent,
treatments of which can be found in books of Erdelyi et al \cite[Sect. 1.10-1.12]{ErdH53},
 Lauren\v{c}ikas and Garunk\v{s}tis \cite{LG02},
Srivastava and Choi \cite[Chap. 2]{SC01},  Kanemitsu and Tsukada \cite[Chaps. 3-5]{KT07}
and Chakraborthy, Kanemitsu and Tsukada\cite[Chap. 3]{CKT10}.

%
%
\subsection{Present work}\label{sec15}

  From the   viewpoint of earlier work, the main point of this paper is to
make an explicit study of 
the multivaluedness of the analytic continuation of the Lerch transcendent,
and to determine the effects of this analytic continuation on its other algebraic structures.
\begin{enumerate}
\item[(i)]
 We determine explicit formulas
for the monodromy functions and their behavior under
specialization. On the conceptual side, these formulas illuminate a new way in which
the the non-positive integer values $s= -n \le 0$ are ``special values'' of the associated functions,
namely  they  are distinguished points
 in the $s$-parameter space in the sense that  these are the unique values where all  monodromy
 functions vanish identically, see Theorem \ref{th24}.
 \item[(ii)]
  It is well known that the  special values $\zeta(-n)$ ($n \ge 0$)
 at negative integers are rational numbers whose arithmetic properties allow   
 $p$-adic interpolation \footnote{The recipe of Kubota and Leopoldt \cite{KL64} permitting
  interpolation requires  that suitable correction factors be applied to
 the values $\zeta(-n)$, related to the Euler factor at the prime $p$. }
 which leads to the construction of $p$-adic $L$-functions.
  In  Section \ref{sec6}  we show that at these special values $s=-n$  one can
 recover  information from nearby nonsingular strata values $z\Phi(-n, z, c)$
(taking limits  $c \to 0^{+}$) that is sufficient to
 interpolate  $p$-adic $L$-functions; this is achieved using periodic zeta function values.
\item[(iii)]
From the viewpoint of polylogarithms
and  iterated integrals,  we show that under specialization
this Lerch transcendent  provides a complete set of
solutions to a one-parameter Fuchsian deformation of the polylogarithm
differential equation in the parameter $c$,
and we determine its monodromy representation.
This deformation of the polylogarithm may in future shed interesting light
on its behavior.
\end{enumerate}

The results of this paper suggest  that further study be made of the limiting structure of functional
equations in a neighborhood of  the polylogarithm point $c=1$. As noted in Section \ref{sec13} there are two sources of functional
equations, which relate these functions for different values of the $s$-parameter.
The three-term and  four-term functional equations come from number theory,
and  represent a generalization of
the  functional
equations for the Riemann zeta function, relating functions with values $s$ and $1-s$.
These functional equations break down at the integer parameter values $c=m \ge 1$.
But exactly at those $c$-values these functions  satisfy many additional
functional equations which relate function values at various $s$-parameter values shifted by integers,
some of which have geometric meaning.
It seems of interest to determine how these 
additonal functional equations deform in the $c$-parameter.

The extra variables  in the Lerch transcendent potentially make visible
new  connections between these number-theoretic and geometric  viewpoints.
The variable $z$ in the Lerch transcendent, added to the Hurwitz
zeta function variable $c$, gives it the  property of
satisfying a linear partial differential equation,
together with raising and lowering operators, whose form
connects to mathematical physics.
On the number theory side, the  Lerch transcendent  may potentially  yield  new information about 
the Hurwitz zeta function and the Riemann zeta function, even though these 
functions live on singular strata.   
This potentially may occur  by explicit limiting processes (for certain parameter ranges),
using also regularization methods, 
and perhaps  through analysis of the indirect influence of its monodromy.
In section \ref{sec9} we suggest a number of   other directions
for further work.

%
%
%

\section{Summary of main results}\label{sec2}
\setcounter{equation}{0}


 We obtain
    the analytic continuation and monodromy functions for the Lerch transcendent
    as a multivalued function defined on the complex $3$-fold
\beql{N201}
\sN :=
\{ (s, z, c) \in \CC \times (\PP^1(\CC)\smallsetminus \{0,1, \infty\}) \times
(\CC \smallsetminus \ZZ)\}.
\eeq
\noindent
The  universal cover
$\tsN$ of $\sN$
can be identified with  homotopy classes $[\gamma]$
of paths $\gamma$ in $\sN$, and we refer to a point
$[\gamma] \in \tsN$, where the curves start  at the fixed base point
$\bx'_0 =(\frac{1}{2}, -1, \frac{1}{2})$ and end at a point $\gamma(1)$ lying
a  point $(s, z, c) \in \sN$.
The point $[\gamma] \in \tsN$ is sometimes
written as $(s, z, c, [\gamma])$ to stress the end point $(s, z, c)$ of $\gamma$ in $\sN$
of the path $[\gamma].$.
The universal covers for other manifolds containing $\sN$ and for covers of the
extended manifold $\sN^{\sharp}$ are defined similarly.

Our notation used here for paths $\gamma$  generalizes the notation  used in part II, which was
restricted to be a loop having $\gamma(0)= \gamma(1) = \bx_0=(\frac{1}{2}, \frac{1}{2}, \frac{1}{2}),$
with associated homotopy class  $[\gamma] \in \pi_1(\sM, \bx_0)$.
In part II we  wrote $Z(s, a, c, [\gamma])$  to denote the function element
centered at the endpoint $\gamma(1)$, with $(s,a, c)$ denoting  local coordinates in a
neighborhood of the endpoint  of the  loop $\gamma$. Reaching
the point $(s,a,c)$ from $\gamma(1)$ can be thought of as following
an  additional path $\gamma^{'}$ from
$\gamma^{'}(0)=\gamma(1)$ to  $\gamma^{'}(1) = (s, a, c)$ that remains  in a simply connected region obtained by
cutting the manifold $\{ (s, a, c) \in  \CC \times (\CC \smallsetminus \ZZ) \times (\CC\smallsetminus \ZZ)$\} along the
lines $\{ a= m+it: t \le 0\}$ for $m \in \ZZ$ and similarly in the $c$-variable.
In this paper  $\gamma$ denotes a path, to be thought of as the analogue of the composed path
$\gamma \circ \gamma^{'}$, paths being composed left to right, as in Hatcher \cite[p. 26]{Hat01}.
Thus $\gamma$ need not be a closed path.


\begin{theorem}\label{th21} {\em (Lerch Transcendent Analytic Continuation)}

(1) The Lerch transcendent $\Phi(s, z,c)$ on $\sN$ analytically continues to a
single-valued holomorphic function $\tZ([\gamma])= \tZ(s, z, c, [\gamma])$ on the universal
cover $\tsN$ of $\sN$.

(2) The function
$\tZ(s, z, c, [\gamma])$ becomes single-valued on a  two-step solvable regular (i.e. Galois)  covering manifold
$\tsN^{solv}$ of $\sN$, which can be taken to be the manifold
fixed by the second commutator subgroup  of $\pi_1(\sN, \bx'_0).$
\end{theorem}

This is established in Section \ref{sec3}, where
Theorem~\ref{th51} and Theorem~\ref{th52}  give more detailed statements,
which imply the result above.
In particular we show that all
monodromy functions vanish identically on
a certain normal subgroup $\Gamma^{'}$ of $\pi_1(\sN, \bx'_0)$
 which contains the second derived
subgroup (second commutator subgroup) $(\pi_1( \sN, \bx'_0))^{'' }$.

Our next result shows that the singularities at $c \in \ZZ_{\ge 1}$ are removable,
giving an analytic continuation to a solvable covering of $\sN^{\#}$,
as follows (cf. Theorem~\ref{Nth23}).


\begin{theorem}\label{th22} {\em (Lerch Transcendent Extended Analytic Continuation)}

(1) The Lerch transcendent $\Phi(s, z,c)$ analytically continues to a
single-valued holomorphic function $\tZ(s, z, c, [\gamma])$ on the universal
cover $\tsN^{\sharp}$ of the manifold
\beql{N202}
\sN^{\sharp} :=
\{ (s, z, c) \in \CC \times (\PP^1(\CC)\smallsetminus \{0,1, \infty\}) \times
(\CC \smallsetminus \ZZ_{\le 0})\}.
\eeq

(2) The function
$\tZ(s, z, c, [\gamma])$ becomes single-valued on a  $2$-step solvable covering manifold
$\tsN^{\sharp, solv}$ of $\sN^{\sharp}$,
which can be taken to be  the covering that is fixed by the second commutator subgroup of
$\pi_1(\sN^{\sharp}, \bx'_0).$
\end{theorem}

In Section \ref{sec4} we  observe that the Lerch transcendent
and its  analytic continuation
satisfies two differential-difference equations
and a linear partial differential equation,
as follows (cf. Theorem~\ref{th61}).


\begin{theorem}\label{th23}
{\em (Lerch Transcendent Differential-Difference Operators)}

(1) The Lerch transcendent $\Phi(s, z, c)$ satisfies
the differential-difference equations
\beql{N231}
 \left( z \frac{\partial}{\partial z} + c\right) \Phi(s, z, c)= \Phi(s-1, z, c),
\eeq
and
\beql{N232}
  \frac{\partial}{\partial c} \Phi(s, z, c) = -s \Phi(s+1, z, c).
\eeq
These differential-difference equations are also  satisfied by
the analytic continuation \\
$\tZ (s, z, c, [\gamma])$
of the Lerch transcendent $\Phi(s, z, c)$ on the universal cover
$\tilde{\sN}$.

(2) The Lerch transcendent $\Phi(s, z, c)$ satisfies the linear
partial differential equation
\beql{N233}
 \left( z \frac{\partial}{\partial z} + c \right) \frac{\partial}{\partial c} \Phi(s,z,c)= s \Phi(s, z,c).
\eeq
The
analytic continuation $\tilde{Z}(s, z, c, [\gamma])$
satisfies this equation on the universal cover
$\tilde{\sN}$ of $\sN$.

(3) For each $[\tau ] \in \pi_1 ( \sN , \bx'_0 )$ its associated
 monodromy function
  $M_{[\tau]} ( \tZ )(s, z, c, [\gamma])$
of the Lerch transcendent satisfies on
$\tilde{\sN}$ the two differential-difference equations
and the linear partial differential equation (\ref{N233}).
\end{theorem}

 As explained in Section \ref{sec31}, the fundamental group $\pi_1(\sN,\bx'_0 )$ is the product of
 $\pi_1(\PP^1(\CC) \smallsetminus \{0, 1, \infty\}, -1)$, a free group on two generators
 $[Z_0]$ and $[Z_1]$, and $\pi_1(\CC \smallsetminus \ZZ, \frac{1}{2})$,
 a free group on generators $[Y_n]$ for $n \in \ZZ$. We next  determine the structure of
the vector spaces $\sW_s$, for fixed $s \in \CC$,  spanned by all the branches of the
multivalued analytic continuation of the Lerch transcendent,
over a neighborhood of a given point $(s, z, c) \in \sN$,
specifying  a set of generators for these spaces.
We define the space $\sW_s$  to be a (generally infinite)
direct sum of one-dimensional vector spaces given
by particular monodromy generators, see Section \ref{sec4} . There is a generic basis for
$s \not\in\ZZ$ and for $s \in \ZZ$ there are linear relations among the generators,
effectively reducing their number.
The following  result is established as Theorem~\ref{th62}.

\begin{theorem}\label{th24}
{\em ( Lerch Transcendent Monodromy Space)}
The Lerch transcendent monodromy space $\sW_s$ at $s$ depends on the parameter
$s \in \CC$ as follows.
\begin{itemize}
\item[(i)] (Generic case)
If $s \not\in \ZZ$, then $\sW_s$ is an infinite-dimensional vector space,
and has as a basis the set of functions
\beql{N243}
\{ M_{[Z_0]^{-k} [Z_1] [Z_0]^k}^s (\tZ ) : k \in \ZZ \} \cup
\{ M_{[Y_n]}^s ( \tZ ) : n \in \ZZ \} \cup \{ \tZ^s\} ~.
\eeq
\item[(ii)]
If $s = m \in \ZZ_{> 0}$, then $\sW_m$ is an infinite-dimensional
vector space, and has as a basis the set of functions
\beql{N242}
\{
M_{[Z_0]^{-k} [Z_1 ] [Z_0]^k}^m ( \tZ ) : k \in \ZZ \} \cup \{ \tZ^s\}~.
\eeq

\item[(iii)]
If $s = -m  \in \ZZ_{\le 0}$, then all Lerch transcendent monodromy functions
vanish identically, i.e.
\beql{N241}
M_{[\tau]}^{-m} ( \tZ ) =0 \quad\mbox{for all}\quad
[\tau ] \in \pi_1 ( \sN , \bx'_0 ) ~.
\eeq
Thus $\sW_{-m} = \CC\tZ^s$ is a one-dimensional vector space.
\end{itemize}
\end{theorem}

In Sections \ref{sec5}-\ref{sec8} we specialize the Lerch transcendent
variables $(s,z,c) $ to cases where  $s=m$ is an integer.
These are exactly the cases where the monodromy functions satisfy
``non-generic'' linear relations. We show that as functions
of the two complex variables $(z,c)$ we obtain further analytic
continuation into the singular strata of the three-variable analytic
continuation given in Section \ref{sec3}. Here the  singular strata  correspond to
$z=0, 1$ and/or  $c \in \ZZ_{\le 0}$. For convenience we state results
in terms of the {\em extended polylogarithm} $Li_{m}(c, z)= z \Phi(m, c, z)$.

In Section \ref{sec5} we treat the case where $s=-m \le 0$ is a non-positive integer.
Here we may note that
 $Li_{-m}(c, z)$ satisfies the ordinary differential equation of order $m+1$ in
the $c$-variable
\beql{N110c}
\frac{d^{m+1}}{dc^{m+1}}Li_{-m}(c, z) =0,
\eeq
which is independent of $z$. It
 implies that $Li_{-m}(c, z)$ is necessarily a polynomial in $c$ of degree at most $m$, having
coefficients which are functions of $z$.
It is  known  that these coefficients are rational functions of $z$;
this  follows from the observation of
Apostol \cite{Ap51} in 1951
 that the function
$Li_{-m}(z,c)$  extends
to a rational function in the $(z, c)$ variables.
 Here we determine formulas for these rational functions  using
  the differential equation, as follows (Theorem ~\ref{th63}).

\begin{theorem}\label{th25}
{\em ($c$-Deformed Negative Polylogarithms)}
For $s = -m \in \ZZ_{\le 0}$ the function $Li_{-m} (z,c)$ analytically continues to  a rational
function of $z$ and $c$ on $\PP^1 (\CC ) \times \PP^1 ( \CC )$.
Here $Li_0 (z, c) = zq_0(z)$ and
\beql{251}
Li_{-m} (z,c) = z \left( \sum_{k=0}^m {\binom{m}{k}}  c^k q_{m-k} (z)\right)\, ,
\quad\mbox{for}\quad m \ge 1,
\eeq
in which the $q_m(z)$ are rational functions of $z$  given by
$q_0 (z) = \frac{1}{1-z}$
and
\beql{253}
q_{m+1} (z) = z \frac{d}{dz} (q_m (z)) \quad\mbox{for}\quad m \ge 0 ~.
\eeq
\end{theorem}

We also determine recursion
relations for these rational functions, and show they have a reflection
symmetry $z^{m+1} r_{m}(\frac{1}{z}) = r_m(z)$. (Theorem ~\ref{Zle44}).
The rational function $Li_{-m}(c, z)$ takes well-defined values on the Riemann sphere for all
 $(z, c) \in \PP^{1}(\CC) \times \PP^1(\CC)$, and thus extends to the
 singular strata  regions given by the complex hyperplanes $c \in \ZZ_{ \le 0} $, resp.
 $z=0$. On
the  singular stratum  $z=0$ these functions take finite values,
 but on the `singular stratum  $z=1$, they have
a nontrivial polar part, and take the constant value $\infty \in \PP^1(\CC)$.

In Section \ref{sec6} we consider the double specialization when $s=-m \le 0$ is
a nonpositive integer and $c=0$. We show the function then
agrees with the analytically continued value
of the periodic zeta function $F(a, s) := \sum_{n=1}^{\infty} e^{2 \pi i na}{n^{-s}}$,
when $z=e^{2\pi i a}$ lies on the unit circle (Theorem \ref{th43}).


\begin{theorem}\label{th26}
{\em (Periodic Zeta Function Special Values)}

(1) For $z \in \PP^1(\CC) \smallsetminus \{0, 1, \infty\}$ and $s=-m \in \ZZ_{\le 0}$  there holds
\beql{261}
Li_{-m} (z, 0) = \lim_{c \to 0^{+}} Li_{-m}(z, c) = \lim_{c \to 0^{+}} z \Phi(-m, z, c),
\eeq
where the limit is taken through values of $c$ in $0< \Re(c)<1$.

 (2) For $0 < \Re(a)  < 1$ the periodic zeta function
$F(a,s) = \sum_{n=1}^\infty  \frac{e^{2 \pi ina}}{n^s}$
analytically continues to an entire function of $s$.
In particular, for $s= -m \in \ZZ_{\le 0}$ 
there holds
\beql{262}
F(a, -m) =   e^{-2 \pi i a} Li_{-m}(e^{2 \pi i a}, 0)= q_m(e^{2 \pi i a}).
\eeq
\end{theorem}

This equality \eqn{261}  is non-trivial because it involves a limiting procedure,
since the point $c=0$ lies in a singular stratum
of the analytic continuation of the Lerch zeta function given in part II.
This equality  permits one
to construct $p$-adic L-functions by interpolation from values of the Lerch zeta function.
In contrast, it appears that one cannot recover values of the Hurwitz zeta function
at $s=-m$  directly from the Lerch transcendent
by such a limiting procedure, letting $a \to 0^{+}$ or $a \to 1^{-}$;
these limits do not exist. \medskip

In Section \ref{sec7}  we treat the one-variable specialization where  $s=m \geq 1$ is a positive
integer. We state results in terms of  the function of two variables $Li_{m}(z, c)$,
observing that it satisfies the (slightly different)  linear PDE
\beql{277}
   \left( z \frac{\partial}{\partial z}  \frac{\partial}{\partial c}+ (c-1) \frac{\partial}{\partial c}\right) Li_{m}(z, c)=
-m Li_{m}(z, c).
\eeq
We obtain  the following result as Theorem~\ref{th64}.

\begin{theorem}\label{th27}
{\em ($c$-Deformed Polylogarithm Analytic Continuation)}
For each positive integer $s=m \ge 1$ the function $Li_m(z,c)$ has a
meromorphic continuation in  two variables $(z,c)$ to the
universal cover
of $~(\PP^1 (\CC ) \smallsetminus \{0,1,\infty \} ) \times \CC$.
For fixed $\tilde{z}$ on the universal cover,
 this function  is meromorphic as a function of  $c \in \CC$,
with its singularities consisting of  poles of exact order $m$ at each of the points $c \in \ZZ_{\le 0}.$
\end{theorem}

\noindent This result gives an analytic continuation  of $Li_{m}(z, c)$ to negative integer $c$,
which are points that  fall in the  singular strata outside of the analytic continuation
given in Section \ref{sec3}

In Section \ref{sec8}  we consider the double specialization, in which $s= m \ge 1$ is
a positive integer and $c$ is fixed.  Here  we consider
the  functions $Li_{m}(z,c)$ of one variable $z$ and show they
have monodromy functions which are a deformation of the monodromy
of the polylogarithm $Li_{m}(z)$, viewing $c$ as a deformation parameter,
with $c=1$ giving the polylogarithm- this value of $c$ is inside the three-variable
analytic continuation given in \S3.
Namely  we observe that $F(z):=Li_{m}(z, c)$ satisfies an ordinary differential
equation $D_{m+1}^c F(z) =0$, where
$$
D_{m+1}^c = z^2 \frac{d}{dz}
 \left( \frac{1-z}{z}\right)\left( z \frac{d}{dz} + c-1\right)^m \in \CC[z, \frac{d}{dz}].
$$
In Theorem~\ref{Nth81} we obtain the following result.

\begin{theorem}\label{th28}
{\em ($c$-Deformed Polylogarithm Ordinary Differential Equation)}\\
Let $m \in \ZZ_{\ge 0}$ and let $c \in \CC$
be fixed.

(1) The function
$ F(z)= Li_{m}(z,c)$ satisfies the  ordinary differential equation
\beql{281}
D_{m+1}^c F(z) =0,
\eeq
 where
$D_{m+1}^c \in \CC[z, \frac{d}{dz}]$ is the linear ordinary differerential operator
\beql{282}
D_{m+1}^c := z^2 \frac{d}{dz} \left(\frac{1-z}{z} \right)\left( z \frac{d}{dz} + c-1\right)^m
\eeq
 of order $m+1$.

(2) The operator $D_{m+1}^c$  is a Fuchsian operator for all $c \in \CC$.
For each $c \in \CC$ its singular points on the Riemann sphere are all regular and are
 contained in the set  $\{0, 1, \infty\}$.

(3) A basis of solutions of $D_{m+1}^c$  for $c \in \CC \smallsetminus \ZZ_{\le 0}$
is,
for $z \in \CC \smallsetminus \{ (-\infty, 0] \cup [1, \infty) \}$,
 given by
\beql{283}
\sB_{m+1,c} := \{ Li_{m}(z,c), z^{1-c}(\log z)^{m-1}, z^{1-c}(\log z)^{m-2}, \cdots, z^{1-c}  \}.
\eeq
(For $c \in \ZZ_{\le 0}$ the function $Li_{k,c}(z)$
is not well defined.)

(4) A basis of solutions of $D_{m+1}^c$ for $c= -k \in \ZZ_{\le 0}$ is,
for  $z \in \CC \smallsetminus \{ (-\infty, 0] \cup[1, \infty) \}$,
given by
\beql{284}
\sB_{m+1,c}^{\ast} := \{ Li_{m}^{*}(z,-k), z^{1-c}(\log z)^{m-1}, z^{1-c}(\log z)^{m-2}, \cdots, z^{1-c}  \},
\eeq
in which
\beql{285}
Li_{m}^{*}(z, -k) := \sum_{{n=0}\atop{n \ne k}}^{\infty} \frac{z^{n+1}}{(n-k)^m} +
\frac{1}{m!} z^{k+1}(\log z)^m.
\eeq
\end{theorem}

We next study the monodromy representation of the fundamental group \\
$\pi_1(\PP^1(\CC) \smallsetminus \{0, 1, \infty\}, -1) = \langle [Z_0], [Z_1]\rangle$
on the multivalued solutions of this differential equation.
The  associated monodromy representation is finite-dimensional, of
dimension $m + 1$   independent of $c \in \CC$.
We show that  the image of the
monodromy representation lies in
a Borel subgroup of $GL(m + 1, \CC)$, and lies in
a  unipotent subgroup exactly
when $c \in \ZZ.$ The result splits into two cases,
one for the ``non-singular'' strata  values of $c$ and the other for the
singular strata values $c \in \ZZ_{\le 0}.$ This first case is
given in Theorem~\ref{thN82}, as follows.

\begin{theorem}\label{th29}
{\em ($c$-Deformed Polylogarithm Monodromy-Nonsingular Case)}\\
For each integer $s=m \ge 1$ and each $c \in \CC \smallsetminus \ZZ_{\le 0}$,
the monodromy action for $D_m^c$ of $
\pi_1(\PP^1(\CC) \smallsetminus \{0, 1, \infty\}, -1) = \langle [Z_0], [Z_1]\rangle
$
acting on the basis
$$
\sB_{m+1,c} := \{ Li_{m}(z,c), z^{1-c}(\log z)^{m-1}, z^{1-c}(\log z)^{m-2}, \cdots, z^{1-c}  \}
$$
is given by
\beql{291}
\rho_{m,c} ([Z_0]) := \left(
\begin{array}{cccccc}
1 & 0 & 0 & \cdots & 0 & 0 \\
0 & e^{- 2 \pi ic} & \frac {2 \pi i}{1!}e^{-2 \pi ic} & \cdots &
\frac{(2 \pi i)^{m - 2}}{(m - 2)!}e^{-2 \pi ic}  & \frac{(2 \pi i)^{m - 1}}{(m - 1)!}e^{-2 \pi ic} \\
\vdots & \vdots & \vdots & ~ & \vdots & \vdots \\
0 &  0 &0 & \cdots & e^{-2\pi i c} & \frac{2\pi i}{1!}e^{- 2\pi i c}  \\
0 & 0 & 0 & \cdots & 0 & e^{-2 \pi ic} \\
\end{array}
\right) \,
\eeq
and
\beql{292}
\rho_{m,c} ([Z_1]) := \left( \begin{array}{cccccc}
1 & -2\pi i & 0 & \cdots & 0 & 0 \\
0 & 1 & 0 & \cdots & 0 & 0 \\
\vdots & \vdots & \vdots & \vdots & \vdots & \vdots \\
0 & 0 & 0 & \cdots & 1 & 0 \\
0 & 0 & 0 & \cdots & 0 & 1
\end{array}
\right) \, .
\eeq
The image of $\rho_{m,c}$ falls in a Borel subgroup of $GL(m+1, \CC)$. The image
 is unipotent if $c \in \ZZ_{>1}$, and is quasi-unipotent if $c \in \QQ \smallsetminus \ZZ.$
\end{theorem}

The special case $c=1$ of this result corresponds to the polylogarithm case considered by Ramakrishnan
\cite{Ra81}, \cite{Ra82}, \cite{Ra89}. (We remark at the end of \S8
on the issue of reconciling  our formulas with those
of Ramakrishnan.) \\

We also obtain the monodromy in the  singular strata cases $c \in \ZZ_{\le 0}$,
outside the analytic continuation in \S3. Here
the monodromy representation exhibits a discontinuous jump from
the ``non-singular'' strata values,
that seems  unresolvable by a
 change of basis of $\langle [Z_0], [Z_1]\rangle$ (cf. Theorem~\ref{thN83}).

\begin{theorem}\label{th210}
{\em ($c$-Deformed Polylogarithm Monodromy-Singular Case)}\\
For each integer $s=m \ge 1$ and each $c \in  \ZZ_{\le 0}$,
the monodromy action for $D_m^c$ of the homotopy group
$
\pi_1(\PP^1(\CC) \smallsetminus \{0, 1, \infty\}, -1 )= \langle [Z_0], [Z_1] \rangle
$
 acting on the basis
$$
\sB_{m+1,c}^{\ast} := \{ Li_{m}^{*}(z,c), z^{1-c}(\log z)^{m-1}, z^{1-c}(\log z)^{m-2}, \cdots, z^{1-c}  \}
$$
is given by
\beql{296}
\rho_{m,c} ([Z_0]) := \left(
\begin{array}{cccccc}
1 &  \frac{2 \pi i}{1!}& \frac{(2\pi i)^2}{2!} & \cdots &  \frac{(2 \pi i)^{m - 1}}{(m - 1)!} &
 \frac{(2 \pi i)^{m}}{m!}\\
0 & 1 & \frac {2 \pi i}{1!} & \cdots &
 \frac{(2 \pi i)^{m - 2}}{(m - 2)!} & \frac{(2 \pi i)^{m - 1}}{(m - 1)!} \\
\vdots & \vdots & \vdots & ~ & \vdots & \vdots \\
0 & 0 & 0 & \cdots & 1 & \frac{2\pi i}{1!} \\
0 & 0 & 0 & \cdots & 0 & 1 \\
\end{array}
\right)
\eeq
and
\beql{297}
\rho_{m,c} ([Z_1]) := \left( \begin{array}{cccccc}
1 & -2\pi i & 0 & \cdots & 0 & 0 \\
0 & 1 & 0 & \cdots & 0 & 0 \\
\vdots & \vdots & \vdots & \vdots & \vdots & \vdots \\
0 & 0 & 0 & \cdots & 1 & 0 \\
0 & 0 & 0 & \cdots & 0 & 1
\end{array}
\right) ~.
\eeq
\end{theorem}


\paragraph{\bf Remark on notation for logarithms.} We
will need  two different versions  of the restriction
of the logarithm to a single-valued function on a cut plane. The
{\em principal branch} $\log u$ makes a cut on the negative real axis, and
the {\em semi-principal branch} $\llog u$ makes a cut on the positive real axis.
In both cases the cut line is connected to the upper half-plane, so that
$$\log 1 = \llog 1=0, ~~~\log (-1) = \llog(-1) = \pi i.$$
The two branches differ in the values assigned to the logarithm in the
lower half-plane.
However if one starts from a base point in the upper half-plane of $u$,
the  multivalued extensions
of these two functions agree.
We use  $\Re(s)$ and $\Im(s)$ for
the real and imaginary parts of a complex variable $s$.
Finally, compositions of paths  are read left to right, so $g \circ h$ with $g(1)= h(0)$ means follow path $g$ then
follow path $h$, see Hatcher \cite[p. 26]{Hat01}.

%
%
%

\section{Analytic continuation of Lerch transcendent}\label{sec3}
\setcounter{equation}{0}

In this section we analytically continue
the Lerch transcendent $\Phi(s,z,c)$ in
the three complex variables $(s,z,c)$ to
a multivalued function and  compute its monodromy functions, using
the results of part II.

%
%
%

\subsection{Analytic continuation of Lerch zeta function}\label{sec31}

We recall some details of the analytic continuation of the Lerch zeta function
given in part II \cite{LL2}. This analytic continuation first extended the
Lerch zeta function to a single-valued function of three complex variables on
 the {\em extended fundamental polycylinder}
$$
\tilde{\Omega}:= \{ s : s \in \CC \} \times \{a :  0 < \Re (a) < 1 \} \times
\{c:  0 < \Re (c) < 1 \} ,
$$
This region is   simply connected, and is invariant under the $\ZZ/4\ZZ$-symmetry of the
functional equation $(s, a, c) \mapsto (1-s, 1-c, a)$.

To understand the multivalued structure,
for fixed $s \in \CC$  we  introduce the ``elementary functions'' 
$$
\phi_n^s(a, c) = \begin{cases}
\frac{e^{-2\pi i na}}{(n-c)^s}
& \quad{if} \quad  n \ge 1, \\
& \\
 \frac{e^{-2\pi i na}}{ (c-n)^s}
 &  \quad{if} \quad  n \le 0. \\
\end{cases}
 $$
 These functions are well-defined in the extended fundamental polycylinder, taking the principal branch $\log z$ of
 the logarithm,
 with
 $w^s := \exp (s \log w)$, for $\Re(w)>0$.
  When analytically continued in the $c$-variable, $\phi_n^s$  has a branch point of infinite order at $c= n$,
which  corresponds to the  branch locus  $\sV(c= n) := \{ (a, c) : \, a \in \CC, \,  c = n\}$  inside  $ \CC^2$.

We also introduce  the ``elementary functions'' 
$$
\psi_n^s(a, c) :=
\begin{cases}
e^{-2 \pi i a c} \frac{e^{2 \pi i nc} }{(n-a)^{s-1}}
& \quad \mbox{if} \quad n \ge 1,\\
& \\
e^{- 2 \pi i ac}  \frac{e^{2 \pi i n c}}{(a-n)^{s-1}}
& \quad \mbox{if} \quad n \le 0.\\
\end{cases}
$$
When analytically continued in the $a$-variable, $\psi_n^s$  has a branch point of infinite order at $a= n$,
which corresponds to the branch locus  $\tilde{\sV} (a=n) := \{ (a,c): a=n, \,c \in \CC \}$ inside  $\CC^2$.
When we remove all the branch loci we obtain (for fixed $s \in \CC$) the punctured complex surface
$\{ (a, c) \in (\CC\smallsetminus \ZZ) \times  ( \CC \smallsetminus \ZZ)$\}.

The Lerch zeta function has two expansions  in  terms of the elementary functions. The first
\begin{equation}\label{EF1}
\zeta(s, a ,c) = \sum_{n=0}^{\infty} \phi_{-n}^s(a, c)
=\sum_{n=-\infty}^{0} \phi_{n}^s (a, c).
\end{equation}
Here the  right side  of this formula converges conditionally for  $0 < \Re(s)<1$
when  $a, c$ are real with  $0< a, c <1$ (and does not converge for $s$ outside this strip).
We obtain a second expression in terms of elementary functions  using the functional equation for
the Lerch zeta function $Z(s, a, c, [pt])$ on $\tilde{\Omega}$, which is
\begin{equation} \label{EF2}
\zeta(s, a, c) = \sum_{n \in \ZZ} c_n(s) \psi_n^s(a, c),
\end{equation}
which has coefficients $c_n(s)$ given by
\begin{equation}
c_n(s) := (2\pi)^{s-1} \Gamma(1-s) \begin{cases}
e^{-\frac{\pi i (1-s)}{2}} &  \quad \mbox{if} \quad n \ge 1,\\
e^{\frac{\pi i (1-s)}{2}} &  \quad \mbox{if} \quad n \le 0. \\
\end{cases}
\end{equation}
The  right side of \eqref{EF2} converges conditionally for $0< \Re(s)<1$
for real $0< a, c <1$ with $0 < \Re(s) <1$.
More generally the coefficients $c_n(s)$ are meromorphic functions, with simple poles at $s\in \ZZ_{\ge 1}$.

The multivaluedness of the Lerch zeta function is exactly that  inherited from the
individual multivalued members of the right sides of \eqref{EF1} and \eqref{EF2},
as we now explain.
To obtain  the analytic continuation in three variables, we fix as base point
 $\bx_0 := (\frac{1}{2}, \frac{1}{2}, \frac{1}{2}) \in \tilde{\Omega}$
 and consider  a closed path  $\gamma$ starting from the base point on the base manifold
 $
 \sM = \{(s, a, c) \in \CC \times (\CC \smallsetminus \ZZ) \times
 (\CC \smallsetminus \ZZ)\}.
 $
 We view the analytic continuation on the universal cover $\tilde{\sM}$ starting
 from the base point. This
 multivalued analytic continuation defines a sequence of
 single-valued analytic functions
 $Z( s, a, c, [\gamma])$
  defined on the extended fundamental polycylinder  $\tilde{\Omega}$,
 which depend only on the homotopy class  $[\gamma] \in \pi_1(\sM, \bx_0)$.
 Here the trivial loop $[pt]$ gives
 $$
 Z(s, a, c, [pt]) = \zeta(s, a, c)
 $$
  The full analytic continuation outside $\tilde{\Omega}$ is obtained by extending these functions to
 a simply connected domain covering the  whole region $\sM$
 by making  a series of cuts in the $a$-plane and $c$-plane separately.

 We recall from part II generators for the homotopy group $\pi_1( \sM, \bx_0)$.
 In the $a$-variable, the generators
for $\pi_1(\CC \smallsetminus \ZZ, a= \frac{1}{2})$ are $\{ [X_n]: n \in \ZZ\}$ in
which $X_n$ denotes
a path from base point $a= \frac{1}{2}$ that lies
entirely in the upper half-plane to the point $a= n + \epsilon i$, followed by
a small counterclockwise oriented loop of radius $\epsilon$ around the
point $a=n$, followed by return along the path.
The generators $[X_n]$ are pictured in Figure \ref{fg53b}.

%
%
\begin{figure}[H]
\centering
\scalebox{1} 
{
\begin{pspicture}(0,-1.61)(7.51,1.61)
\psline[linewidth=0.02cm](1.2,1.6)(1.2,-1.6)
\psline[linewidth=0.02cm](0.0,0.0)(4.8,0.0)
\psdots[dotsize=0.12](1.2,0.0)
\psdots[dotsize=0.12](2.8,0.0)
\usefont{T1}{ptm}{m}{n}
\rput(4.71,1.285){$X_n$}
\usefont{T1}{ptm}{m}{n}
\rput(0.97,-0.255){$0$}
\usefont{T1}{ptm}{m}{n}
\rput(2.81,-0.275){$1$}
\psdots[dotsize=0.12](4.4,0.0)
\psdots[dotsize=0.12](2.0,0.0)
\usefont{T1}{ptm}{m}{n}
\rput(1.99,-0.365){$a=\frac{1}{2}$}
\usefont{T1}{ptm}{m}{n}
\rput(4.41,-0.295){$2$}
\psline[linewidth=0.02cm](4.8,0.0)(5.3,0.0)
\psline[linewidth=0.02cm,linestyle=dotted,dotsep=0.03cm](5.3,0.0)(5.7,0.0)
\psline[linewidth=0.02cm](5.7,0.0)(7.5,0.0)
\psdots[dotsize=0.12](6.7,0.0)
\usefont{T1}{ptm}{m}{n}
\rput(6.69,-0.255){$a = n$}
\psarc[linewidth=0.04](6.7,0.0){0.7}{96.70984}{83.659805}
\psline[linewidth=0.04](6.8,0.7)(6.8,0.96)(1.96,0.96)(1.96,-0.02)
\psline[linewidth=0.04](6.6,0.7)(6.6,0.8)(2.04,0.8)(2.04,0.0)
\psline[linewidth=0.08cm,linestyle=none]{->}(3.3,0.8)(3.9,0.8)
\psline[linewidth=0.08cm,linestyle=none]{->}(5.1,0.96)(3.2,0.96)
\psline[linewidth=0.08cm,linestyle=none]{->}(6.26,-0.54)(6.48,-0.66)
\psline[linewidth=0.08cm,linestyle=none]{->}(7.36,0.32)(7.18,0.52)
\end{pspicture}
}
\caption{Generators $[X_n]$ of $\pi_1(\CC \smallsetminus \ZZ, a=\frac{1}{2})$  in the $a$-plane}
\label{fg53b}
\end{figure}

In the $c$-variable, the homotopy group
$\pi_1 ( \CC \smallsetminus \ZZ , c=\frac{1}{2}) $
has a set of generators $ \{ [Y_n] : n \in \ZZ \}$,
in which $[Y_n]$ denotes
a path from base point $c= \frac{1}{2}$ that lies
entirely in the upper half-plane to the point $c= n + \epsilon i$, followed by
a small counterclockwise oriented loop of radius $\epsilon$ around the
point $c=n$, followed by return along the path.
The generators $[Y_n]$ are pictured in Figure \ref{fg53}.

%
%
\begin{figure}[H]
\centering
\scalebox{1} 
{
\begin{pspicture}(0,-1.61)(7.51,1.61)
\psline[linewidth=0.02cm](1.2,1.6)(1.2,-1.6)
\psline[linewidth=0.02cm](0.0,0.0)(4.8,0.0)
\psdots[dotsize=0.12](1.2,0.0)
\psdots[dotsize=0.12](2.8,0.0)
\usefont{T1}{ptm}{m}{n}
\rput(4.71,1.285){$Y_n$}
\usefont{T1}{ptm}{m}{n}
\rput(0.97,-0.255){$0$}
\usefont{T1}{ptm}{m}{n}
\rput(2.81,-0.275){$1$}
\psdots[dotsize=0.12](4.4,0.0)
\psdots[dotsize=0.12](2.0,0.0)
\usefont{T1}{ptm}{m}{n}
\rput(1.99,-0.365){$c=\frac{1}{2}$}
\usefont{T1}{ptm}{m}{n}
\rput(4.41,-0.295){$2$}
\psline[linewidth=0.02cm](4.8,0.0)(5.3,0.0)
\psline[linewidth=0.02cm,linestyle=dotted,dotsep=0.03cm](5.3,0.0)(5.7,0.0)
\psline[linewidth=0.02cm](5.7,0.0)(7.5,0.0)
\psdots[dotsize=0.12](6.7,0.0)
\usefont{T1}{ptm}{m}{n}
\rput(6.69,-0.255){$c = n$}
\psarc[linewidth=0.04](6.7,0.0){0.7}{96.70984}{83.659805}
\psline[linewidth=0.04](6.8,0.7)(6.8,0.96)(1.96,0.96)(1.96,-0.02)
\psline[linewidth=0.04](6.6,0.7)(6.6,0.8)(2.04,0.8)(2.04,0.0)
\psline[linewidth=0.08cm,linestyle=none]{->}(3.3,0.8)(3.9,0.8)
\psline[linewidth=0.08cm,linestyle=none]{->}(5.1,0.96)(3.2,0.96)
\psline[linewidth=0.08cm,linestyle=none]{->}(6.26,-0.54)(6.48,-0.66)
\psline[linewidth=0.08cm,linestyle=none]{->}(7.36,0.32)(7.18,0.52)
\end{pspicture}
}
\caption{Generators $[Y_n]$ of $\pi_1(\CC \smallsetminus \ZZ, c=\frac{1}{2})$  in the $c$-plane}
\label{fg53}
\end{figure}

The monodromy of the individual elementary functions is simple.
In the case of the $c=n$ singularity a counterclockwise loop traversed $k$ times in the $c$-plane
( with homotopy type $[Y_n]^k$) sends  $\phi_n^s(a, c) \mapsto e^{ 2\pi ik s} \phi_n^s(a, c)$,
while for the $a=n$ singularity a counterclockwise loop traversed $k$ times
in the $a$-plane (corresponding to $[X_n]^k$)
sends  $\psi_n^s(a, c) \mapsto e^{ 2\pi i k(1-s)} \psi_n^s(a, c)$.
In part II    (\cite[Theorem 4.5]{LL2}) we establish that,
the multivaluedness comes directly from the individual terms in \eqref{EF1} and \eqref{EF2}.
We obtain on $\tilde{\Omega}$ that
\begin{equation}\label{EFB1}
 Z(s, a, c, [Y_n]^k) = \zeta(s, a, c) +
 \begin{cases}
 0 & \quad \mbox{if} \quad n \ge 1, \\
 &\\
 (e^{2 \pi i k s} -1) \phi_n^s(a, c) & \quad \mbox{if} \quad n \le 0.
 \end{cases}
 \end{equation}
 For the loops $[X_n]$ we obtain  formulas
 that\footnote{The formula is expressed in a different form in
 \cite[Theorem 4.5]{LL2} using gamma function identities to remove the poles  in $c_n(s)$ at positive integers $s$,
 compare \eqref{519ab} below.}
 can be shown to be equivalent to
\begin{equation}\label{EFB2}
 Z(s, a, c, [X_n]^k) = \zeta(s, a, c) +
 (e^{-2 \pi i k s} -1)c_n(s) \psi_n^s(a, c)  .
  \end{equation}

 Various consequences of the analytic continuation of the Lerch zeta function are deducible
 from \eqn{EFB1} and \eqn{EFB2}:
 \begin{enumerate}
 \item
 $ Z(s, a, c, [\gamma]) = \zeta(s, a, c)$
 whenever $[\gamma]$ is in the commutator subgroup
 $$
 {\bf D}^{(1)}( \pi_1(\sM, \bx_0)):= [ \pi_1(\sM, \bx_0), \pi_1(\sM, \bx_0)]
 $$
 of $\pi_1(\sM, \bx_0)$ (\cite[Theorem 4.6]{LL2}).  Thus the function $Z(s, a, c, [\gamma])$
 is single valued on the maximal abelian cover $\tilde{\sM}^{ab}$  of $\sM$.
 The multivaluedness  on this cover is then described entirely in terms of the winding
 numbers $k_n([\gamma])$
  {(resp. $k_n^{*}([\gamma])$)} of the path $\gamma$ around the
 {sub}manifolds $\tilde{\sV}( a = -n)$
 {(resp.  ${\sV}( c= n)$) while holding $s$ constant}.
  \item
  For  $s= -m  \le 0$ a
nonpositive integer, all
the monodromy  functions vanish identically, i.e. for all $n, k \in \ZZ$, one has
$$
Z(-m, a, c, [X_n]^k) = Z(-m, a, c, [Y_n]^k )= \zeta(-m,a, c)
$$
(\cite[Theorem 7.2]{LL2}).
\item
For positive integers $s=m \ge 1 $  the $[Y_n]$-monodromy vanishes identically.
However  the  $[X_n]$ monodromy does not vanish for  $s= m \ge 1$ because
the functions $c_n(s)$ have simple poles
offsetting the term $(e^{2 \pi i k s} -1)$. In fact one can show that
the monodomy functions  for $m \ge 1$ are:
\begin{equation}\label{EFC1}
Z(m, a, c, [X_n]^k) = \zeta(m, a, c) +
\begin{cases}
- \frac{ k}{m!}e^{\frac{-\pi i m}{2}} (2 \pi)^{m}  \psi_n^m(a, c) &  \mbox{if} \quad m \ge 1,\\
&\\
\frac{ k}{m!}e^{\frac{\pi i m}{2}} (2 \pi)^{m}  \psi_n^m(a, c) &  \mbox{if} \quad m \le 0.
\end{cases}
\end{equation}
 The cases $s=m \ge 1$ correspond to  polylogarithm parameter values.
 \end{enumerate}
 A picture of the covering manifold structure over $\sM$ is given in Figure \ref{fg51} in Section \ref{sec32}.

%
%
%

\subsection{Homotopy generators for Lerch transcendent}\label{sec32}

We regard the Lerch transcendent variables $(s,z,c)$ as lying on the manifold
\beql{502}
 \sN := \{ (s, z, c) \in \CC \times (\PP^1 ( \CC) \smallsetminus \{ 0,1,
\infty \}) \times (\CC \smallsetminus \ZZ) \}.
 \eeq
 which is covered by  the manifold
 $\sM := \{ (s, a, c) \in \CC \times (\CC\smallsetminus \ZZ) \times (\CC \smallsetminus \ZZ)\},$
with covering map
$\pi': \sM \to \sN$ given by $\pi ' (s, a, c)= (s, e^{ 2 \pi i a}, c)$
preserving the
complex-analytic structure.
The covering group
of deck transformations is $\ZZ$. In addition $\sN$  has
universal cover $\tilde{\mathcal N}$
equal to the universal cover $\tsM$.
 Figure \ref{fg51} below exhibits the set of covering maps
relating $\tsM$, $\tsM^{ab}$, $\sM$ and $\sN$ and the associated
covering groups for each factor.

%
%
\begin{figure}[H]
\centering
\mbox{
\xymatrix%
@R+.5pc
@M=10pt
{
\tilde{\mathcal{M}} \ar[d]_{\pi_{\bf D}}^{\mathrlap{\qquad\textstyle G(\tilde{\mathcal{M}} / \tilde{\mathcal{M}}^{ab}) = {\bf D}^{(1)}(\pi_1(\mathcal{M}, {\bf x}_0))}} \\
\tilde{\mathcal{M}}^{ab} \ar[d]_{\pi_{ab}}^{\mathrlap{\qquad\textstyle G(\tilde{\mathcal{M}}^{ab}/\mathcal{M}) = \mathbb{Z}^\infty \oplus \mathbb{Z}^\infty}}\\
\mathllap{\mathbb{C} \times (\mathbb{C} \smallsetminus \mathbb{Z}) \times (\mathbb{C} \smallsetminus \mathbb{Z}) ={}}\mathcal{M} \ar[d]_{\pi'}^{\mathrlap{\qquad\textstyle G(\mathcal{M}/\mathcal{N}) = \mathbb{Z}}} \\
\mathllap{\mathbb{C} \times (\mathbb{P}^1(\mathbb{C}) \smallsetminus \{ 0, 1, \infty\} ) \times (\mathbb{C} \smallsetminus \mathbb{Z}) ={}}\mathcal{N}
}
}
\caption
{
Covering manifolds and automorphism groups.
}
\label{fg51}
\end{figure}

To obtain an analytic continuation $\tilde{Z}$ of $\Phi(s,z,c)$
to $\tsN$ we use paths starting from the  base point
\beql{206aa}
\bx'_0=(s, z, c) := (\frac{1}{2}, -1, \frac{1}{2}) \in \sN,
\eeq
 which is the image
 under the covering map $\pi'(\bx_0)= \bx'_0$
of the base point
\beql{206bb}
\bx_0 = (s, a, c) := ( \frac {1}{2} , \frac {1}{2}, \frac{1}{2} ) \in \sM
\eeq
used in part II. A path ${\gamma}$ in $\sN$ with basepoint $\bx'_0$
lifts via $\pi '$ to a unique path $\tilde{\gamma}$ in $\sM$ with
basepoint $\bx_0$, and we set
\beql{506}
\tilde{Z}([{\gamma}]) =  \tilde{Z}(s, z, c, [{\gamma}]) := Z(s, a, c, [\tilde{\gamma}]) ~.
\eeq
This agrees with
$\Phi (s,z,c)$ in a small neighborhood of $\bx'_0$, so
effects its analytic continuation to the universal cover
$\tilde{\sN} \equiv \tsM$. Theorem 4.2
of part II  (\cite{LL2}) now shows
that $\tilde{Z}$ is single-valued on $\tsM^{ab}$.

\begin{figure}[H]
\subfloat
{
\scalebox{1} 
{
\begin{pspicture}(0,-1.61)(6.41,1.61)
\psline[linewidth=0.02cm](3.4,1.6)(3.4,-1.6)
\psline[linewidth=0.02cm](0.0,0.0)(6.4,0.0)
\psdots[dotsize=0.12](3.4,0.0)
\psdots[dotsize=0.12](5.0,0.0)
\psdots[dotsize=0.12](1.8,0.0)
\psellipse[linewidth=0.04,dimen=middle](3.0,0.0)(1.2,0.8)
\psline[linewidth=0.08cm,linestyle=none]{->}(3.5,0.8)(2.9,0.8)
\psline[linewidth=0.08cm,linestyle=none]{->}(2.68,-0.8)(3.08,-0.8)
\usefont{T1}{ptm}{m}{n}
\rput(2.57,1.045){$Z_0$}
\usefont{T1}{ptm}{m}{n}
\rput(1.45,-0.235){$-1$}
\usefont{T1}{ptm}{m}{n}
\rput(3.17,-0.255){$0$}
\usefont{T1}{ptm}{m}{n}
\rput(5.01,-0.275){$1$}
\end{pspicture}
}
}\qquad
\subfloat
{
\scalebox{1} 
{
\begin{pspicture}(0,-1.61)(6.41,1.61)
\psline[linewidth=0.02cm](3.4,1.6)(3.4,-1.6)
\psline[linewidth=0.02cm](0.0,0.0)(6.4,0.0)
\psdots[dotsize=0.12](3.4,0.0)
\psdots[dotsize=0.12](5.0,0.0)
\psdots[dotsize=0.12](1.8,0.0)
\usefont{T1}{ptm}{m}{n}
\rput(5.39,0.965){$Z_1$}
\usefont{T1}{ptm}{m}{n}
\rput(1.45,-0.235){$-1$}
\usefont{T1}{ptm}{m}{n}
\rput(3.17,-0.255){$0$}
\usefont{T1}{ptm}{m}{n}
\rput(5.01,-0.275){$1$}
\psbezier[linewidth=0.04](1.8,0.0)(2.2,0.6)(2.8,1.1)(3.58,1.16)(4.36,1.22)(4.98,0.72)(5.36,0.36)(5.74,0.0)(5.52,-0.46)(5.1,-0.62)(4.68,-0.78)(4.42,-0.62)(4.22,-0.32)(4.02,-0.02)(3.96,0.42)(3.56,0.56)(3.16,0.7)(2.4,0.4)(1.8,0.0)
\psline[linewidth=0.08cm,linestyle=none]{->}(3.88,0.32)(3.96,0.18)
\psline[linewidth=0.08cm,linestyle=none]{->}(4.76,0.84)(4.64,0.9)
\end{pspicture}
}
}
\caption{Generators $[Z_0], [Z_1]$ of $\pi_1 (\PP^1 \smallsetminus \{0,1,\infty\}, -1)$ in the $z$-plane}
\label{fg52}
\end{figure}
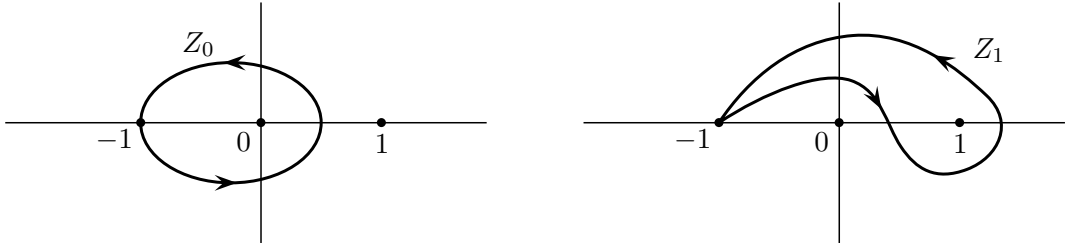

We now specify a set of generators $\sG'$ of  $\pi_1(\sN, \bx'_0)$.
The group $\pi_1 (\sN , \bx'_0)$ has $\pi_1 (\sM , \bx_0)$ as a normal
subgroup with quotient group $\ZZ$. In particular it contains more
closed loops than $\pi_1 ( \sM , \bx_0)$, e.g. a path with
endpoints $(s, a,c)$ and $(s, a+1,c)$ in $\sM$ projects to a
closed path on $\sN$. Since $\sN$ is a product manifold, we have
\beql{507}
 \pi_1 (\sN , \bx'_0 ) \simeq \pi_1 ( \PP^1 \smallsetminus \{ 0,1,\infty \} , -1) \times
\pi_1 ( \CC \smallsetminus \ZZ , \frac{1}{2} ) ~.
\eeq
Now $\pi_1 ( \PP^1 \smallsetminus \{0,1, \infty \}, -1 )$ is a
free group on two generators $[Z_0]$ and $[Z_1]$.
Our choice of loops $Z_0$ and $Z_1$ is  pictured in Figure \ref{fg52}.

%
%

\noindent Here we choose  $Z_0$
to be a loop in the $z$-plane based at $z=-1$ which
is a simple closed curve enclosing $z=0$ counterclockwise
and not intersecting the line $L_1^{-}= \{ 1-it: t \ge 0\} $
(so  not enclosing $z=1$.)
We choose {$Z_1$} to be a simple closed curve enclosing $z=1$ counterclockwise, and not
intersecting the half-line $L_0^{-} = \{ -it: t \ge 0 \}$.
 The loops $Z_0$ and
$Z_1$ may be viewed as closed paths in $\sN$ which hold the variables
$s= \frac{1}{2}$ and $c= \frac{1}{2}$ constant.
One may
also  define a homotopy class $[Z_{\infty}]$ so that $[Z_0][Z_1] [Z_{\infty}]= [\mbox{Id}]$
 is trivial. A representative of this class is  a clockwise oriented
loop enclosing both $z=0$ and $z=1$.

For  the $c$-variable in
{$\CC \smallsetminus \ZZ$}
we retain  the generators $[Y_n]$ given in Section \ref{sec31},
and we obtain the full  set of generators
\beql{509}
\sG' := \{ [Z_0 ], [Z_1] \} \cup \{ [Y_n] : n \in \ZZ\}
\eeq
for $\pi_1(\sN, \bx'_0)$.

We next determine the lifts of the loops $[Z_0]$
and $[Z_1]$ to the $a$-plane,
and gives other  relevant information on the relation of $\sM$ and $\sN$,
as follows.


\begin{lemma}\label{le30}

(1)
The lift $\widetilde{Z_0}$  to $(\sM, \bx_0)$
 of the loop $Z_0$ on  $(\sN, \bx'_0)$
 is a non-closed path with $c= s = \frac{1}{2}$ fixed and whose projection to the $a$-plane connects the base point $a=\frac{1}{2}$ to $a= \frac{3}{2}$
while remaining in the upper half-plane $\{ a: \Im(a) >0\}$.

(2) The lift $\widetilde{Z_1}$ to $(\sM, \bx_0)$ of $Z_1$ on $(\sN, \bx_0')$ has $c= s = \frac{1}{2}$ fixed
and is homotopic to the closed loop  $ X_0$ on $(\sM, \bx_0)$.

 (3) The image group $\hH_1:= (\pi')_{\ast} ( \pi_1(\sM, \bx_0))$ in  $\pi_1(\sN, \bx'_0)$
 is given by
\beql{519b}
\hH_1 := \langle [Z_0]^{-n} [Z_1] [Z_0]^n, ~[Y_n], ~~n \in \ZZ \rangle.
\eeq

 (4) The image group $\hH_1$ contains the commutator subgroup of $\pi_1(\sN, \bx'_0)$, i.e.
 \beql{519c}
(\pi_1(\sN, \bx'_0))' :=  [ \pi_1(\sN, \bx'_0), \pi_1(\sN, \bx'_0)]  \subset \hH_1~.
\eeq
In particular, $\hH_1$ is a normal subgroup of
 $\pi_1(\sN, \bx'_0)$, and $\sM$ is an abelian Galois covering of $\sN$.
 \end{lemma}

\begin{figure}[H]
\centering
\subfloat
{
\scalebox{1} 
{
\begin{pspicture}(0,-1.61)(5.6,1.61)
\psline[linewidth=0.02cm](1.2,1.6)(1.2,-1.6)
\psline[linewidth=0.02cm](0.0,0.0)(4.8,0.0)
\psdots[dotsize=0.12](1.2,0.0)
\psdots[dotsize=0.12](3.4,0.0)
\usefont{T1}{ptm}{m}{n}
\rput(3.85,0.685){$\tilde{Z}_0$}
\usefont{T1}{ptm}{m}{n}
\rput(0.97,-0.255){$0$}
\usefont{T1}{ptm}{m}{n}
\rput(3.41,-0.295){$1$}
\psdots[dotsize=0.12](2.3,0.0)
\usefont{T1}{ptm}{m}{n}
\rput(2.29,-0.375){$a = \frac{1}{2}$}
\usefont{T1}{ptm}{m}{n}
\rput(4.51,-0.375){$\frac{3}{2}$}
\psline[linewidth=0.02cm](4.8,0.0)(5.3,0.0)
\psdots[dotsize=0.12](4.5,0.0)
\psline[linewidth=0.06cm](2.3,0.0)(3.0,0.0)
\psline[linewidth=0.06cm](3.8,0.0)(4.5,0.0)
\psarc[linewidth=0.06](3.4,0.0){0.4}{355.10092}{184.64508}
\psline[linewidth=0.08cm,linestyle=none]{->}(2.5,0.0)(2.8,0.0)
\psline[linewidth=0.08cm,linestyle=none]{->}(4.0,0.0)(4.3,0.0)
\psline[linewidth=0.08cm,linestyle=none]{->}(3.12,0.36)(3.5,0.42)
\end{pspicture}
}
}\qquad
\subfloat
{
\scalebox{1} 
{
\begin{pspicture}(0,-1.81)(5.02,1.81)
\psline[linewidth=0.02cm](2.88,1.8)(2.88,-1.8)
\psline[linewidth=0.02cm](0.88,0.0)(4.88,0.0)
\psdots[dotsize=0.12](2.88,0.0)
\psdots[dotsize=0.12](4.48,0.0)
\usefont{T1}{ptm}{m}{n}
\rput(4.31,1.365){$Z_0$}
\usefont{T1}{ptm}{m}{n}
\rput(2.67,-0.275){$0$}
\usefont{T1}{ptm}{m}{n}
\rput(4.71,0.205){$1$}
\psdots[dotsize=0.12](1.28,0.0)
\psarc[linewidth=0.04](2.88,0.0){1.6}{14.036243}{345.96375}
\psarc[linewidth=0.04](4.48,0.0){0.4}{93.57633}{266.86365}
\usefont{T1}{ptm}{m}{n}
\rput(0.51,-0.255){$z = -1$}
\psline[linewidth=0.08cm,linestyle=none]{->}(2.0,1.36)(1.66,1.04)
\psline[linewidth=0.08cm,linestyle=none]{->}(3.56,-1.58)(4.02,-1.14)
\psline[linewidth=0.08cm,linestyle=none]{->}(3.96,-0.08)(4.2,0.3)
\end{pspicture}
}
}
\caption{  (a) Loop  $\tilde{Z}_0$ in $a$-plane,  (b) Projection of $\tilde{Z}_0$ in $z$-plane.}
\label{fg54}
\end{figure}
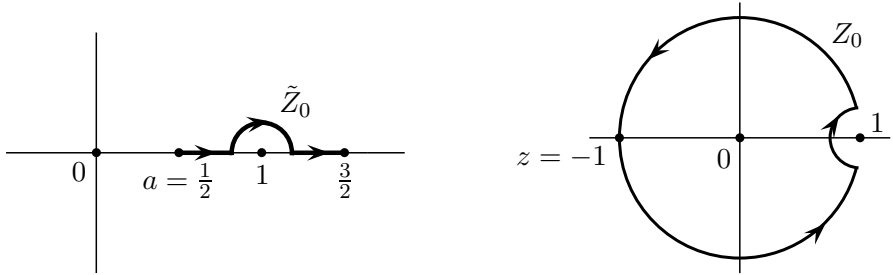

 \paragraph{\em Proof.}
 (1).  This is established by  checking that images under $\pi'$ of certain paths
 in the $a$-plane result in loops homotopic to $[Z_0]$ (resp. $[Z_1]$) in the
 $z$-plane. For $\widetilde{Z}_0$ we consider a path from
$a= \frac{1}{2}$ to the point $a= \frac{3}{2}$ that consists of
 a line segment in $a$-plane, from $a=\frac{1}{2}$
 to $a= \frac{3}{2}$, except for a small clockwise-oriented half-circle in
$Im(a) >0$ centered at $a=1$, made to detour around the point $a=1$.
{When projected to the $z$-plane,
the curve}
proceeds from $z=-1$ in a counterclockwise
circle of radius $1$ around $z=0$, with an indentation near $z=1$ to
leave $z=1$ outside the loop. This image is clearly homotopic to $Z_0$.
These are pictured in Figure \ref{fg54}.

%
%

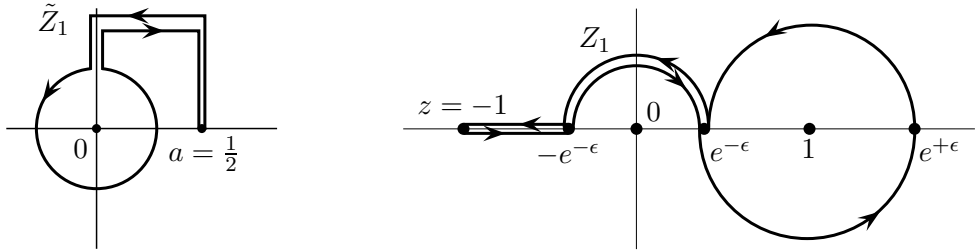
\begin{figure}[H]
\centering
\subfloat
{
\scalebox{1} 
{
\begin{pspicture}(0,-1.61)(4.18,1.61)
\psline[linewidth=0.02cm](1.38,1.6)(1.38,-1.6)
\psline[linewidth=0.02cm](0.18,0.0)(3.78,0.0)
\psdots[dotsize=0.12](1.38,0.0)
\psdots[dotsize=0.12](2.78,0.0)
\usefont{T1}{ptm}{m}{n}
\rput(0.81, 1.45){$\tilde{Z}_1$}
\usefont{T1}{ptm}{m}{n}
\rput(1.17,-0.275){$0$}
\usefont{T1}{ptm}{m}{n}
\rput(2.81,-0.335){$a = \frac{1}{2}$}
\psarc[linewidth=0.04](1.38,0.0){0.8}{95.710594}{84.55967}
\psline[linewidth=0.04](1.46,0.78)(1.46,1.3)(2.74,1.3)(2.74,0.0)
\psline[linewidth=0.04](1.3,0.78)(1.3,1.5)(2.82,1.5)(2.82,-0.02)
\psline[linewidth=0.08cm,linestyle=none]{->}(0.96,0.8)(0.66,0.36)
\psline[linewidth=0.08cm,linestyle=none]{->}(1.88,1.3)(2.28,1.3)
\psline[linewidth=0.08cm,linestyle=none]{->}(2.28,1.5)(1.78,1.5)
\end{pspicture}
}
}\qquad
\subfloat
{

\scalebox{1} 
{
\begin{pspicture}(0,-1.605)(8.36,1.605)
\psline[linewidth=0.01cm](3.1,1.6)(3.1,-1.6)
\psline[linewidth=0.01cm](0.0,0.0)(7.6,0.0)
\psdots[dotsize=0.16](3.1,0.0)
\usefont{T1}{ptm}{m}{n}
\rput(2.55,1.185){$Z_1$}
\usefont{T1}{ptm}{m}{n}
\rput(3.33,0.225){$0$}
\psdots[dotsize=0.16](5.4,0.0)
\psdots[dotsize=0.16](0.8,0.0)
\psdots[dotsize=0.16](4.0,0.0)
\psdots[dotsize=0.16](2.2,0.0)
\psline[linewidth=0.04cm](0.8,0.06)(2.2,0.06)
\psline[linewidth=0.04cm](0.8,-0.06)(2.2,-0.06)
\psarc[linewidth=0.04](3.1,0.0){0.84}{0.0}{180.0}
\psarc[linewidth=0.04](3.1,0.02){0.96}{0.0}{180.0}
\usefont{T1}{ptm}{m}{n}
\rput(5.39,-0.255){$1$}
\usefont{T1}{ptm}{m}{n}
\rput(4.35,-0.275){$e^{-\epsilon}$}
\usefont{T1}{ptm}{m}{n}
\rput(2.19,-0.315){$-e^{-\epsilon}$}
\usefont{T1}{ptm}{m}{n}
\rput(0.79,0.305){$z=-1$}
\psline[linewidth=0.08cm,linestyle=none]{->}(1.12,-0.06)(1.36,-0.06)
\psline[linewidth=0.08cm,linestyle=none]{->}(1.7,0.06)(1.54,0.06)
\psline[linewidth=0.08cm,linestyle=none]{->}(3.48,0.9)(3.36,0.94)
\psline[linewidth=0.08cm,linestyle=none]{->}(3.6,0.68)(3.76,0.54)
\psdots[dotsize=0.16](6.8,0.0)
\psarc[linewidth=0.04](5.37,-0.03){1.43}{179.22227}{1.7357045}
\psarc[linewidth=0.04](5.43,0.01){1.37}{0.0}{180.0}
\psline[linewidth=0.08cm,linestyle=none]{->}(5.92,-1.4)(6.34,-1.1)
\psline[linewidth=0.08cm,linestyle=none]{->}(5.18,1.44)(4.8,1.22)
\usefont{T1}{ptm}{m}{n}
\rput(7.13,-0.275){$e^{+\epsilon}$}
\end{pspicture}
}
}
\caption{    (a) Loop  $\tilde{Z}_1$ in $a$-plane,  (b) Projection of $\tilde{Z}_1$ in $z$-plane.}
\label{fg55}
\end{figure}

(2). We assert that the homotopy class
{of the projection of $\widetilde{Z}_1$} to the
$z$-plane equals
{that of $X_0$}, see Figure \ref{fg55}.
To verify this, note that the  path $X_0$  given in part II  
is  a closed path in the $a$-plane that first
moves vertically from the base point
$a=\frac{1}{2}$ to $a= \frac{1}{2} + \frac{i \epsilon}{2 \pi},$ (for small enough $\epsilon$)
then moves horizontally to $a= \frac{i \epsilon}{2\pi} $, then moves in a counterclockwise loop
of radius $\frac{\epsilon}{2\pi}$ around $a=0$ back to $a= \frac{i \epsilon}{2 \pi}$, and finally
returns to $a= \frac{1}{2}$ following the original path.
One may verify
that,
{when projected to the $z$-plane}, the image of ${X_0}$
moves along the
$z$-axis from $z= -1$ to $z= -e^{ -\epsilon}$, then proceeds in a clockwise
half-circle to $z= e^{- \epsilon}$.
Next,  the image of the counterclockwise loop in the $a$-plane around $a=0$ is
$$z = \exp( -\epsilon e^{i \theta}), 0 \le \theta \le 2 \pi,$$
which  is a nearly circular path that
encircles $z=1$ counterclockwise, reaching at $\theta=\pi$ the point $z=e^{\epsilon} >1$
on the real axis, with the second half of its path from $\theta= \pi$ to
$\theta=2\pi$ being the reflection of the first half in the real axis. Then it returns
to $z= -1$ along the outgoing path. This path is clearly
homotopic to $Z_1$, whence $\widetilde{[Z_1]} = [X_0]$.
This is pictured in Figure \ref{fg55}.

%
%

(3).  Recall that $\pi_1 ( \sM , \bx_0)$
 has generating set
\beql{519aa}
\sG := \{ [X_n]: n \in \ZZ\} \cup \{ [Y_n]: n \in \ZZ \},
\eeq
in which the homotopy class  $[X_n]$ is given by a path $X_n$  in the $a$-plane with
basepoint $a= \frac{1}{2}$, holding $s= \frac{1}{2}, c=\frac{1}{2}$ fixed throughout,
that traverses a line segment in the upper-half plane
to $a= n + \epsilon i,$ followed by a counterclockwise loop of radius $\epsilon$
around $a=n$, followed by return along the line segment; similarly
 $[Y_n]$ is given by a path $Y_n$  in the $c$-plane with basepoint $c= \frac{1}{2}$,
 holding $s=\frac{1}{2}$ and $a=\frac{1}{2}$ fixed throughout,
that traverses a line segment in the upper-half plane
to $c= n + \epsilon i,$ followed by a counterclockwise loop of radius $\epsilon$
around $c=n$, followed by return along the line segment.
Extending the argument of (2) we find that
\beql{519}
(\pi ' )_{\ast}([X_{n}]) = [Z_0 ]^{n} [Z_1] [Z_0]^{-n} ~, ~~n\in \ZZ.
\eeq
This follows since the path $[Z_0]^n [Z_1] $ lifted to $\sM$ first moves from $a= \frac{1}{2}$ to
the point $a= n+ \frac{1}{2}$, then encircles
the point $a=n$ once counterclockwise and returns to $a= n+ \frac{1}{2}$
and $[Z_0]^{-n}$ then returns to the base point
$\bx_0=(\frac{1}{2}, \frac{1}{2}, \frac{1}{2} )$ in $\sM$.
We also trivially have
$$
(\pi ' )_{\ast}([Y_{n}])= [Y_{n}], ~~~n\in \ZZ,
$$
since the projection is constant.
It follows that  the image group $\hH_1 := (\pi')_{\ast}( \pi_1(\sM, \bx_0))$ has
generating set  \eqn{519b}.

(4).
To verify the inclusion \eqn{519c}, note first
that since  both generators $[Z_0]$ and $[Z_1]$ commute with all $[Y_n]$
in $\pi_1(\sN, \bx'_0)$ we have
\beql{519d}
(\pi_1(\sN, \bx'_0))'=\langle [Z_0]^{k} [Z_1]^l [Z_0]^{-k} [Z_1]^{-l},
[Y_n]^k [Y_p]^l [Y_n]^{-k}[Y_p]^{-l},
k,l, n, p \in \ZZ \rangle.
\eeq
By (2), $[Z_1] \in \hH_1$, whence all commutators $[Z_0]^{k} [Z_1]^l [Z_0]^{-k} [Z_1]^{-l}$
are in $\hH_1$, and since all $[Y_n] \in \hH_1$ we see that all generators of
$(\pi_1(\sN, \bx'_0))'$ are in $\hH_1$ and the inclusion \eqn{519c} follows.

Finally, since all subgroups of a group that contain its
commutator subgroup are normal, we conclude  from \eqn{519c} that $\hH_1$ is a normal subgroup of
$\pi_1(\sN, \bx'_0)^{'}$. Thus $\sM$ is an abelian covering of $\sN$. $~~~\Box$

%
%
%

\subsection{Multivalued  continuation of  Lerch transcendent}\label{sec33}

The Lerch transcendent  is defined by \eqn{101} as an analytic function of
three variables on the region
$$
 \sDD_0  = \{ s: \Re(s) > 1\} \times \{ z: |z| < 1 \} \times \{c: 0 < \Re(c) < 1\}.
$$
on which it satisfies $\Phi(s, e^{2 \pi i a}, c) = \zeta(s, a, c)$.
We now show it analytically continues to a single-valued function on the
 simply-connected region
 \beql{WL510}
  \sDD := \{ s: s \in \CC\} \times \{ z: z \in \CC \smallsetminus \RR_{\ge 0} \} \times \{c : 0 < \Re(c) < 1 \} \subset \sN
 \eeq
 of $\sN$.
 The region $\sDD$  is the image under
 the covering map $\pi'$ of the  {\em extended fundamental polycylinder}
\beql{WL510b}
\tilde{\Omega}:= \{ s : s \in \CC \} \times \{a :  0 < \Re (a) < 1 \} \times
\{c:  0 < \Re (c) < 1 \}  \subset \sM ~
\eeq
and the restriction
$\pi': \tilde{\Omega} \to \sDD$  is a biholomorphic map.
 Theorem 2.1 of part II analytically continued  $\zeta(s, a,c)$ to $\tilde{\Omega}$,
 which gives the analytic continuation of $\Phi(s, z, c)$ to $\sDD$, defined via
\beql{510c}
\Phi(s, e^{2 \pi i a}, c) =  \zeta(s, a, c)  ~~\mbox{for}~~ (s, a, c) \in \tilde{\Omega}.
\eeq
We regard $\sDD$ as embedded in the universal cover
$\tilde{\sN} \equiv \tilde{\sM}$ by lifting it to $\tilde{\Omega}$ in $\sM$
 followed by the embedding of
$\tilde{\Omega}$ in $\tilde{\sM}$.\\

To describe the multivalued nature of the analytic continuation of
the Lerch transcendent $\tilde{Z}$, we recall some definitions from part II.
%
\begin{defi}\label{Nde21}
{\rm Let $f: \tsN \to \CC$ be a continuous function on the
universal cover $\tsN$ of a manifold $\sN$ and let $[\tau ] \in\pi_1 ( \sN , \bx'_0 )$
 be a homotopy class. The operator $\sQ_{[\tau]} $ takes the function $f$ to the
  {\em $[\tau]$-translated function} $\sQ_{[\tau]} (f): \tsN \to \CC$
defined by
 \beql{Z204}
\sQ_{[\tau]} (f) ( [\gamma ] ) := f([\tau \gamma  ]) \,,
\eeq
where
$\gamma$ is an arbitrary path with basepoint $\gamma (0) = \bx'_0$ and
$ \tau \gamma$ is the composed path (first follow $\tau$ and then  $\gamma$,
as in \cite[p. 26]{Hat01}).
}
\end{defi}
%
\begin{defi}\label{Nde22}
{\rm
Let $f: \tsN \to \CC$ be a continuous function and let
$[\tau ] \in \pi_1 (\sN, \bx'_0)$ be a homotopy class.
The {\em monodromy function} $M_{[\tau]} (f): \tsN \to \CC$
of $f$ at $[\tau ]$ is defined by
\beql{Z205A}
M_{[\tau]} (f) := ( \sQ_{[\tau]} - I ) (f) \, .
\eeq
That is, for all paths $\gamma : [0,1] \to \sN$ with base point
$\gamma (0) = \bx'_0$,
\beql{Z206}
M_{[\tau]} (f) ([\gamma]) :=
f([\tau \gamma  ]) - f ([\gamma ]) ~.
\eeq
}
\end{defi}

In the sequel  we will need two different branches of the logarithm, defined
as follows.
We let $\log z$ denote the principal branch of the logarithm, cut along
the negative real axis, with the negative real axis itself viewed as belonging
to the upper half-plane, so $\log 1=0, \log (-1) = \pi i, \log(-i) = - \frac{\pi i}{2}$.
The  {\em semi-principal branch} $\llog z$ of the logarithm is
defined on the complex plane cut
along the positive real axis, whose value at $z=-1$ is $\pi i$,
and with the positive real axis connected to the upper half-plane, so $\llog(1)=0$,
$\llog (-1) = \pi i, \llog (-i) = \frac{3 \pi i}{2}$.

We now  specify certain functions
that  will appear as monodromy functions of the Lerch transcendent.
For each integer $n$, we now define the function
\beql{WL511}
f_n(s,z,c) := \left\{
\begin{array}{rll}
e^{\pi i (s - 1)} e^{2 \pi inc} z^{-c}~( n - \frac{1}{2 \pi i} \llog z)^{s-1}
 ~ & {\rm if} & n \geq 1, \\
 &  & \\
 e^{2 \pi inc} z^{-c}~(\frac{1}{2 \pi i} \llog z - n)^{s-1}  ~~~~~~~
& {\rm if} & n \leq 0, \\
\end{array}
\right.
\eeq
on the simply-connected domain $\sDD$.
Here   $z^{-c} = e^{- c \llog z}$, and $a= \frac{1}{2 \pi i} \llog z$,
but  in \eqn{WL511} we evaluate the term
$(a-n)^{s-1}:= e^{(s-1) \log (a-n)}$ using the principal branch
of the logarithm, noting that $a-n$ (resp. $n-a$) always has
positive real part when $(s, a, c) \in \tilde{\Omega}$.
The  semi-principal branch $\llog z$ in the formula above
 is  needed to apply on  the domain $z \in \CC \smallsetminus \RR_{\ge 0}$ used
 in $\sDD$.
The function \eqn{WL511} then extends to a function $f_n ([\gamma])$ on $\tilde{\sN}$ by analytic
continuation.


\begin{theorem}\label{th51}
{\em (Lerch Transcendent Monodromy Formulas)}
The Lerch transcendent $\Phi (s,z,c)$ analytically
continues in three complex variables
to a single-valued function $\tilde{Z}= \tilde{Z}(s,z,c, [\gamma]) $ on the universal
cover
 $\tilde{\sN}$ of
 $$
 \sN =\{ (s,z,c) \in \CC \times ( \PP^1 ( \CC )\smallsetminus \{ 0,1, \infty \} )
  \times ( \CC \smallsetminus \ZZ)\}.
  $$
  The monodromy functions of $\tilde{Z}$ for the generating set $\sG'$  of
$\pi_1 ( \sN , \bx'_0 )$ are as follows.

\begin{itemize}
\item[(i-a)]
In the domain  $(s,z,c) \in \sDD \subset \tsN$ given in \eqn{WL510} the
monodromy function for the generators
$[Z_0], [Z_1]$ are
\beql{WL512}
M_{[Z_0]} ( \tilde{Z}) (s, z, c) \equiv 0 ~,
\eeq
and
\beql{WL513}
M_{[Z_1]} ( \tilde{Z}) (s,z,c) = -
\frac{(2 \pi )^s}{\Gamma (s)} e^{\frac{i \pi s}{2}} f_0 (s,z,c) ~.
\eeq
\item[(i-b)] On the domain $\sDD$
the functions $f_p (s,z,c)$ for a fixed $p \in \ZZ$ have monodromy functions
\beql{WL514}
M_{[Z_0]^k} (f_p)(s, z, c) = f_{p-k}(s, z, c) - f_p(s, z,c)~~~~for~~~k \in ~\ZZ,
\eeq
and
\beql{WL515}
M_{[Z_1]^k} (f_p)(s,z,c) \equiv 0 ~~~~for~~~k \in ~\ZZ.
\eeq
In addition
\beql{WL515b}
M_{[Y_n]^k} (f_p)(s,z,c) \equiv  0 ~~~~for~~~k, n \in ~\ZZ.
\eeq

\item[(ii-a)]
On the domain $\sDD$, the monodromy functions for
the generators $[Y_n]$ for all $n \in \ZZ$ are
\beql{WL515a}
M_{[Y_n]} (\tilde{Z}) (s,z,c) =\left\{ \begin{array}{lll}
0 & \mbox{if} & n \ge 1 ~, \\
(e^{- 2 \pi is} -1) z^{-n} (c-n)^{-s} &
\mbox{if} & n \le 0 ~.
\end{array}\right.
\eeq
\item[(ii-b)]
For a path  $\gamma$ in $\sN$ from $\bx'_0$
to an endpoint  falling in the multiply-connected region
$$
\sN_s= \{s\} \times ( \PP^1 (\CC ) \smallsetminus \{ 0,1,
\infty \} ) \times ( \CC \smallsetminus \ZZ ),
$$
there holds
\beql{WL516}
M_{[Y_n]^{-1}} (\tilde{Z} ) ([\gamma])= - e^{2 \pi is}
M_{[Y_n]} (\tilde{Z})([\gamma])
 \eeq
 and
 \beql{WL517} M_{[Y_n]^{\pm k}} (\tilde{Z} )([\gamma])
 =\frac{e^{\mp 2 \pi isk} -1}{e^{\mp 2 \pi is} -1}
M_{[Y_n]^{\pm 1}} (\tilde{Z})([\gamma]) ~.
 \eeq
\end{itemize}
\end{theorem}

\paragraph{\bf Remarks.}
(1) The hypothesis  ``on the domain $\sDD$''  above
means to view $(s, z, c, [\gamma]) \in \sDD$  as the endpoint of
a path $\gamma$ based at $\bx'_0=(\frac{1}{2}, -1, \frac{1}{2})$ that
remains entirely inside $\sDD$.

(2) The formulas given in (i-a), (i-b) above are sufficient to compute the monodromy
functions for all elements of the subgroup
$\hG_Z := \langle [Z_0], [Z_1] \rangle  $
of $\pi_1 ( \sN , \bx'_0)$; these include all monodromy functions not already
determined in part II.

\paragraph{\em Proof.}
(i-a). We use the fact shown earlier
that $(\pi')_{\ast}([X_m]) = [Z_0]^{-m} [Z_1][Z_0]^m$ for $m \in \ZZ$,
which gives
$$
M_{[X_m]} (Z) = M_{ [Z_0]^{-m} [Z_1][Z_0]^m}(\tilde{Z}).
$$
We can now apply the formula
\beql{519ab}
M_{[X_m]} (Z) (s,a,c) = \left\{
\begin{array}{rll}
-\frac{(2 \pi)^s e^{\frac{\pi i s}{2}}}{\Gamma (s)}
e^{\pi i (s - 1)}(m - a)^{s-1} e^{- 2 \pi ic (a-m)} ~ & {\rm if} & m \geq 1, \\
 & & \\
-\frac{(2 \pi)^s e^{\frac{\pi i s}{2}}}{\Gamma (s)} (a-m)^{s-1}
e^{- 2 \pi ic (a-m)} ~~~~~  & {\rm if} & m \leq 0, \\
\end{array}
\right.
 \eeq
 of Theorem 4.1
of part II. This
immediately gives, for $(s, z, c) \in \sDD,$ on choosing $m=0$, that
\beql{520}
M_{[Z_1]} ( \tilde{Z} )(s, z, c) =  -\frac{(2 \pi )^s   e^{\frac{\pi is}{2}}       }{\Gamma (s)}
z^{-c} \left(\frac{1}{2 \pi i}\llog z \right)^{s-1} ~,
\eeq
 in which $z^{-c} = e^{-c \llog z}$.

The  Lerch zeta function
$\zeta (s,a,c) = \sum_{n=0}^\infty e^{2 \pi ina} (n+c)^{-s}$
is defined by the right side as an analytic function  of three variables in the region
$$
\sU^+ := \{s: \Re (s) > 0 \} \times \{a: \Im (a) > 0 \} \times \{c: \Re (c) > 0 \} ~
$$
and on this region it is invariant under  $a \mapsto a+1$.
Recall from part II \cite[Theorem 4.1 proof]{LL2} that this function
analytically continues  to a single-valued function in the  larger domain
$$
\sO_1 := \{ s: \Re (s) > 0 \} \times
\{ a: a \in \sA_L \} \times \{ c: \Re (c) > 0 \} ~,
$$
in which $\sA_L= \CC \setminus \{ M_k: k \in \ZZ\}$
with $M_k = \{ a=k-it: t \ge 0\}$.
We can find a lifting of the path in the class $[Z_0]$  in
$\sM$ that is a path in $\sO_1$ from
$\bx_0 = (\frac{1}{2}, \frac{1}{2}, \frac{1}{2} )$ to
$\bx_1 = ( \frac{1}{2} , \frac{3}{2} , \frac{1}{2} )$ in $\sM$,
which has $\Im (a) > 0$
everywhere except at the endpoints, so lies in $\sU^{+}$ except at the endpoints.
Therefore the equality $\zeta (s,a,c) = \zeta (s,a+1, c)$ holds
near both endpoints $\bx_0$ and $\bx_1$ of this path
(in the upper half plane $\Im (a) > 0$)
and we obtain
\beql{527}
M_{[Z_0]} ( \tilde{Z} ) \equiv 0 ~,
\eeq
since it is identically zero in a small open disk containing $\bx_0'.$

(i-b) The calculation of the monodromy functions
of $f_p (s,z,c)$ for the loops
$[Z_0]^k$, $[Z_1]^k$ and $[Y_n]^k$ is straightforward, except that the case of
$[Z_0]^k$ where $p+k$ and  $p$  have opposite signs requires some care.

(ii-a), (ii-b)
Clearly $(\pi')_{\ast} ([Y_n]) = [Y_n]$,
the closed loops $[Y_n]^k$ in $\sM$ based at $\bx_0$ project to closed loops $[Y_n]^k$
based at $\bx'_0$ in $\sN$.
The monodromy formulas in Theorem 4.1(ii) of part II immediately apply
to yield \eqn{WL515a}--\eqn{WL517}.~~~$\bsq$\\

%
%
%

\subsection{Conditions for vanishing monodromy}\label{sec34}

We study the action of
the subgroups
\beql{527bb}
\hG_Z := \langle [Z_0] , [Z_1] \rangle
\quad\mbox{and}\quad
\hG_Y
:= \langle [Y_n ]: n \in \ZZ \rangle
\eeq
of $\pi(\sN, \bx'_0)$ acting on monodromy functions. We show
the actions commute, and further analyze their actions to
show that certain monodromy functions  vanish identically.


\begin{theorem}\label{th52}
{\em (Lerch Transcendent Vanishing Monodromy)}
Suppose that $[\tau ] \in \pi_1 ( \sN , \bx'_0)$ satisfies
\beql{528}
[\tau ] := [S_1]^{\ep_1} \cdots [S_m]^{\ep_m}
\eeq
with each $[S_i] \in \{ [Z_0 ], [Z_1] \} \cup \{ [Y_n] : n \in \ZZ \}$,
and with each $\ep_j$ equal to one of $\{ \pm 1 \}$.
Define $[\tau_Z]$  (resp. $[\tau_Y]$) to be the element of
$\pi_1 (\sN , \bx'_0)$ obtained by setting all
generators $[Y_n]$ (resp. $[Z_0]$ and $[Z_1]$) in \eqn{528} to be
the identity,
so that
$$[\tau_Z ] \in \hG_Z = \langle [Z_0 ], [Z_1] \rangle \quad\mbox{and}\quad
[\tau_Y] \in \hG_Y :=\langle [Y_n ] : n \in \ZZ \rangle.
$$
Then the following hold.

(1) For all paths $\gamma$ based at $\bx'_0$,
the monodromy function $M_{[\tau]} (\tilde{Z} )([\gamma]) $
satisfies
\beql{529}
M_{[\tau]} ( \tilde{Z} ) ([\gamma])= M_{[\tau_Z ]} (\tilde{Z} ) ([\gamma])+
M_{[\tau_Y]} (\tilde{Z} ) ([\gamma])~.
\eeq

(2) For all paths $\gamma$ based at $\bx'_0$,
the monodromy function $M_{[\tau_Y]} (\tilde{Z} ) ([\gamma])$ satisfies
\beql{530}
M_{[\tau_Y]} (\tilde{Z} ) ([\gamma])=
\sum_{n\in \ZZ} M_{[Y_n]^{k(n)}  }  (\tilde{Z} ) ([\gamma])~,
\eeq
in which $k(n) \in \ZZ$ is the sum of exponents of $[Y_n]$ occurring in
$[\tau]$.

(3) The monodromy function
$M_{[\tau]} ( \tilde{Z} )$ vanishes identically for any
$[\tau ]$ in the subgroup
$\Gamma'$ of $\pi_1(\sN, \bx'_0)$ defined by
\beql{531}
\Gamma^{'} := [ \hH_0: \hH_0],
\eeq
in which
\beql{532}
\hH_0 := \langle [Z_0]^{m} [Z_1] [Z_0]^{-m} , [Z_0]^{m}[Y_n][Z_0]^{-m} :~ m,~ n \in \ZZ \rangle.
\eeq
The group $\Gamma'$ is a normal subgroup
of $(\pi_1(\sN, \bx'_0))$ that contains its
  second commutator subgroup
$(\pi_1(\sN, \bx'_0))^{'' }$. In particular its quotient
$ \pi_1(\sN,\bx'_0)/ \Gamma'$ is a two-step solvable group.
\end{theorem}

\paragraph{\em Proof.}
In this proof the argument is carried out for each homotopy class  $[\gamma]$
separately.  We regard it as  fixed, and
so  abbreviate $M_{[S]}(f) [\gamma])$ to $M_{[S]}(f)$, throughout.

We know that the generators $[Z_0]$ and $[Z_1]$ of $\pi_1(\sM, \bx'_0)$
each commute with all generators $[Y_n]$, but this  is not  in itself sufficient to
imply \eqn{529}.
The extra facts to be used are that the $[Z_j]$-generators annihilate
all $[Y_n]$-monodromy functions,
and the $[Y_n]$-generators annihilate all $[Z_j]$-monodromy functions.
More precisely, since
$z^{-n} (c-n)^{-s}$ is single-valued in the
$z$-variable
for fixed $s$ and $c$ (for $n \le 0$), it follows from Theorem \ref{th51}
that
$$
M_{[Z_j]^{\pm 1}} (M_{[Y_n]^{\pm k}}(\tZ) ) \equiv 0
$$
for $j \in \{0,1\}$ and all  $k, n \in \ZZ$. Also
since $f_p (s,z,c)$ is single-valued in the
$c$-variable for fixed $s$ and $z$,
it follows from Theorem \ref{th51} that
$$
M_{[Y_n]^{\pm 1}} (f_p) \equiv0
$$
for all  $p \in \ZZ$.
These relations give for all
$[\tau_Y ] \in \hG_Y$ and $[\tau_Z ] \in \hG_Z$ that
\beql{533}
M_{[Z_j]^{\pm 1}} (M_{[\tau_Y]} (\tilde{Z} )) \equiv 0
\quad\mbox{for}\quad
j \in \{0,1 \} ~,
\eeq
\beql{534}
M_{[Y_n]^{\pm 1}} (M_{[\tau_Z]} ( \tilde{Z} ))  \equiv 0
\quad\mbox{for}\quad
n \in \ZZ ~.
\eeq

 (1) We now prove \eqn{529} by induction on the length $m$ of the
formula \eqn{528}. The base case $m=1$ is clearly true. For the
induction step, write $[\tau ] = [S'] [\tau ']$, where $[S'] =[S_1]^{\ep_1}$
and $[\tau '] = [S_2]^{\ep_2} \cdots [S_m]^{\ep_m}$.
Lemma 4.4 of part II gives
$$
M_{[\tau]} (\tilde{Z}) = M_{[S']} (\tilde{Z})  +M_{[\tau']} (\tilde{Z}) +
M_{[S']} (M_{[\tau']} ( \tilde{Z} )) ~.
$$
Applying the induction hypothesis to the terms on the right side
gives
 \beql{535}
  M_{[\tau]} (\tZ) = M_{[S'_Z]} (\tZ) + M_{[S'_Y]}
(\tZ ) + M_{[\tau'_Z]} (\tZ ) + M_{[\tau'_Y]} ( \tZ ) + M_{[S']}
(M_{[\tau'_Z]} (\tZ) + M_{[\tau'_Y]} ( \tZ )) ~.
\eeq
Next, the
relations \eqn{533} and \eqn{534} imply that
 \beql{536} M_{[S']}
(M_{[\tau'_Z]} ( \tZ ) +M_{[\tau'_Y]} ( \tZ )) =
M_{[S'_Y]} (M_{[\tau'_Y]} (\tZ )) + M_{[S'_Z]} (M_{[\tau'_Z]} (\tZ )) ~,
\eeq
using the fact that $[S']$ equals one of $[S'_Y]$ and $[S'_Z ]$
and the other term on the right side of \eqn{536} is identically
zero since $M_{[I]} (f) =0$ for all functions $f$. Substituting
\eqn{536} into \eqn{535},
 and applying
 Lemma 4.4 of part II
 twice to the right side yields
\begin{eqnarray}\label{537}
M_{[\tau]} (\tZ ) & = &
M_{[S'_Z][\tau'_Z]} (\tZ ) + M_{[S'_Y][\tau'_Y]} (\tZ ) \nonumber \\
& = & M_{[\tau_Z]} (\tZ ) + M_{[\tau_Y]} ( \tZ ) ~,
\end{eqnarray}
which completes the induction step, yielding \eqn{529}.

(2) The subgroup $\hG_Y \subseteq \pi_1 ( \sN , \bx'_0)$ is the image under
$\pi ' : \sM \to \sN$ of the subgroup \\
$\langle [Y_m ] : m \in \ZZ \rangle$ of $\pi_1 (\sM , \bx_0)$.
The map $\pi'$
is injective when restricted to this subgroup, hence \eqn{530} is a direct consequence
of the formula
\beql{537a} M_{[\tau]} (Z) = \sum_{[S] \in \sG}
M_{[S]^{k(S)}} (Z) ~,
\eeq
given in Theorem 4.6  of part II when restricted to the subgroup $\hG_Y$.\\

(3) We observe that
 the homomorphism from $\pi_1(\sN, \bx'_0) $ to $\ZZ$ which maps $[\tau ] \in\pi_1(\sN, \bx'_0)$
given by \eqn{528} to the sum of the exponents of $[Z_0]$
occurring in \eqn{528} has kernel $\hH_0$ generated by
 \beql{WL533}
\hH_0 : = \langle [Z_0]^{k} [Z_1] [Z_0]^{-k} , [Z_0]^{k} [Y_n][Z_0]^{-k},~k,n \in \ZZ \rangle.
\eeq
That is, we have the exact sequence
$$
0 \to \hH_0 \to \pi_1(\sN, \bx'_0) \to  \ZZ \to 0.
$$
Furthermore every element $[\tau ]$ in
$\pi_1(\sN, \bx'_0)$ can be written as $[\tau ] = [\sigma] [Z_0]^k$ for some
$[\sigma]$ in $\hH_0$ and $k \in \ZZ$. Note that the inclusions
\beql{533b}
[ \pi_1(\sN, \bx'_0),  \pi_1(\sN, \bx'_0)] \subset \hH_1 \subset \hH_0.
\eeq
follow using \eqn{519b} and \eqn{519c}.

We next claim that the the monodromy functions
$M_{[\tau]} ( \tZ )$ vanish identically for $[\tau ] \in \Gamma'= [\hH_0, \hH_0]$.
 To prove the claim,  Theorem 4.6
of part II shows that the
monodromy functions $M_{[\tau]} (Z)$ vanish identically for
$[\tau]$ in the commutator subgroup $\Gamma= (\pi_1(\sM, \bx_0))^{'}$.
Consequently all monodromy functions of $\tilde{Z}$ vanish on the
image group $\Gamma^{'} := (\pi')_{\ast}(\Gamma)$.
Now   $(\pi')_{\ast}(\pi_1(\sM, \bx_0)) = \hH_0$, so  we conclude
\beql{534c}
\Gamma^{'} = (\pi')_{\ast}(\Gamma)= [ \hH_0, \hH_0],
\eeq
which is \eqn{531}.\\

We next observe that $\hH_0$ is a normal subgroup of $\pi_1(\sN, \bx'_0)$
because we have the inclusion of the commutator subgroup
$$
(\pi_1(\sN, \bx'_0))^{'} \subset \hH_1 \subset \hH_0.
$$
It follows that $\Gamma^{'}=[\hH_0,\hH_0]$ is a normal subgroup of $ \pi_1(\sN, \bx'_0)$
because it is a characteristic subgroup of the normal subgroup $\hH_0$
of the group (cf. Robinson \cite[1.5.6 (3)]{Rob82}).
In addition the inclusion  $(\pi_1(\sN, \bx'_0))^{'} \subset \hH_1$
implies that
\beql{534e}
(\pi_1(\sN, \bx'_0))^{'' } \subset  [\hH_1, \hH_1] \subset [\hH_0, \hH_0] = \Gamma^{'}.
\eeq
It follows that $\pi_1(\sN, \bx'_0)/ \Gamma'$ is solvable in two steps, with
$ \Gamma^{'}   ~ \triangleleft  ~\hH_1 ~ \triangleleft ~  \pi_1(\sN, \bx'_0)$.
 ~~~$\bsq$\\

The formulas  in Theorem~\ref{th52} together with those of
Theorem~\ref{th51} suffice to evaluate any monodromy function
$M_{[\tau]} (\tilde{Z} )$ for $[\tau] \in \pi_1 ( \sN , \bx'_0)$.
Indeed, for  $[\tilde{\tau}] \in \pi_1 ( \sN , \bx'_0)$, we may write
$[\tilde{\tau}] = [\tau] [Z_0]^n$ for some
$[\tau] \in \hH_0$ and some integer $n$.
Since $M_{[Z_0]^n} ( \tilde{Z} ) =0$ by Theorem \ref{th51}, we have
$$
M_{[\tilde{\tau}]} ( \tZ ) = M_{[\tau]} ( \tZ ) = M_{[\tau_Z]} ( \tZ ) + M_{[\tau_Y]} ( \tZ ),
$$
where the last equality follows from \eqn{529}.
Note that $[\tau_Y] = [\tilde{\tau}_Y]$ so that
$M_{[\tau_Y]} ( \tZ )$ is given by \eqn{530}.
The element $[\tau_Z ]$ lies in the group
\beql{555a}
\hH_{Z} := \langle [Z_0]^{k} [Z_1] [Z_0]^{-k}, k \in \ZZ \rangle.
\eeq
As $M_{[\sigma]} ( \tZ )$ vanishes for $[\sigma]$ in the commutator subgroup
$\hH_Z^{'} = [\hH_Z, \hH_Z]$, we find
$$M_{[\tau_Z]} ( \tZ ) = \sum_{k \in \ZZ}
M_{([Z_0 ]^{k} [Z_1] [Z_0]^{-k} )^{h(k)}} ( \tZ ) ~,
$$
in which $h(k)$ is the sum of exponents of
$[Z_0]^{k} [Z_1] [Z_0]^{-k}$ occurring in $[\tau_Z ]$.
Recall that for $(s,z,c) \in \Delta$, we have, using \eqn{519aa},
\begin{eqnarray}\label{WL536}
M_{[Z_0]^{k} [Z_1] [Z_0]^{-k}} (\tZ ) (s,z,c) & = &
 - \frac{(2 \pi )^s e^{\frac{\pi is}{2}}}{\Gamma (s)} f_{k}(s,z,c) ~,
\\
\label{WL537}
M_{([Z_0]^{k} [Z_1] [Z_0]^{-k} )^{-1}} (\tZ ) (s, z,c)& = &
- e^{-2 \pi is} M_{[Z_0]^{k} [Z_1] [Z_0]^{-k}} (\tZ )(s, z, c) ~,
\end{eqnarray}
and
\beql{WL538}
M_{([Z_0]^{k} [Z_1] [Z_0]^{-k})^{\pm j}} (\tZ )(s, z, c) =
\frac{e^{\pm 2 \pi ijs}-1}{e^{\pm 2\pi is}-1} M_{([Z_0]^{k} [Z_1][Z_0]^{-k})^{\pm 1}} (\tZ)(s, z, c)
~\mbox{for}~k \ge 1\,.
\eeq
These formulas give a way to evaluate any $M_{[\tau]} ( \tZ )(s, z, c)$ explicitly.\\

\paragraph{{\em Proof of Theorem ~\ref{th21}.}}
The result  follows from Theorem~\ref{th51} and Theorem \ref{th52}. $~~~\Box$\\

%
%
%

\subsection{Extended analytic continuation}\label{sec35}

We  now obtain the extended analytic continuation of the Lerch transcendent,
to $\sN^{\#}$, the manifold obtained from $\sN$
by gluing in  the regions ${\sV} (c= n) :=\{ (a, c) :  a \in \CC, c = n\}$  for   $c=n \in \ZZ_{\ge 1}$.

\begin{theorem}\label{Nth23}
{\em (Lerch Transcendent Extended Analytic Continuation )}
The Lerch transcendent $\Phi (s,z,c)$ analytically
continues to a single-valued function $\tilde{Z}= \tilde{Z}(s, z, c, [\gamma])$
 on the universal
cover
 $\tilde{\sN}^{\#}$ of
 $$
 \sN^{\#} =\{ (s,z,c) \in \CC \times \left( \PP^1 ( \CC )\smallsetminus \{ 0,1, \infty \} \right)
  \times \left( \CC \smallsetminus \ZZ_{\le 0} \right) \}.
 $$
All monodromy functions $M_{[\tau]}(\tilde{Z})$ vanish identically for
$[\tau]$ in the second commutator subgroup $(\pi_1( \sN^{\#}, \bx'_0))'' $
of $\pi_1( \sN^{\#}, \bx'_0)$. Thus the analytic continuation becomes  single-valued on a
covering manifold of $\sN^{\#}$ which has  a two-step solvable group of
deck transformations.
\end{theorem}

\paragraph{Proof.} The
fact that the possible singularities at
$c=n$ are removable follows almost immediately  from the corresponding result in
Part II \cite[Theorem 2.3]{LL2}.  It remains only to  check that the additional monodromy functions
in Theorems~\ref{th51} and \ref{th52} remain  holomorphic at points $(s, a, c)$
with  $c= n \ge 1$, in $\tilde{\sN}^{\#}$. This
 is manifested from their form given in  Theorem \ref{th51}, i.e. the only locations where they
are possibly not holomorphic is  points $(s, z,c )$ with $z= 0, 1$ or   $c \in \ZZ_{\le 0}$.

The manifold  $\tilde{\sN}^{\#}$ has a smaller fundamental group than $\sN$,
which is obtained
as a quotient of $\pi_1(\sN, \bx'_0)$ by setting all
 generators $\{ [Y_n]: n \ge 1\}$ equal to the identity. However
the vanishing of all relevant monodromy for $M_{[Y_n]}(f)$ for $n \ge 1$, for $f= \tilde{Z}$
or $f= M_{[\tau]}(\tilde{Z})$ a monodromy function,
allows the conclusion that all monodromy functions vanish identically for
$[\tau]$ in the second commutator subgroup $(\pi_1( \sN^{\#}, \bx'_0))'' $.
This follows  from the corresponding assertion of Theorem~\ref{th52}
for $(\pi_1(\sN, \bx'_0))'' $.
$~~~\bsq$ \\


\begin{rem}\label{remark38}
{\em
The {\em k-th lower central series group}
 ${\bf \gamma}_{(k)}(\hG)$ of a group $\hG$ is defined recursively by
\beql{Z122}
{\bf \gamma}_{(k)}(\hG) := [{\bf \gamma}_{(k-1)}(\hG)~ :~\hG].
\eeq
Each quotient group $\hG/ {\bf \gamma}_{(k)}(\hG)$ is a
nilpotent group. It can be shown that the group $\pi_1(\sN, \bx'_0) / \Gamma '$ is
not nilpotent, where $\Gamma'= [ \hH_0, \hH_0]$ is the group given in Theorem~\ref{th52}.
It follows  that the lower central series group
$\gamma_{(k)} (\pi_1(\sN, \bx'_0)) \not\subseteq \Gamma '$ for all $n \ge 1$. Thus
the quotient group $\pi_1(\sN, \bx'_0)/ \Gamma '$ falls outside the framework
studied by Deligne \cite{De89}.
}
\end{rem}

%
%
%

\section{Differential-difference operators and monodromy functions}\label{sec4}
\setcounter{equation}{0}

As discussed in part II, the  Lerch zeta function $\zeta(s,a,c)$ satisfies two differential-difference
equations. We introduce the operators
\beql{599a}
\hD_{L}^{-} =  \frac{1}{2 \pi i} \frac{\partial}{\partial a} + c   \quad \text{and} \quad \hD_{L}^{+} =  \frac{\partial}{\partial c}.
\eeq
Then we have
\beql{600a}
\hD_{L}^{-}\, \zeta(s, a, c) = \zeta(s-1, a ,c),
\eeq
and
\beql{600b}
  \hD_{L}^{+} \,\zeta(s, a, c) = -s \zeta(s+1, a ,c).
 \eeq
 In consequence the analytically continued Lerch zeta function
 satisfies the  linear partial differential equation
 \beql{600c}
\hD_L \,\zeta(s, a, c, [\gamma]) = -s \zeta(s, a, c, [\gamma])
\eeq
in which
\beql{600f}
\hD_L := \hD_L^{-} \hD_{L}^{+} = \frac{1}{2 \pi i} \frac{\partial}{\partial a} \frac{\partial}{\partial c} + c\frac{\partial}{\partial c}.
\eeq
The Lerch transcendent $\Phi (s,z,c)$ then satisfies suitable differential-difference equations
 and a linear partial differential equation inherited from these.
The substitution $z = e^{2 \pi ia}$ transforms the
Lerch differential operator to the {\em polyzeta operator}
\beql{601a}
\he{D}_{\Phi} := \hD_{\Phi}^{-} \hD_{\Phi}^{+} = \left( z \frac{\pt}{\pt z} + c \right) \frac{\pt}{\pt c} ~.
\eeq
This lifts to an operator on functions on the universal cover $\tilde{\sN} \equiv \tsM$
which is {\em equivariant} with respect to the group $\sG$ 
of diffeomorphisms of $\tsN$ that preserve the projection from
$\tilde{\sN}$ to $\sN$; the group $\sG$  is isomorphic to  $\pi_1 (\sN, \bx'_0 )$.


\begin{theorem}\label{th61}
{\em (Lerch Transcendent Differential-Difference Operators)}

(1) The analytic continuation $\tZ (s, z, c, [\gamma])$
of the Lerch transcendent $\Phi(s, z, c)$ on the universal cover
$\tilde{\sN}$
 satisfies
the two differential-difference equations
\beql{601}
\left( z \frac{\partial}{\partial z} + c\right) \tilde{Z}(s, z, c, [\gamma]) = \tilde{Z}(s-1, z, c, [\gamma_{-}]),
\eeq
and
\beql{602}
 \frac{\partial}{\partial c} \tilde{Z}(s, z, c, [\gamma]) = -s \tilde{Z}(s+1, z, c, [\gamma_{+}]),
\eeq
in which $[\gamma_{+ }] $ and $[\gamma_{-}]$ denote  paths in $\tilde{\sN}$ which first traverse
$\gamma$ and then traverse a path from the endpoint of $\gamma$ that changes
the $s$-variable only, moving from  $s$ to $s \pm 1$, respectively.

(2) The analytic continuation $\tilde{Z}(s, z, c, [\gamma])$ on $\tilde{\sN}$
satisfies the linear partial differential equation
\beql{603}
\he{D}_{\Phi} (\tilde{Z})(s,z,c,[\gamma] ) = -s \tilde{Z}(s, z, c, [\gamma] ) ~,
\eeq
where
$\he{D}_{\Phi} = z \frac{\pt}{\pt z} \frac{\pt}{\pt c} + c \frac{\pt}{\pt c}$,

(3) For each $[\tau ] \in \pi_1 ( \sN , \bx'_0 )$, the
  Lerch transcendent monodromy function \\
  $M_{[\tau]} ( \tZ )(s, z, c, [\gamma])$
satisfies on $\tilde{\sN}$ the two differential-difference equations
and the linear differential equation above.
\end{theorem}

\paragraph{\em Proof.}
(1), (2) These results for  $\tZ (s, z, c, [\gamma])$  follow directly from Theorem 5.1
of part II by a change of variable $z= e^{2 \pi i a}$. Here we have
$$
z  \frac{\partial}{\partial z}= \frac{1}{2 \pi i} \frac{\partial}{\partial a}.
$$
Thus in the $(s,z, c)$-variables,
the  corresponding differential-difference
equations satisfied by the Lerch transcendent are
\beql{604a}
\left( z \frac{\partial}{\partial z} + c\right) \Phi(s, z, c)= \Phi(s-1, z, c)
\eeq
and
\beql{604b}
 \frac{\partial}{\partial c} \Phi(s, z, c) = -s \Phi(s+1, z, c),
\eeq
and the corresponding  differential equation lifts to
that given in \eqn{603} by analytic continuation.

(3) The
monodromy functions satisfy   the same differential-difference equations
and differential equation because
the differential operators $z\frac{\partial}{\partial z}$ and $\frac{\partial}{\partial c}$
are equivariant with respect to the covering map from $\tilde{\sN}$ to $\sN$.
~~~$\bsq$\\

%
%
%

We now study  the restricted Lerch transcendent monodromy functions
$M_{[\tau']}^s ( \tZ )$  obtained by holding the variable
$s$ fixed, in a fashion analogous to \S7
 of part II \cite{LL2}.


 \begin{defi} \label{Nde31}
 {\em
The {\em Lerch transcendent
monodromy space} $\sW_s$ is the $\CC$-vector
space
spanned by all functions $\sQ_{[\tau']}^s ( \tZ )$ (as
$[\tau']$ varies) regarded as functions on the (simply connected) submanifold
\beql{603a}
 \tsN_s := (\pi)^{-1}(\{s\} \times (\PP^{1} (\CC)\smallsetminus \{0, 1, \infty \} ))
  \times ( \CC \smallsetminus \ZZ) ~\subset~  \tsN,
\eeq
where $\pi: \tilde{\sN} \to \sN$ is the universal covering map.
Here $\tZ^s := \sQ_{[id]}(\tZ)$ is the restriction of $\tZ$ to $\tsN_s$.
This
vector space is a direct  sum of one dimensional spaces:
 \beql{604}
 \sW_s := \CC[\tZ^s] \bigoplus \left( \sum_{ [\tau'] \in \pi_1 ( \sN , \bx'_0)}  \CC [M_{[\tau']}^s ( \tZ )] \right) ,
 \eeq
consisting of all finite linear combinations of the given
countable generating set of
vectors.
}
\end{defi}

These vector spaces have the following properties.
\begin{theorem}\label{th62}
{\em ( Lerch Transcendent Monodromy Space)}
The Lerch transcendent monodromy space $\sW_s$ at $s$ depends on the parameter
$s \in \CC$ as follows.
\begin{itemize}
\item[ (i)] (Generic Case)
If $s \not\in \ZZ$, then $\sW_s$ is an infinite-dimensional vector space,
and has as a basis the set of functions
\beql{607}
\{ M_{[Z_0]^{k} [Z_1] [Z_0]^{-k}}^s (\tZ ) : k \in \ZZ \} \cup
\{ M_{[Y_n]}^s ( \tZ ) : n \in \ZZ \} \cup \{ \tZ^s\} ~.
\eeq

\item[ (ii)] (Positive Integer Case)
If $s = m \in \ZZ_{> 0}$, then $\sW_m$ is an infinite-dimensional
vector space, and has as a basis the set of functions
\beql{606}
\{
M_{[Z_0]^{k} [Z_1 ] [Z_0]^{-k}}^m ( \tZ ) : k \in \ZZ \} \cup \{ \tZ^m\}~.
\eeq

\item[ (iii)] (Nonpositive Integer Case)
If $s = -m  \in \ZZ_{\le 0}$, then all Lerch transcendent monodromy functions
vanish identically, i.e.
\beql{605}
M_{[\tau]}^{-m} ( \tZ ) =0 \quad\mbox{for all}\quad
[\tau ] \in \pi_1 ( \sN , \bx'_0 ) ~.
\eeq
Thus $\sW_{-m} = \CC[ \tZ^{-m}]$ is a one-dimensional vector space.
\end{itemize}
\end{theorem}

\paragraph{\em Proof.}
We establish these cases in reverse order.

(iii)~Theorem \ref{th51} shows that for $s = -m \in \ZZ_{\le 0}$ the
monodromy functions of all generators $[S] \in \sG'$ and their
inverses vanish identically.
In the case of $M_{[Z_1]}^{-m} ( \tZ )$ this is because
$\frac{1}{\Gamma (-m)} =0$.
This carries over to all $[\tau' ] \in \pi_1 ( \sN , \bx'_0)$ by
induction in the length of a word expressing $[\tau' ]$ in terms
of the generators. Only the original function $\tZ$ remains.

(ii)~For $s= m \in \ZZ_{\ge 1}$, the monodromy of all
$\{M_{[Y_n]}^m ( \tZ ) : n \in \ZZ \}$ vanishes identically.
Theorem \ref{th52} implies  that the space $\sW_m$ is spanned by $\tZ^s = \sQ_{[id]}(\tZ)$
together with $\{M_{[Z_0]^{k} [Z_1] [Z_0]^{-k}}^m ( \tilde{Z} ) : n \in \ZZ \} \cup \{\tilde Z^m\}$.
The formulas in Theorem \ref{th51} (i) then indicate that
each $M_{[Z_0]^{k} [Z_1] [Z_0]^{-k}}^m ( \tilde{Z} )$ is a nonzero multiple
of $f_{n} |_{s=m}$,
where $f_n$ at $(s, z, c) \in \Delta$ is given by \eqref{WL511}.
Thus
 \beql{N309}
 \sW_m :=\CC[\tZ^s] \bigoplus \left( \oplus_{n \in \ZZ}  \CC [f_n |_{s=m} ].\right)
 \eeq
Each finite subset of the functions $f_n (s, z, c)|_{s = m}$ together with  $\tZ^s$,
 is easily checked to be linearly
independent in a small neighborhood of the point
$(m,-1, \frac{1}{2} )$, so (ii) follows.

(i)~The proof in the generic case parallels that of Theorem 7.1
of part II (\cite{LL2}). It also makes use of the independence formula
\eqn{529} of Theorem \ref{th52}.  ~~~$\bsq$\\

\begin{rem}\label{remark44}
  {\em
  Theorem~\ref{th62} shows that values  $s \in \ZZ$  are
``special values''  in the sense that  the monodromy functions satisfy
 non-generic linear relations at these values. We note the coincidence that these same points $s\in \ZZ$
 are ``special values'' in the sense of number theory, in  that the values of these functions at these points
  encode important arithmetic information, discussed below.
  Non-positive integer $s=-m \le 0$ are especially interesting
because all monodromy functions vanish identically: therefore the values at
these points are  well-defined on the base manifold $\sN$, without having to lift
to any covering manifold. This observation strengthens the
observation made in part II of the vanishing property of monodromy for
the Lerch zeta function, because we have a larger set of monodromy functions.
The ``well-definedness'' property of these values
seems  particularly significant in that these  values contain arithmetic
information:  $p$-adic $L$-functions
can be obtained by $p$-adic interpolation through values at these points.
More precisely, results in \S6 establish:
\begin{enumerate}
\item[(i)]
The values of the periodic zeta function $F(a, s) = \sum_{n=1}^{\infty} \frac{e^{2\pi i na}}{n^s}$
at negative integers $s= -m \le 0$,
corresponding to a
singular strata degenerations of the  Lerch transcendent, are  recoverable
as limits of  non-singular Lerch transcendent values in this paper.
\item[(ii)]  The classical  $p$-adic $L$-functions
 can be constructed by $p$-adic interpolation starting from values
of $F(a, -m)$.
\end{enumerate}
}
\end{rem}


\begin{rem}\label{remark45}
{\em
Theorem~\ref{th62} treats the monodromy functions
as multivalued functions of  two variables $(z, c)$.
If the variable $c$ is  also fixed, so that
the monodomy functions depend on only the  variable $z$, then
in Section \ref{sec8}  we show that
the monodromy vector spaces become finite-dimensional when $s=m \in\ZZ_{>0}$ is a
positive integer, see
Theorem~\ref{Nth81}.
}
\end{rem}

%
%
%

\section{Specialization of Lerch transcendent: $s$ a non-positive integer }\label{sec5}
\setcounter{equation}{0}

In this and the next section
we specialize
  $s=-m $ ($m  \ge 0$)
to be  a nonpositive integer.
In Theorem \ref{th24}  (3) we observed that the values   $s=-m$
are distinguished  by the fact
that all monodromy functions vanish identically as functions of $(z, c)$. At these values of $s$
a great simplification occurs and the resulting two-variable functions
are rational functions of $(z, c)$.
These rational functions of two variables
were determined  by
Apostol \cite[Section 3] {Ap51} in 1951. They were later  studied
by Kanemitsu, Katsurada and Yoshimoto \cite[Section 4]{KKY00},
who obtained various formulas for them, cf. their Theorem 6.

Specificially, for $z=-m$ ($M \ge 0)$ the modified function
\beql{N500a}
Li_{-m}(z,c) := z \Phi(-m, z, c)= \sum_{n=0}^{\infty} (n+c)^m z^{n+1},
\eeq
is a $c$-deformed polylogarithm of negative integer order.
As mentioned above, it is  known that  the functions   $Li_{-m}(z, c)$ meromorphically
continue  as rational functions of two variables $(z, c) \in \CC \times \CC$.
These rational functions  automatically give a
meromorphic continuation in the $(z, c)$-variables to all  integer points $c \in \ZZ$;
thus including the  singular strata points
$c=-n \le 0$ outside the analytic continuation
of part II (\cite{LL2}). \\

In this section we  determine recursions for
these rational  functions and deduce various symmetry
properties they exhibit (Theorem ~\ref{Zle44}) .
We begin with the following expression for $c$-deformed negative polylogarithms
$Li_{-m}(z,c)$.\\

\begin{theorem}\label{th63}
{\em ($c$-Deformed Negative Polylogarithms)}
For $s = -m \in \ZZ_{\le 0}$ the function $Li_{-m} (z,c)$ analytically continues to  a rational
function of $z$ and $c$ on $\PP^1 (\CC ) \times \PP^1 ( \CC )$.
Here $Li_0 (z, c) = zq_0(z)$ and
\beql{615}
Li_{-m} (z,c) = z \left( \sum_{k=0}^m
{\binom{m}{k}}  c^k q_{m-k} (z)\right)  ~,
\quad\mbox{for}\quad m \ge 1,
\eeq
in which the $q_m(z)$ are rational functions given by
\beql{616}
q_0 (z) = \frac{1}{1-z} ~,
\eeq
and
\beql{617}
q_{m+1} (z) = z \frac{d}{dz} (q_m (z)) \quad\mbox{for}\quad m \ge 0 ~.
\eeq
\end{theorem}

\paragraph{\em Proof.}
The case $m = 0$ with $c \ne 0, \infty$ follows immediately
from \eqn{101b},
$$
Li_0(z, c) =z \left( \sum_{n=0}^{\infty} z^n\right) = \frac{z}{1-z}.
$$
Since this function is independent of $c$,
it extends trivially
to all $c \in \PP^1(\CC)$ as $Li_0 (z, c)$.
We next use the identity
$$
\Phi(s-1, z, c) = \left( z \frac{\partial}{\partial z} + c \right) \Phi (s, z, c)
$$
which yields, taking $s=-m$,
\beql{N318}
Li_{-m-1} (z,c) =z  \left( z \frac{\partial}{\partial z} + c \right) \left(z^{-1}
Li_{-m} (z,c)\right).
\eeq
This holds for $c \ne 0, \infty$ when $z$ is in a sufficiently
small circular neighborhood of
the origin (with radius depending on both $c$ and $m$).
The fact that the function $Li_{-m}(z, c)$ has the form \eqn{615} for
all $m \ge 1$ with $q_j(z)$ as described is proved by  induction on
$m \ge 1$ using  \eqn{N318}, taking the differential equation \eqn{617}
as the definition of $q_{m+1}(z)$. The differential equation
\eqn{617} shows that each $q_m(z)$ is a rational function of $z$. ~~~$\bsq$\\

In the remainder of this section we study these ``special values'' with $s = -m \in \ZZ_{\le 0}$ in more
detail.
One easily shows by induction on $m$ that the rational function $q_m (z)$
has the form
\beql{620}
q_m (z) = \frac{r_m (z)}{(1-z)^{m+1}},
\eeq
in which $r_m(z) \in \ZZ [z]$ is a monic polynomial of degree $m$.
For $m \ge 1$ we have $r_m (0) =0$ and $r_m (1) = m!$.
The first few values of
$r_m(z)$ are given in the following table.
\begin{table}[htb]
$$
\begin{array}{|c|c|} \hline
m & r_m (z) \\ \hline
1 & z \\
2 & z^2 + z \\
3 & z^3 + 4z^2 + z \\
4 & z^4 + 11z^3 + 11z^2 +z \\
5 & z^5 + 26z^4 + 66z^3 + 26 z^2 + z \\ \hline
\end{array}
$$
\caption{Values of $r_m (z)$.}
\label{ta61}
\end{table}

To evaluate the rational functions $q_m(z)$ we use the following result, which is
due to Apostol \cite[(3.1)]{Ap51}.

\begin{lemma}\label{Zle41}
Let
\beql{Z45a}
G(z, c; u) := \sum_{m=0}^\infty Li_{-m} (z,c) \frac{u^m}{m!}
\eeq
be the exponential generating function for the functions
$\{Li_{-m}(c, z): m \ge 0\}$.
This series converges absolutely on
the region $\{ (z, c, u):~|z| < e^{-|u|},  c \in \CC\}$, and is there given by
\beql{Z46}
G(z, c; u) =  \frac{ze^{cu}}{1-ze^u} ~.
 \eeq
 The  right side of \eqn{Z46} gives a meromorphic continuation of  $G(z, c;u)$ to $(z,c, u) \in \CC^3$.
\end{lemma}

\paragraph{\em Proof.}
For $|z| < 1$ we have
$$Li_{-m} (z,c) = \sum_{k=0}^\infty z^{k+1} (k+c)^m,$$
where the series converges absolutely for all $c \in \CC$.
If in addition
$|z| <  e^{-|u|}$, then we have the bounds
\begin{eqnarray*}
\sum_{m=0}^\infty \sum_{k=0}^\infty |z|^{k+1} (k+ |c| )^m \frac{|u|^m}{m!}
& = &
\sum_{m=0}^\infty \sum_{k=0}^\infty |z|^{k+1}
\frac{((k+|c|)|u|)^m}{m!} \\
& \le  & \sum_{k=0}^\infty |z|^{k+1} e^{(k+|c|)|u|} = e^{|cu|-|u|}
\left( \frac{1}{1- |z|e^{|u|}} \right) < \infty ~.
\end{eqnarray*}
Thus when $|z| < e^{-|u|}$, we have
\begin{eqnarray*}
\sum_{m=0}^\infty Li_{-m} (z,c) \frac{u^m}{m!} & = &
\sum_{m=0}^\infty \left( \sum_{k=0}^\infty z^{k +1}(k+c)^m \right) \frac{u^m}{m!} \\
& = & \sum_{k=0}^\infty z^{k+1} \left( \sum_{m=0}^\infty
\frac{((k+c)u)^m}{m!} \right)
\\
& = & \sum_{k=0}^\infty z^{k+1} e^{(k+c) u} = \frac{ze^{cu}}{1-ze^u} ~,
\end{eqnarray*}
which is \eqn{Z46}.~~~$\bsq$\\

We  now establish
 various properties of the rational functions
$q_m (z)$ appearing in these special values.


\begin{theorem}\label{Zle44}
{\em (Properties of $q_m(z)$)}
For $m \ge 0$ the rational functions
$q_m (z) = \frac{r_m (z)}{(1-z)^{m+1}}$, where $r_m (z)$
is a monic polynomial of degree exactly $m$,  having $r_m (1) = m!$.
These polynomials have the following properties.

(i) {\em (Laurent expansion)}. The Laurent expansion
$q_m (z) = \sum_{k=0}^{m+1} \frac{a_{m,k}}{(1-z)^k}$ around $z=1$ is
given by
\beql{Z421}
r_m (z) = \sum_{k=0}^{m+1} a_{m,k}(1-z)^{m+1-k}
~, \eeq
 in which  $a_{m, 0} = 0$ for $m \geq 0$,  and for $1 \le k \le m+1$ we have
\beql{Z422}
a_{m,k} = (-1)^m \sum_{l=0}^{k-1} (-1)^l {\binom{k-1}{l}} (l+1)^m
~. \eeq

(ii) {\em (Reflection symmetry)}. For $m \ge 1$,
\beql{Z423}
 z^{m+ 1} r_m \left( \frac{1}{z} \right) = r_m (z) ~.
  \eeq

(iii) {\em (Recursion)}. The polynomials $r_m(z)$ satisfy
\beql{Z424}
r_m (z) = z \sum_{j=1}^m {\binom{m}{j}} r_{m-j} (z) (1-z)^{j-1} ~.
\eeq
\end{theorem}

\paragraph{\em Proof.}
We specialize  the $c$ variable to $c=0$ in Theorem ~\ref{th63} to obtain
\beql{Z47}
Li_{-m} (z,0) = zq_m (z).
\eeq
Letting $c=0$ in Lemma~\ref{Zle41} and dividing by $z$ then
yields the exponential generating function
\beql{Z48}
\sum_{m=0}^\infty q_m (z) \frac{u^m}{m!} = \frac{1}{1-ze^u}.
\eeq

(i) Theorem ~\ref{th63}
implies that the  rational function  $q_m(z)$ has  the form $\sum_{k = 0}^{m+1} \frac {a_{m,k}}{(1 - z)^k}$.
 Substituting this expression for
$q_m(z)$ into  the left hand side
of (\ref{Z48})
and interchanging the order of summation over $m$ and $k$ yields
\begin{eqnarray}\label{leftLaurent}
\sum_{m=0}^{\infty} q_m(z)\frac{u^m}{m!} =
\sum_{m=0}^{\infty}\sum_{k = 0}^{m+1} \frac {a_{m,k}}{(1 -
z)^k}\frac{u^m}{m!} = \sum_{k=1}^{\infty} \frac{1}{(1-z)^k}
\left( \sum_{m=k-1}^{\infty} a_{m,k}\frac{u^m}{m!}\right).
\end{eqnarray}
We express the right hand side of (\ref{Z48}) also as an infinite
series in powers of $\frac{1}{1-z}$ and $u$, valid for $|1 -
e^{-u}| < |1 - z|$:
\begin{eqnarray}\label{rightLaurent}
\frac{1}{1 - ze^u} &=& \frac{e^{-u}}{e^{-u} - z} =
\frac{e^{-u}}{(1-z) - (1 - e^{-u})} = \frac{1}{1-z}\cdot \frac
{e^{-u}}{1 - \frac{1-e^{-u}}{1-z}} \nonumber \\
&=& \frac{1}{1-z}\cdot e^{-u} \sum_{l =0}^{\infty}(\frac{1-e^{-u}}{1-z})^l = \sum_{k=1}^{\infty}
\frac{1}{(1-z)^k}\cdot e^{-u}(1 - e^{-u})^{k-1} \nonumber \\
&=& \sum_{k=1}^{\infty}
\frac{1}{(1-z)^k} \left(\sum_{l=0}^{k-1}(-1)^l{\binom{k-1}{l}}e^{-(l+1)u}\right)
\nonumber \\
&=& \sum_{k=1}^{\infty} \frac{1}{(1-z)^k}
\left(\sum_{l=0}^{k-1}(-1)^l{\binom{k-1}{l}}\left( \sum_{m=0}^{\infty}(-1)^m(l+1)^m\frac{u^m}{m!}\right)\right).
\end{eqnarray}
Comparing the coefficients of $\frac{1}{(1-z)^k}\frac{u^m}{m!}$ in
(\ref{leftLaurent}) and (\ref{rightLaurent}) gives the expression
for $a_{m,k}$.

(ii) The proof of this identity also uses (\ref{Z48}). More
precisely, replacing $z$ by $ \frac{1}{z}$ in (\ref{Z48}), we
obtain
\begin{eqnarray*}
\sum_{m = 0}^{\infty} q_m(\frac{1}{z}) \frac{u^m}{m!} = \frac{1}{1
- \frac{1}{z}e^u} = \frac{z e^{-u}}{ze^{-u} - 1} = 1 - \frac{1}{1
- ze^{-u}} = 1 - \sum_{m = 0}^{\infty} q_m(z) \frac{(-u)^m}{m!}.
\end{eqnarray*}
This gives the relation
\begin{eqnarray*}
q_m(\frac{1}{z}) = (-1)^{m+1}q_m(z) \quad {\rm for} \quad m \ge 1.
\end{eqnarray*}
Since $$q_m(z) = \frac{r_m(z)}{(1-z)^{m+1}}~,$$ the above identity
can be rewritten as $$ \frac{z^{m+1}
r_m(\frac{1}{z})}{(z-1)^{m+1}} = \frac{r_m(\frac{1}{z})}{(1 -
\frac{1}{z})^{m+1}} = (-1)^{m+1}\frac{r_m(z)}{(1-z)^{m+1}}~,$$
which proves (ii).

(iii) Dividing both sides of (\ref{Z424}) by $(1-z)^{m+1}$, we
convert (\ref{Z424}) to an equivalent form
\begin{eqnarray}\label{Z424'}
q_m(z) = \frac{z}{1-z} \sum_{j = 1}^m {\binom{m}{j}} q_{m-j}(z)
\quad {\rm for} \quad m \ge 1.
\end{eqnarray}
When $m = 1$, the right hand side is $\frac{z}{1-z} q_0(z) =
\frac{z}{(1-z)^2} = q_1(z)$. We shall prove (\ref{Z424'}) by
induction on $m$. Assume it holds for some $m \ge 1$. Rewrite the
identity in this case as
\begin{eqnarray}\label{Z424"}
(1-z) q_m(z) = z \sum_{j = 1}^m {\binom{m}{j}} q_{m-j}(z).
\end{eqnarray}
Differentiating both sides of (\ref{Z424"}) and using the identity
(\ref{617}) that $q_{n+1}(z)= z \frac{\partial}{\partial z}( q_n(z))$ for $n \ge 0$, we arrive
at
\begin{eqnarray*}
- q_m(z) + (1-z)\frac{\partial}{\partial z}( q_m(z)) &=& \sum_{j = 1}^m {\binom{m}{j}}
q_{m-j}(z) + z\sum_{j = 1}^m {\binom{m}{j}} \frac{\partial}{\partial z} q_{m-j}(z) \\
&=& \sum_{j = 1}^m {\binom{m}{j}} q_{m-j}(z) + \sum_{j = 1}^m
{\binom{m}{j}} q_{m-j+1}(z)\\
&=& q_0(z) + \sum_{j = 2}^m {\binom{m+1}{j}} q_{m-j+1}(z) +
m q_m(z).
\end{eqnarray*}
In other words,
$$(1-z) \frac{\partial}{\partial z}( q_m(z)) = \sum_{j = 1 }^{m+1} {\binom{m+1}{j}}
q_{m+1-j}(z).$$ Therefore
$$
 q_{m+1}(z) = z \frac{\partial}{\partial z}(q_m(z))=
 \frac{z}{1-z}\sum_{j =1 }^{m+1} {\binom{m+1}{j}},
$$
as desired.
$~~~\bsq$ \\

%
%
%
\section{Double specialization: Periodic zeta function}\label{sec6}
\setcounter{equation}{0}

We now consider the double specialization of the Lerch transcendent setting $s=-m$
a  non-positive integer and setting  $c=0$. The  importance of these
``special values'' $s= -m$ of the $s$-parameter  is that values of zeta and $L$-functions
at these points have arithmetic significance;
they encode  information about the arithmetic structure of number fields.
We have the problem that the  $c$-parameter  value $c=0$ lies on the  singular stratum outside the analytic
continuation of $\Phi(s, z, c)$ given in Section \ref{sec3}.
Our object is to show that sufficiently many of these ``special values'' 
can be  recovered by a limiting process from ``regular stratum'' values of
the Lerch transcendent so as to carry out the number-theoretic construction
of $p$-adic $L$-functions.
We shall   take a limit as $c \to 0$,
and to do this we set   $z =e^{2 \pi i a}$ and will also suppose that
$0 < \Re(a) < 1, $ so that $z \in \CC \smallsetminus \RR_{>0}$.

Concerning the number theoretic significance of the
 values $s=-m$ ($m \ge 0$) is that on a singular stratum of the Lerch transcendent
 giving the   Riemann zeta
function
these values are the important rational numbers
$$
\zeta(-k) = - \frac{B_{k+1}}{k+1},
$$
where the $B_k$ are Bernoulli numbers. We follow the convention
on Bernoulli numbers that they are defined by
$$
\frac{t}{e^t-1} = \sum_{k=1}^{\infty} \frac{B_k}{k!}t^k.
$$
More generally the special values $L(-m, \chi)$ of Dirichlet $L$-functions,
for $m \ge 0$ are arithmetically important algebraic numbers. These values
 can be expressed  as linear combinations of either
Hurwitz zeta function values $\zeta(-m, \frac{a}{d})= \Phi(-m, 1, \frac{a}{d})$ or,
alternatively,  of periodic
zeta function values $P(\frac{a}{d}, -m)= \Phi(-m, e^{2 \pi i \frac{a}{d}}, 1)$ taken at
values  which are roots of unity in the  $z$ variable.

The arithmetic information in the  ``special values'' $L(-m, \chi)$
include $p$-adic regularities captured by   interpolating these values  $p$-adically
to obtain $p$-adic $L$-functions.  Two  different interpolation methods to
construct $p$-adic $L$-functions  are known.
A original construction  of Kubota and Leopoldt \cite{KL64} in 1964 used  interpolation of  Hurwitz zeta function values;
it is presented in Washington \cite[Sect. 5.2, Theorem 5.11]{Wa97}.
A second approach was given in 1977 by
 Morita ~\cite{Mo77}.
 which uses interpolation of periodic zeta function values,
 at the points $(z, a)$ with $a= \frac{j}{p^k}$
having $0 < j < p^k$ and with $z= e^{\frac{2 \pi i j}{p^k}}$. More information on
the periodic zeta function approach to $p$-adic $L$-functions is given in
Amice and Fresnel~\cite{AF72} and Naito~\cite{Na82}.

We use  the periodic zeta function approach to the ``special values'', where limits exist,
rather than the Hurwitz zeta function approach, where the limits do not
exist, see Remark \ref{remark62} below.
The following result shows that it is
possible to take  a limit as $c \to 0^{+}$ of  $z\Phi(s, z, c)$, with
 $z=e^{2 \pi i a}$
 to recover the
values at $s=-m$ of the (analytically continued)  periodic zeta function
$$
F(a,s) = \sum_{n=1}^\infty  \frac{e^{ 2 \pi ina } }{n^s}.
$$
 That is, this limit exists, and these values $F(z,s)$ are extractable from
data in  the nonsingular part of the analytic continuation of
$\Phi(s, z, c)$ at $z=-m$, when approaching the singularity.


\begin{theorem}\label{th43}
{\em (Periodic Zeta Function Special Values)}

(1) For $z \in \PP^1(\CC) \smallsetminus \{0, 1, \infty\}$ and $-m \in \ZZ_{\le 0}$  there holds
\beql{N631A}
Li_{-m} (z, 0) = \lim_{c \to 0^{+}} Li_{-m}(z, c) = \lim_{c \to 0^{+}} z \Phi(-m, z, c),
\eeq
where the limit is taken through values of $c$ in $0< \Re(c)<1$.

 (2) For $0 < \Re(a)  < 1$ the periodic zeta function
$F(a,s) = \sum_{n=1}^\infty  \frac{e^{2 \pi ina}}{n^s}$
analytically continues to an entire function of $s$. In particular,
for $s= -m \in \ZZ_{\le 0}$ there holds
\beql{N631}
F(a, -m) =   e^{-2 \pi i a} Li_{-m}(e^{2 \pi i a}, 0)= q_m(e^{2 \pi i a}),
\eeq
\end{theorem}

\paragraph{\em Proof.}
(1) By Theorem ~\ref{th63}  $Li_{-m}(z, c)$ is a rational function of $(z, c)$, giving
the left equality in
\beql{N632}
Li_{-m}(z, 0) = \lim_{c \to 0^{+} }Li_{-m}(z, c),
\eeq
where the limit is taken over values in  $0<\Re(c)< 1$. Now the equality
$Li_{-m}(z, c) = z\Phi(-m, z, c)$ (which holds by
analytic continuation) gives the result \eqn{N631A}.

(2) The analytic  continuation of $F(a, s)$ to an entire function of $s$
for $0 < \Re(a) < 1$ follows directly from the integral representation
$$
F(a,s) = \frac{1}{\Gamma(s)} \int_{0}^{\infty} \frac{1}{1- e^{2\pi i a} e^{-t}} t^{s-1} dt,
$$
since $e^{2\pi i a}$ stays off the nonnegative real axis.

In  \eqn{510c} we have shown that on the extended fundamental polycylinder
$$\tilde{\Omega}= \{ s \in \CC \} \times \{ a: 0 < \Re(a) < 1 \} \times \{ c: 0< \Re(c)< 1\}$$
we have
\beql{N633}
\Phi(s, e^{2 \pi i a}, c) = \zeta(s, a, c).
\eeq
We now suppose $a$ is real and  apply the results from part I (\cite{LL1}) on
the limiting behavior of the Lerch zeta function $\zeta(s,a, c)$ approaching
the boundary of $\Box^{\circ}$.
Namely, for  $\Re(s) < 0$,
Theorem 2.3 (iii) and equation (2.13)
of \cite{LL1} give
 for $0 < a < 1$ that
\beql{N634}
 \lim_{c \to 0^{+}} \zeta_{\ast}(s, a, c) = \zeta_{\ast}(s, a, 0) = F(a, s),
\eeq
in which  $F(a, s)$ is the analytic continuation to $s \in \CC$ of the periodic zeta function
in the $s$-variable. Here $\zeta_{\ast}(s, a, c)$ is the analytic continuation in
the $s$-variable of
$$
\zeta_{\ast}(s, a, c) = \sum_{n+c>0} e^{2 \pi i n a} (n+c)^{-s},
$$
and it satisfies $\zeta_{\ast}(s, a, c) = \zeta(s, a, c)$ for $(a, c) \in \Box^{\circ}.$
A key point is that  the $n=0$ term $c^{-s}$ approaches $0$ as $c \to 0^{+}$,
when $\Re(s)<0$.

Now  \eqn{N632}-\eqn{N634} give, taking $z=e^{2 \pi i a}$ with
$0< a<1$, and $s= -m$ that
\begin{eqnarray*}
F(-m, s) & = & \lim_{c \to 0^{+}} \zeta(-m, a, c)   \\
& = & \lim_{c \to 0^{+}} \Phi(-m, e^{2 \pi i a}, c) \\
& = & e^{- 2 \pi i a} Li_{-m}(e^{e\pi i a}, 0),
\end{eqnarray*}
where the limits are taken over real $0<c<1.$
Finally, the validity of  \eqn{N631} extends from real $a$ to  $0 < \Re(a) < 1$ by uniqueness of analytic continuation in the $a$-variable.
$~~~\bsq$. \\

\begin{rem}\label{remark62}
{\em
The Kubota-Leopoldt \cite{KL64}  construction of
 $p$-adic $L$-functions in 1964 used interpolation of  special values of the Hurwitz zeta function.
However   it is not possible to obtain Hurwitz zeta function special values directly by a limiting process
involving the Lerch transcendent.  Theorem 2.3 of Part I shows that the
Hurwitz zeta function values $\zeta(s, c)$ at negative integers $s=-m$  cannot
be obtained as a limit of values $\Phi(s, z, c)$, as the parameter $a =\frac{1}{2\pi i } log z$ has $a \to 0$ (resp. $a \to 1$).
Under a variable change this limit corresponds to taking $z \to 1$, which is a
 singular stratum value.  As indicated by Theorem \ref{th43} above we can  indirectly access
Hurwitz zeta function values  by expressing them
as a linear combination of periodic zeta function values,
see  Apostol \cite[Theorem 12.3]{Ap76}. The periodic zeta function values
 are obtainable as limiting values using  Theorem~\ref{th43}.
}
\end{rem}

%
%
%
\section{Specialization of Lerch transcendent: $s$ a positive integer}\label{sec7}
\setcounter{equation}{0}

In this section we treat specialization of variables
related to  polylogarithms, where we specialize the $s$-parameter to
be a positive integer value $s=m \ge 1$. We  will state results in terms of
the {\em extended polylogarithm}
\beql{612aa}
Li_{m}(z, c) := z \Phi(m, z, c) = \sum_{k=0}^{\infty} \frac{z^{k+1}}{(k+c)^{m}}.
\eeq
Recall that the  {\em classical polylogarithm} $Li_m (z)$ is defined by
\beql{613}
Li_m (z) := \sum_{k=1}^\infty \frac{z^k}{k^m} ~,
\eeq
for integer  $m \ge 1$,
see Lewin \cite{L1}, \cite{L2}. The function $Li_m (z,c)$  is related to the classical polylogarithm by
taking $c=1$, with
 \beql{614}
 Li_m (z)= Li_m (z,1)  \quad\mbox{for}~~~ m \ge 1 ~.
\eeq
This function is included in the analytic continuation of the Lerch transcendent
in Theorem~\ref{Nth23}.

First note that for fixed $s \in \CC$
as a function of two variables  $F(z, c):= Li_{s}(z, c)$ satisfies
the linear PDE
\beql{615a}
\left( z \left( z \frac{\partial}{\partial z} +c \right)  \frac{1}{z} \right) \frac{\partial}{\partial c}\,Li_{s}(z, c) =
-s \, Li_{s}(z, c),
\eeq
which easily follows from the linear PDE  \eqn{107} satisfied by $\Phi(s, z, c)$.
Using the identity
 $$
z\left( z \frac{\partial}{\partial z} +c \right) \frac{1}{z}= z \frac{\partial}{\partial z} + c-1
$$
and specializing to $s= m \in \ZZ$, we obtain
\beql{615b}
\left( z \frac{\partial}{\partial z}  \frac{\partial}{\partial c}+c \frac{\partial}{\partial c}- \frac{\partial}{\partial c} \ + m\right) Li_{m}(z, c)=0.
\eeq

We now show  that  for $s=m \ge 1$ there is an analytic continuation of $Li_m(z, c)$
in the two variables $(z, c)$ that meromorphically extends
to  all $c \in \CC$, including positive and negative
integer values.  Note that  this continuation in two variables extends to
nonpositive integer values of $c$ that fall in a ``singular
stratum'' outside the three-variable analytic continuation in Theorem~\ref{Nth23}.

\begin{theorem}\label{th64}
{\em ($c$-Deformed Polylogarithm Analytic Continuation)}
For each integer $m \ge 1$ the function $Li_m(z,c)$ has a
meromorphic continuation in  two variables $(z,c)$ to the
universal cover  $\tilde{\CC}_{0, 1,\infty} \times \CC$
of $~(\PP^1 (\CC ) \smallsetminus \{0,1,\infty \} ) \times \CC$.
For fixed $\tilde{z} \in \tilde{\CC}_{01\infty}$,
 this function  is meromorphic as a function of  $c \in \CC$,
with its singularities consisting of  poles of exact order $m$ at each of the points $c \in \ZZ_{\le 0}.$
\end{theorem}

\paragraph{\em Proof.}
As shown in the proof of Theorem 4.1
of part II, the Lerch zeta function $\zeta (s, a, c)$ extends to a
single-valued holomorphic function on any simply connected domains
contained in
 \beql{622}
 \sO_1^\infty :=\{s: \Re (s) > 0 \} \times
\{ a: a \in \CC \smallsetminus \ZZ \} \times \{ c: \Re (c)
> 0 \} ~.
 \eeq
 Since the monodromy functions over $\sO_1^\infty$
are holomorphic in $c$, for each $m \in \ZZ_{> 0},$ the function
$Li_m (z, c)$ defined for $|z| < 1$ and $\Re (c) > 0$ by
\beql{623}
 Li_m (z,c) =\sum_{k=0}^\infty \frac{z^{k+1}}{ (k+c)^{m} }
 \eeq
extends to a holomorphic function on the covering space
$\tilde{\CC}_{0 ,1 ,\infty} \times \{c~:~ \Re(c) > 0 \}.$ To
determine its behavior on the region $\Re(c) > -N,$ rewrite
\eqn{623} as
\begin{eqnarray}\label{624}
Li_m(z,c) & = & \frac{z}{c^m} + \frac{z^2}{(c + 1)^{m}}+ \dots + \frac{z^{N}}{(c + N - 1)^{m}}
+ \sum_{k = N}^\infty \frac {z^{k+1}}{(k + c)^m}  \nonumber \\
& = & \frac{z}{c^m} + \frac{z^2}{(c + 1)^{m}}+ \dots + \frac{z^{N}}{(c + N - 1)^{m}}
+ z^N Li_m (z, c + N).
\end{eqnarray}
The analytic continuation of $Li_m(z,c)$ to
$\tilde{\CC}_{0, 1, \infty} \times \{c~: \Re(c) > -N\}$  is given
by the right side of \eqn{624}, which clearly shows that for a
fixed $\tilde{z} \in \tilde{\CC}_{0, 1, \infty}$, the function
$Li_m (\tilde{z}, c)$ as a function of $c$ has poles of order $m$
at $c= 0, -1, \dots, -N + 1.$
Theorem~\ref{th64} follows on letting $N \to \infty.$ ~~~$\bsq$\\

%
%
%
\section{Double specialization: Deformed polylogarithm}\label{sec8}
\setcounter{equation}{0}

We next study   the extended polylogarithm
$Li_m( z, c)$   when two variables $s=m \in \ZZ_{\ge 1} $ and $c \in \CC$ are fixed,
and only the parameter $z$ varies. This is a double specialization.
We view the $c$-parameter as giving
a deformation of the polylogarithm, and to emphasize this we rewrite it as
\beql{501a}
Li_{m, c}(z) := Li_{m}(c, z) = z \Phi(m, z, c).
\eeq
We observe that this function of $z$ satisfies a linear ordinary differential equation,
of order $m+1$ with rational function coefficients in $\CC(z)$. We show that
 for all $c  \in \CC$,  its
singular points are
at $\{0, 1, \infty\}$,  and that  is of
Fuchsian type. Then   we determine its monodromy as a function of the parameter $c$.
We observe  that the monodromy representation
varies continuously in the parameter $c$ at nonsingular values $c \in \CC
\smallsetminus \ZZ_{\le 0}$, but has discontinuous jumps
at  ``singular set'' values $c \in \ZZ_{\le 0}$.

\begin{theorem}\label{Nth81}
{\em ($c$-Deformed Polylogarithm Ordinary Differential Equation)}
Let $m \in \ZZ_{\ge 0}$ and let $c \in \CC$
be fixed.

(1) The function
$ F(z)= Li_{m,c}(z)$ satisfies the  ordinary differential equation
\beql{N809a}
D_{m+1}^c F(z) =0,
\eeq
 where
$D_{m+1}^c \in \CC[z, \frac{d}{dz}]$ is the linear ordinary differerential operator
\beql{N810}
D_{m+1}^c := z^2 \frac{d}{dz} \left(\frac{1-z}{z} \right)\left( z \frac{d}{dz} + c-1\right)^m
\eeq
 of order $m+1$.

(2) The operator $D_{m+1}^c$  has  singular
points contained in the set  $\{0, 1, \infty\}$ on the Riemann sphere in $z$,
all of which are regular singular points. In particular,  this equation
is a Fuchsian operator for all $c \in \CC$.

(3) A basis of solutions of $D_{m+1}^c$
 for $c \in \CC \smallsetminus \ZZ_{\le 0}$
is, for $z \in \CC \smallsetminus \{ (-\infty, 0] \cup [1, \infty) \}$,
 given by
\beql{N811}
\sB_{m+1,c} := \{ Li_{m,c}(z), z^{1-c}(\log z)^{m-1}, z^{1-c}(\log z)^{m-2}, \cdots, z^{1-c}  \}.
\eeq
(For $c \in \ZZ_{\le 0}$ the function $Li_{k,c}(z)$
is not well defined.)

(4) A basis of solutions of $D_{m+1}^c$  for $c= -k \in \ZZ_{\le 0}$ is,
for  $z \in \CC \smallsetminus \{ (-\infty, 0] \cup[1, \infty) \}$,
 given by
\beql{N812}
\sB_{m+1,c}^{\ast} := \{ Li_{m,-k}^{*}(z), z^{1-c}(\log z)^{m-1}, z^{1-c}(\log z)^{m-2}, \cdots, z^{1-c}  \},
\eeq
in which
\beql{813}
Li_{m,-k}^{*}(z) := \sum_{{n=0}\atop{n \ne k}}^{\infty} \frac{z^{n+1}}{(n-k)^m} +
\frac{1}{m!} z^{k+1}(\log z)^m.
\eeq
\end{theorem}

To prove this result, we use a preliminary lemma, concerning the form of
the equation.

\begin{lemma}\label{Nlem81a}
The operator
\beql{N821}
D_{m+1}^c = \sum_{k=0}^{m+1} a_{m+1, k}^c (z) \frac{d^k}{dz^k} \,
\eeq
in which
$a_{m+1, k}^c(z) \in \CC[z]$ is of degree at most $j+1$, with factorization
\beql{N822}
a_{m+1,k}^c(z) = \left( \alpha_{m+1, k}(c)  z + \beta_{m+1, k}(c) \right)z^k
\eeq
and the coefficients
$\alpha_{m+1, j}(c), \beta_{m+1, j}(c) \in \ZZ[c]$. The top order
coefficient
\beql{N823}
a_{m+1, m+1}^c(z) = (1-z) z^m
\eeq
has coefficients independent of the parameter $c$.
\end{lemma}

\paragraph{\bf Proof.}  By using  the identity
$z^2 \frac{d}{dz} \left( \frac{1-z}{z} \right) = (1-z)z \frac{d}{dz} -1$
we  may rewrite the definition \eqn{N810} as
\begin{eqnarray}\label{N824}
D_{m+1}^c &= &\left( (1-z) z \frac{d}{dz} -1\right) \left( z \frac{d}{dz} + c-1\right)^m
\nonumber \\
&=& (1-z) \left( z \frac{d}{dz} + c-1\right)^{m+1} - c \left( z \frac{d}{dz} + c-1\right)^m.
\end{eqnarray}
Using the Weyl algebra
commutation relation $\frac{d}{dz}\, z = z \frac{d}{dz} +1$ we obtain by induction on $j \ge 0$ that
$$
\left( z \frac{d}{dz}\right)^j = \sum_{k=1}^j b_{j,k} z^k \frac{d^k}{dz^k}
$$
for certain scalars $b_{jk} \in \ZZ$, including $b_{jj}=1.$ From this we obtain
$$
( z \frac{z}{dz} + (c-1))^m = \sum_{k=0}^{m} f_{m,k}(c) z^k \frac{d^k}{dz^k},
$$
in which each $f_{m,k}(c) \in \ZZ[c]$. Combining this with \eqn{N824} yields
\begin{eqnarray}\label{N826}
D_{m+1}^c &=  & (1-z) \left( \sum_{k=0}^{m+1} f_{m+1,k}(c) z^k \frac{d^k}{dz^k}\right)
- c \left( \sum_{k=0}^{m} f_{m,k}(c) z^k \frac{d^k}{dz^k}\right) \nonumber \\
&=&
 \sum_{k=0}^{m+1} \left(-f_{m+1, k}(c) z + f_{m+1, k}(c) -c f_{m,k}(c)\right)z^k \frac{d^k}{dz^k}.
\end{eqnarray}
This yields \eqn{N822} with
$$
\alpha_{m+1, k}(c) := - f_{m+1, k}(c),~~~\beta_{m+1, k}(c) := f_{m+1, k}(c) -c f_{m,k}(c)
$$
and \eqn{N823} follows easily. $~~~\Box$.\\

\paragraph{\bf Proof of Theorem ~\ref{Nth81}.}
(1) The factor of $z^2$ at the front of the differential operator $D_{m+1}^c$ is
included so that it will have polynomial coefficients, and belong to the
Weyl algebra $\CC[z, \frac{d}{dz}]$, as shown in Lemma~\ref{Nlem81a}.
We check, for $m \ge 1$,
$$
\left( z \frac{d}{dz} + c-1\right) Li_{m, c}(z) = Li_{m-1, c}(z),
$$
where
$$
Li_{0}(c, z) = \sum_{n=1}^{\infty} z^n =  \frac{z}{1-z}
$$
is independent of $c$. This implies $D_{m+1}^c Li_{m,c}(z) =0$.

(2) Using Lemma~\ref{Nlem81a} we have
\begin{eqnarray}\label{N829}
\frac{1}{(1-z)z^{m+1}} D_{m+1}^c &= &
\frac{d^{m+1}}{dz^{m+1}} +
\sum_{k=0}^m \frac{a_{m+1, k}^c(z)}{(1-z)z^{m+1}} \frac{d^k}{dz^k} \nonumber \\
&=& \frac{d^{m+1}}{dz^{m+1}} +
\sum_{k=0}^m \frac{\alpha_{m+1, k}(c)z + \beta_{m+1, k}(c)}{(1-z)z^{m+1-j}} \frac{d^k}{dz^k}.
\end{eqnarray}
We now apply standard criteria for identifying  singular points and determining
if they are regular, given in
Coddington and Levinson \cite[Chap. 4, Theorems 5.1, 6.1, 6.2]{CL55}.
The finite singular points can only be at poles of the coefficients $c_k(z)$ of
$\frac{1}{(1-z)z^{m+1}} D_{m+1}^c
= \sum_{k=0}^{m+1} c_k(z) \frac{d^k}{dz^k}$, which can
occur only at   $z=0, 1$. The condition
for a singular point of first kind (which implies regular singular point) at $z_0 \in \CC$
is that, for $0 \le j \le m+1$, the order of
the pole at $z_0$ of the coefficient $c_k(z)$ is at most $m+1-k$. Now
\eqn{N829} shows that this condition holds at $z=0,1$. The
necessary and sufficient condition for at most a  regular singular
point at $z= \infty$ is that, for $0 \le k \le m+1$, each
 $c_k(z)$ has a zero of order at least $m+1-k$ at $z=\infty$.
This clearly holds in \eqn{N829} as well. Thus the singular points are always a  subset
of $\{0,1, \infty\}$ and each point is either nonsingular or
else a regular singular point is regular, for all $c \in \CC.$ Thus the
differential operator  is Fuchsian for all
$c \in \CC$.
(In fact $z=1$ is not a singular point for $m=0$ and all $c \in \CC$, but
when $m \ge 1$  all three points
$0, 1, \infty$ are genuine
singularities for all values of $c$.)

(3) For integer $k \ge 0$ we have
$$
\left( z \frac{d}{dz} + c- 1\right) z^{1-c} (\log z)^k = k z^{1-c} (\log z)^{k-1} .
$$
It follows that $\{ z^{1-c} (\log z)^j : 0 \le j \le m-1\}$ are annihilated by
$\left( z \frac{d}{dz} + c-1\right)^m$ hence by $D_{m+1}^c$. This shows that the
$m+1$ functions listed above in $\sB_{m+1,c}$ are all in the solution space
of $D_{m+1}^c$. It remains to check that they are linearly independent over
$\CC$. The functions $z^{1-c}(\log z)^k$ are well-defined solutions
for all $c \in \CC$, they are linearly independent from powers of the logarithm,
and none of
them have a singularity at $z=1$.
For $c \in \CC \smallsetminus \ZZ_{\le 0}$ the
function  $Li_{m,c}(z)$ is a well-defined solution,
and it  does have a singularity
there since it diverges approaching this point along the real axis $z \to 1^{+}$.
We conclude that these $m+1$ functions are a basis of
solutions of the differential equation \eqn{N809a}.

(4) The values $c =-k \in \ZZ_{\le 0}$ are singular strata values of the
three-variable analytic continuation.
To see that $D_{m+1}^c\left( Li_{m, -k}^{\ast}(z)\right) =0$, observe that
$$
(z \frac{d}{dz} -k-1)^m ( \sum_{{n=0}\atop{n \ne k}}^{\infty} \frac{z^{n+1}}{(n-k)^m})=
\sum_{{n=1}\atop{n \ne k+1}}^\infty z^n,
$$
while
$$
(z \frac{d}{dz} -k-1)^m \left(z^{k+1} (\log z)^m \right) = m! z^{k+1}
$$
Combining these results gives
$$
(z \frac{d}{dz} -k-1)^m Li_{m, -k}^{\ast}(z) = \sum_{n=1}^{\infty} z^n = \frac{z}{1-z},
$$
which gives the result. Linear independence over $\CC$ is verified as in (3).
$~~~\Box$.\\

We now compute the monodromy representation of this differential equation.
We view $Li_{m, c}(z) := z \Phi(m, z, c)$ as a
multivalued
function of the variable $z$ on the domain $ \PP^1 ( \CC ) \smallsetminus \{0,1, \infty\}$.
The fundamental group
$$
\hGZ := \pi_1( \PP^1 ( \CC ) \smallsetminus \{0,1, \infty\}, -1)= \langle [Z_0 ], [Z_1] \rangle
$$
is a free group on two generators. It
acts on the universal cover $\tilde{\CC}_{0,1,\infty}$ by deck transformations.
We treat the ``non-singular'' case that $c \in \CC \smallsetminus \ZZ_{\le 0}$
and the ``singular'' case $c \in \ZZ_{\le 0}$ separately.


\begin{defi}~\label{de51}
{\em Let $\sW_{m+1,c}$ denote the $(m+1)$-dimensional complex vector space of
functions on
$\PP^1(\CC) \smallsetminus \{0,1,\infty\}$
spanned by the vectors in $\sB_{m+1, c}$  for $c \in \CC \smallsetminus \ZZ_{\le 0}$ (resp. $\sB_{m+1, c}^{*}$ if $c \in \ZZ_{\le 0}$),
making cuts along the real axis $(-\infty, 0)$ and $(1, \infty)$.
The {\em monodromy representation} of $Li_{m,c} (z)$ is the induced action
\beql{N815a}
\rho_{m,c} : \hGZ \to Aut (\sW_{m+1,c} ) \simeq  GL(m+1, \CC) ~.
\eeq
The matrix representation on $GL(m+1, \CC)$ is obtained by viewing
the entries of $\sB_{m+1, c}$ as a
$1 \times (m+1)$ column vector.
}
\end{defi}

We may alternatively view $\sW_{m+1, c}$  as a $(m+1)$-dimensional vector
bundle over the manifold $\sX:= \PP^{1}(\CC) \smallsetminus \{ 0, 1, \infty\}.$

\begin{theorem}\label{thN82}
{\em ($c$-Deformed Polylogarithm Monodromy-Nonsingular Case)}\\
For each integer $m \ge 1$ and each $c \in \CC \smallsetminus \ZZ_{\le 0}$,
the monodromy action on the basis
$$
\sB_{m+1,c} := \{ Li_{m,c}(z), z^{1-c}(\log z)^{m-1}, z^{1-c}(\log z)^{m-2}, \cdots, z^{1-c}  \}
$$
is given by
\beql{N831}
\rho_{m,c} ([Z_0]) := \left(
\begin{array}{cccccc}
1 & 0 & 0 & \cdots & 0 & 0 \\
0 & e^{- 2 \pi ic} & \frac {2 \pi i}{1!}e^{-2 \pi ic} & \cdots &
\frac{(2 \pi i)^{m - 2}}{(m - 2)!}e^{-2 \pi ic}  & \frac{(2 \pi i)^{m - 1}}{(m - 1)!}e^{-2 \pi ic} \\
\vdots & \vdots & \vdots & ~ & \vdots & \vdots \\
0 &  0 &0 & \cdots & e^{-2\pi i c} & \frac{2\pi i}{1!}e^{- 2\pi i c}  \\
0 & 0 & 0 & \cdots & 0 & e^{-2 \pi ic} \\
\end{array}
\right) \,
\eeq
and
\beql{N832}
\rho_{m,c} ([Z_1]) := \left( \begin{array}{cccccc}
1 & -2\pi i & 0 & \cdots & 0 & 0 \\
0 & 1 & 0 & \cdots & 0 & 0 \\
\vdots & \vdots & \vdots & \vdots & \vdots & \vdots \\
0 & 0 & 0 & \cdots & 1 & 0 \\
0 & 0 & 0 & \cdots & 0 & 1
\end{array}
\right) \, .
\eeq
The image of $\rho_{m,c}$ falls in a Borel subgroup of $GL(m+1, \CC)$. The image
 is unipotent if $c \in \ZZ_{>1}$, and is quasi-unipotent if $c \in \QQ \smallsetminus \ZZ.$
\end{theorem}

\paragraph{\bf Proof.}
Around $z=0$ the function $Li_{m, c}(z)$ is analytic, so it has no monodromy.
This gives the first row in the matrix \eqn{N831}. We also have, for all $c \in \CC$,
\begin{eqnarray}\label{N833}
M_{[Z_0]} ( \frac{1}{n!} z^{1-c} (\log z)^n) &= & \frac{1}{n!}\left(
e^{-2\pi i (1-c)} z^{1-c}(\log z + 2 \pi i)^n - z^{1-c} (\log z)^n \right)
\nonumber \\
&=& e^{-2\pi i c} z^{1-c} \left( \sum_{j=1}^n \frac{1}{n!} \left( {{n}\atop {j}} \right) (2 \pi i)^j (\log z)^{n-j}\right)
\nonumber \\
&=& e^{-2\pi i c} z^{1-c}
\left( \sum_{j=1}^n \frac{(2 \pi i)^j}{j!} \left( \frac{1}{(n-j)!} (\log z)^{n-j} \right) \right).
\end{eqnarray}
This gives the remaining rows in \eqn{N831}.

Around $z=1$, we note first that the  $b_{n,c}(z) =  \frac{1}{n!} z^{1-c} (\log z)^n$ in $\sB_{m+1, c}$
are analytic at $z=1$, which give zero entries  above the diagonal in the last $m$ rows of \eqn{N832}.
It  remains to evaluate the monodromy of $Li_{m,c}(z)$.
We use the analytic continuation of the Lerch transcendent given in Section \ref{sec3}.
Note that
Theorem \ref{th51} shows that when $s=m \in \ZZ$, the monodromy of
$\tZ$ along $[Y_n]$ vanishes for all $n \in \ZZ$,
hence we may view the restricted function $\tZ^m(z, c, [\gamma])$ as a function on
$\tilde{\CC}_{01 \infty} \times (\CC \smallsetminus \ZZ_{\le 0} )$.
 We denote by $\tZ_{m,c}(z)$ ( resp.
 $M_{[\tau]}^{m,c}(\tilde{Z})$)  the restrictions of $\tZ^m(z, c, [\gamma])$ ( resp.
 $M_{[\tau]}^m ( \tZ)$) to $\tilde{\CC}_{0,1, \infty} \times \{c\}$.
 We also treat these as functions of one variable, $z$, on the universal cover
$\tilde{\CC}_{0,1, \infty}$.
By \eqn{WL536} we have
\beql{N835}
M_{[Z_0]^{-k} [Z_1][Z_0]^k}^{m,c}( \tZ) = \alpha(m)f_{k}^{m,c}(z),
\eeq
in which the constants
\beql{N836}
 \alpha(m) := - \frac {(2\pi)^m e^{\frac {\pi i m}{2}}}{\Gamma (m)} = - \frac{(2 \pi i)^m}{ (m - 1)!},
 \eeq
 and the functions
\beql{N837}
f_{k}^{m,c} (z):= e^{2 \pi i k c} z^{-c}( \frac{1}{2\pi i}\mbox{Log} ~z - k)^{m-1},
\eeq
are  the restrictions of $f_{k}(s,z,c)$ in \eqn{WL511} to the
domain $\{k \} \times \tilde{\CC}_{01 \infty} \times \{c\}.$ Choosing $k=0$ yields,
for $z$ in the unit interval $(0,1)$,
\begin{eqnarray} \label{N838}
zM_{[Z_1]}^{m,c}(\tZ) &= & -\frac{(2\pi i)^m}{(m-1)!} e^{-2\pi ikc} z^{1-c}(\frac{1}{2 \pi i \,Log \,z})^{m-1}
\nonumber \\
&=& - 2\pi i \left( \frac{1}{(m-1)!} e^{-2\pi ikc} z^{1-c} (\log z)^{m-1}\right),
\end{eqnarray}
since $Log ~z= \log z$ on $(0,1)$. This gives the first row in \eqn{N833}.
~~~$\Box$\\

Now we treat the ``singular'' case $c \in \ZZ_{\le 0}$, and determine  that the monodromy
representation jumps discontinuously between the  ``nonsingular'' and ``singular'' cases.

\begin{theorem}\label{thN83}
{\em ($c$-Deformed Polylogarithm Monodromy-Singular Case)}\\
For each integer $m \ge 1$ and each $c \in \ZZ_{\le 0}$,
the monodromy action on the modified basis
$$
\sB_{m+1,c}^{\ast} := \{ Li_{m,c}^{*}(z), z^{1-c}(\log z)^{m-1}, z^{1-c}(\log z)^{m-2}, \cdots, z^{1-c}  \}
$$
is given by
\beql{N841}
\rho_{m,c} ([Z_0]) := \left(
\begin{array}{cccccc}
1 &  \frac{2 \pi i}{1!}& \frac{(2\pi i)^2}{2!} & \cdots &  \frac{(2 \pi i)^{m - 1}}{(m - 1)!} &
 \frac{(2 \pi i)^{m}}{m!}\\
0 & 1 & \frac {2 \pi i}{1!} & \cdots &
 \frac{(2 \pi i)^{m - 2}}{(m - 2)!} & \frac{(2 \pi i)^{m - 1}}{(m - 1)!} \\
\vdots & \vdots & \vdots & ~ & \vdots & \vdots \\
0 & 0 & 0 & \cdots & 1 & \frac{2\pi i}{1!} \\
0 & 0 & 0 & \cdots & 0 & 1 \\
\end{array}
\right)
\eeq
and
\beql{N842}
\rho_{m,c} ([Z_1]) := \left( \begin{array}{cccccc}
1 & -2\pi i & 0 & \cdots & 0 & 0 \\
0 & 1 & 0 & \cdots & 0 & 0 \\
\vdots & \vdots & \vdots & \vdots & \vdots & \vdots \\
0 & 0 & 0 & \cdots & 1 & 0 \\
0 & 0 & 0 & \cdots & 0 & 1
\end{array}
\right) ~.
\eeq
In these cases the image of $\rho_{m,c}$
 falls in a unipotent subgroup of $GL(m+1, \CC)$.
\end{theorem}

\paragraph{\bf Proof.}
Recall that for $c=-k \in \ZZ_{\le 0}$,
$$
Li_{m,-k}^{*} = \sum_{{n=0}\atop{n \ne k}}^{\infty} \frac{z^{n+1}}{(n-k)^m} +
\frac{1}{m!} z^{k+1}(\log z)^m.
$$
We compute the monodromy by taking a scaling limit approaching the singular point.
Letting  $c= -k + \vep$ for positive $\vep$, we assert that
\beql{N844}
Li_{m,-k}^{\ast} = \lim_{\vep \to 0^{+} }
Li_{m,c} - \frac{1}{\vep^m} \left( \sum_{j=0}^{m-1}\frac{\vep^j}{j!} z^{1-c} (\log z)^j \right).
\eeq
To verify this, we use the expansion
$$
z^{1-c} = z^{k+1} z^{-\vep} = z^{k+1} \sum_{k=0}^{\infty} \frac{(-1)^k}{k!} (\vep \log z)^k.
$$
Now we assert that the  expansion in powers of $\vep$ is
\begin{eqnarray*}
\frac{1}{\vep^m} \left( \sum_{j=0}^{m-1}\frac{\vep^j}{j!} z^{-\vep} (\log z)^j \right)
&=&
\frac{1}{\vep^m}  \sum_{j=0}^{m-1} \left(\sum_{k=0}^{\infty}
\frac{1}{j!} \frac{(-1)^k}{k!}\vep^{j+k}(\log z)^{j+k}  \right)\\
&=& \frac{1}{\vep^m} - \frac{1}{m!}(\log z)^m + O \left(\vep\right).
\end{eqnarray*}
Here the  last step follows using the identities, if $1 \le n \le m-1$,
\beql{N843}
\sum_{{j+k= n}\atop{0 \le j \le m-1}}  \frac{1}{j!} \frac{(-1)^k}{ k!}= \frac{1}{n!} (1-1)^n =0
\eeq
while if $n=m$,
$$
\sum_{{j+k= m}\atop{0 \le j \le m-1}}  \frac{1}{j!} \frac{(-1)^k}{k!}= - \frac{1}{m!}.
$$
This expansion in powers of $\vep$, holding $z$ fixed,  gives \eqn{N844}.

Now for $z=1$ the right side of \eqn{N844} has no monodromy except that coming
from $Li_{m, c}(z)$ hence the limit of the right side yields
$$
M_{[Z_1]} (Li_{m, -k}^{\ast})(z) = - 2\pi i (\frac{1}{(m-1)!} z^{k+1} (\log z)^m.
$$
The formula \eqn{N842} follows since the other basis elements in $\sB^*_{m+1, -k}$
are analytic at $z=1$.

For  $z=0$, we obtain
\begin{eqnarray}\label{N845}
M_{[Z_0]} \left(Li_{m, -k}^{\ast}(z) \right)&=& \lim_{\vep \to 0^{+}} M_{[Z_0]}( Li_{m,-k +\vep}(z))
- \frac{1}{\vep^m} z^{k+1}
\left( \sum_{j=0}^{m-1} \frac{\vep^j}{j!} M_{[Z_0]} (z^{\vep}(\log z)^j)\right)
\nonumber \\
&=& - \frac{1}{\vep^m} z^{k+1} \left( \sum_{j=0}^{m-1} \frac{\vep^j}{j!}
 \left(e^{-\vep(\log z + 2 \pi i)}(\log z+ 2\pi i)^j - e^{-\vep \log z} (\log z)^j\right)\right)
\nonumber \\
&=& - \frac{1}{\vep^m} z^{k+1} \sum_{j=0}^{m-1}
\left( \sum_{k=0}^{\infty} \frac{1}{j!}\frac{(-1)^k}{k!} \vep^{j+k}(\log z + 2\pi i)^{j} \right)
\nonumber \\
&&+ \frac{1}{\vep^m} z^{k+1}\sum_{j=0}^{m-1} \left( \sum_{k=0}^{\infty}
\frac{1}{j!}\frac{(-1)^k}{k!} \vep^{j+k}(\log z)^{j}\right)
\nonumber \\
&=&
\frac{1}{m!}z^{k+1} \left( (\log z + 2 \pi i )^m - (\log z)^m\right) +O (\vep),
\end{eqnarray}
where we again used \eqn{N843} at the last step. Expanding this term gives the first
row in the monodromy matrix \eqn{N841}. The remaining rows follow using
\eqn{N831}.
$~~~\Box$\\

\begin{rem}\label{remark86}
{\em
There is an interpretation of the polylogarithm as  a  variation of Hodge
structure over $\PP^1 (\CC) \smallsetminus \{0,1, \infty\}$,
described in Bloch \cite[p. 278]{Bl91}. This property does cannot extend to irrational parameters
$c \not\in \QQ$ because the monodromy is then not
quasi-unipotent.
}
\end{rem}
\begin{rem}\label{remark87}
{\em
For $c \in \CC \smallsetminus \ZZ_{\le 0}$ the vector space $\sW_{m+1, c}$
is spanned by the  vector space of
functions on
$\tilde{\CC}_{0,1,\infty}$
that  is the direct sum of
all the images of $Li_{m,c}$ under the action of $\hGZ$, i.e.
$$
\sW_{m+1,c} = \bigoplus_{[\tau] \in \hGZ} \CC[Q_{[\tau]}(Li_{m,c})],
$$
in which $Q_{[\tau]}$ are the operators given in  Definition \ref{Nde21}.
Viewed as functions of $z$ on $\tilde{\CC}_{01 \infty}$, using \eqn{N833},
the
vector space spanned by $\{zM_{[Z_0]^{-k} [Z_1][Z_0]^k}^{m,c}(\tZ)~:~ k \in \ZZ\}$
 is the same vector space as that
spanned by $\{ zf_k^{m,c}~:~k \in \ZZ \}$.
By expanding  the products $( \frac{1}{2\pi i}\mbox{Log} ~z - k)^{m-1}$
for each $k$ we see that
this vector space is spanned by $\{ zf_0^{j,c}~:~1 \leq j \leq m\},$ or,
 alternatively, $\{zM_{[Z_1]}^{j,c}(\tZ)~:~1 \leq j \leq m\}.$
 }
\end{rem}

\begin{rem}\label{remark88}
{\em
Comparing the formulas of Theorem \ref{thN83} for $c=1$ with those of
D. Ramakrishnan for the monodromy
of the polylogarithm, in  \cite{Ra82} and \cite[Sect. 4.2 and Sect. 7.6]{Ra89},
we note some discrepancies.   The  formulas for
$[Z_1]$ given in \cite[Sect 4.2]{Ra82}  disagree with ours in a sign,
and those in  \cite[Prop. 7.6.7]{Ra89} we believe
 have a misprint that interchanges $[Z_0]$ and $[Z_1]$.
 (We think our formulas are correct.)
}
\end{rem}

%
%
\section{Further directions}\label{sec9}

There are many directions
for  further investigation; we discuss a few of them.

(1) Can one better understand
 the  nature of the singularities on the  singular strata?
   In particular, the  Riemann zeta function is (formally) obtained as a limit function on a doubly-specialized
 singular stratum.   Part I showed that  there are obstructions to this
 limiting process: for example,
 limiting values approaching the singular stratum $a=1$ in the $(a,c)$-variables exist
 only for $Re(s) >1$.
  It is an   interesting problem  to obtain
 ``renormalized'' limits on singular strata for other ranges of the $s$-variable.
 In part I \cite[Sect. 6]{LL1} the authors  gave a way to do this for the Lerch zeta function for one
 fixed singular stratum, by removing a small number of  divergent terms.
   In this paper  in Section \ref{sec6} we showed  that
 one can extract data approaching  the singular stratrum
 at   $c= 0^{+}$   and negative integer $s$ sufficient to
 reconstruct $p$-adic $L$-functions.

 (2)  The   Lerch zeta function possesses additional discrete symmetries.
One can define an action of  a commuting family of  (two-variable)
``Hecke operators''  in the $(a, c)$ variables on the Lerch zeta function (resp. $(z, c)$-variables for the Lerch transcendent)
for which $\zeta(s, a, c)$ (resp. $\Phi(s, z, c)$)  is a   simultaneous eigenfunction, acting on
various function spaces. Part IV (\cite{LL4}) of this series considers 
such operators on a function space  with real variables.
 These additional  discrete symmetries together with the differential equation \eqref{107b}
  suggest that there should be  an   automorphic interpretation
 of the Lerch zeta function, made in terms of  the related functions $L^{\pm}(s, a, c)$.
 The first author has found such an  interpretation for the real variables form treated in parts I and IV of this series,
showing that symmetrized Lerch functions
with characters are Eisenstein series  on the real Heisenberg group 
quotiented by the integer Heisenberg group (\cite{Lag14a}). It
is an open problem to find  an automorphic interpretation for 
 the complex-analytic version of the Lerch zeta function (reap. Lerch transcendent) treated  in
 part II and this paper.

 (3). One can ask if results of this paper might be interpretable  in the framework of
 an algebraic $D$-module over the  Weyl algebra
 $ \CC[ c, z, \frac{\partial}{\partial c}, \frac{\partial}{\partial z}]$,
     which is associated to a parametric family of linear PDE's $\Delta_{\Phi} - sI$ with
     $$ \Delta_{\Phi} := \frac{1}{2}\Big(\hD_{\Phi}^{+} \hD_{\Phi}^{-} + \hD_{\Phi}^{-}\hD_{\Phi}^{+}) =
      z \frac{\partial}{\partial z}\frac{\partial}{\partial c}  + c \frac{\partial}{\partial c} +\frac{1}{2}I,$$
      with $s$ as eigenvalue parameter.
   Such a reformulation may involve a non-holonomic $D$-module with an infinite set of singularities.

 (4) What are the properties of the extension of polylogarithms 
 under deformation in the $c$-variable? 
As  mentioned in Section \ref{sec13}, one can ask
whether  functional
equations such as the five-term relation for the dilogarithm might
survive in some  fashion under  the $c$-deformation
of the polylogarithm studied in Section \ref{sec8}.
Specifically, the integer points  $c=m \ge 2$ and $s=n \ge  1$
have  maximally unipotent monodromy with
apparent discontinuity in the monodromy matrices.
One can ask whether the functions at these special points 
 satisfy  interesting identities in parallel fashion to the polylogarithms.

The dilogarithm is known to have a single-valued variant,
 the {\em Rogers dilogarithm}, obtained by adding a correction term.
One may wonder if there exists analogous single-valued variant  of
the extended function in the $c$-variable, or at specific integer points $c=m \ge 2$.


 (5) One may investigate  generalizations of the Lerch transcendent in
the direction of an ``elliptic Lerch zeta function'', made in analogy with
work of  Beilinson \cite{BL94} and Levin \cite{Lev97})  on the elliptic
polylogarithm.

(6)   The partial differential operator  $\Delta_{\Phi}$ in \eqn{107b}
in the introduction  can be viewed  as an unbounded operator acting on functions 
    restricted to the domain 
    $$
    T:= \{ (z, c) \in S^1 \times [0,1]\},  \quad \mbox{with} \quad S^1 = \{ |z|=1 \},
   $$ 
 inside the Hilbert space $L^2 \big( T, \frac{dz}{z} \, dc \big)$.
 On this Hilbert space $\Delta_{\Phi}$  can  be shown to be formally skew-adjoint and to have the $xp$-form suggested
 by Berry and Keating (\cite{BK99}, \cite{BK99b})
  as the possible form of a Hilbert-Polya operator encoding the zeta zeros as eigenvalues.
  This operator is obtained from the corresponding operator 
  $ \Delta_L = \frac{1}{2 \pi i} \frac{\partial}{\partial a}\frac{\partial}{\partial c} + c \frac{\partial}{\partial c}+ \frac{1}{2}$ 
  for the Lerch zeta function, which satisfies
  $$
  \Delta_L(s, a, c) = -(s-\frac{1}{2}) \zeta(s, a, c).
  $$
 acting on the Hilbert space $L^2 \big( [0,1]^2, da dc \big)$,
 which is treated in \cite[Sect. 9.2]{Lag15a}.
 
     One may search for natural  skew-adjoint
``boundary conditions''  on the operator $\Delta_{\Phi}$ or $\Delta_L$
which yield operators having spectra on the line $\Re(s-\frac{1}{2}) = 0$.
One such set of boundary conditions will be  presented  in \cite[Sect. 9]{Lag15};
the spectrum of the resulting operator is purely continuous. 
It is an open question whether one can formulate  natural geometric boundary conditions
 on $\Delta_{\Phi}$ that will yield  a  Hilbert-Polya operator for $\zeta(s)$.

\subsection*{Acknowledgments}
 The authors  thank Dinakar Ramakrishnan for conversations
regarding his work on  polylogarithms.
The first author thanks   Peter Scott for discussions and queries on
multidimensional covering manifolds. The authors thank
the reviewers for helpful comments.
This project  was initiated while the first author was at
AT\&T Labs-Research and the second author consulted there;
they thank AT\&T for support. The first author  received
support from the Mathematics Research Center at Stanford University
in 2009-2010. The second author received support from the National Center for
Theoretical Sciences and National Tsing Hua University in Taiwan in 2009-2014.
To these institutions the authors express their gratitude.

\end{document}